\documentclass[12pt,fleqn] {article}
\usepackage{amssymb,amsmath,amsthm}
\usepackage[english]{babel}
\usepackage{graphicx}
\usepackage{a4wide}
\usepackage{enumitem}
\newcommand{\R}{\mathbb{R}}

\newcommand{\bs}{\boldsymbol}

\newcommand{\adb}{\allowdisplaybreaks}

\newtheorem{Theorem}{Theorem}

\newtheorem{Lemma}{Lemma}
\newtheorem{Corollary}{Corollary}

\usepackage{hyperref}
\usepackage{centernot}
\newtheorem{thmx}{Theorem}
\usepackage{natbib}
\setcitestyle{authoryear,open={(},close={)}} %Citation-related commands
\newcommand{\inv}{\frac{1}}

\begin{document}
%\textbf{Editing for ArXiv submission; begun 2022-07-20} %\par
%\hfill \textbf{2019-07-25.1} %\par \textbf{Slight Revn of pp. 1-26, of 2019-07-04, and incorporating Thomas's 2019-07-25 from its p.30 onwards }\par
%\hfill \textbf{Start quoting Fung and Seneta (2017)}
%\begin{center}
%\textbf{Reformulation of Theorem 1}
%\end{center}
%\documentclass[12pt]{article}
\title{Tail asymptotics for the bivariate skew normal in the general case\footnote{Edited version of 2019-07-25.1, begun 2022-07-22}}

\author{Thomas Fung$^{1,}$\footnote{Corresponding author. Email address: thomas.fung@mq.edu.au (Thomas Fung).} and Eugene Seneta$^{2}$ \\ {\small $^{1}$
School of Mathematical and Physical Sciences, Macquarie University, North Ryde, N.S.W. 2109, Australia} \\
{\small $^{2}$
School of Mathematics and Statistics FO7, The University of Sydney, N.S.W. 2006, Australia}}

%\hfill{Editing for ArXiv submission; begun 2022-07-22; of 2019-07-25.1}

\maketitle

\begin{abstract}
The present paper is a sequel to and generalization of \cite{fungTailAsymptoticsBivariate2016} whose main result gives the asymptotic behaviour  as $ u \to 0^{+}$ of $\lambda_L(u) =  P(X_1 \leq F_1^{-1}(u) | X_2
\leq F_2^{-1}(u)),$ when $\textbf{X} \sim SN_2(\bs{\alpha}, R)$ with $\alpha_1 = \alpha_2 = \alpha,$ that is: for the bivariate skew normal distribution in the equi-skew case, where $R$ is the correlation matrix, with off-diagonal entries $\rho,$ and $F_i(x), i=1,2$ are the marginal cdf's of $\textbf{X}$. A  paper of  \cite{berangerModelsExtremalDependence2017} enunciates an upper-tail  version which does not contain the constraint $\alpha_1=\alpha_2= \alpha$ but requires the constraint $0 <\rho <1$ in particular.   The  proof, in their Appendix A.3,  is very condensed.  When  translated to the  lower tail setting of \cite{fungTailAsymptoticsBivariate2016}, we find that when  $\alpha_1=\alpha_2= \alpha$ the exponents of $u$ in the regularly varying function asymptotic expressions do agree, but the slowly varying components, always of asymptotic form $const (-\log u)^{\tau}$, are not asymptotically equivalent.  Our general approach encompasses the case $ -1 <\rho < 0$, and covers all possibilities.
\end{abstract}

\section{Introduction}
The coefficient of lower tail dependence of a random vector $\textbf{Z} =
(Z_1,Z_2)^{\top}$ with marginal inverse distribution functions $F_1^{-1}$ and
$F_2^{-1}$ is defined as
\begin{equation}
\lambda_L = \lim_{u \rightarrow 0^+}\lambda_L(u), \quad \text{where} \quad \lambda_L(u) = P(Z_1 
\leq F_1^{-1}(u) | Z_2
\leq F_2^{-1}(u)).\label{defn:tail dependence}
\end{equation}
%The coefficient of upper tail dependence of a random vector $\textbf{X}$ is defined similarly as
%\begin{equation}
%\lambda_U = \lim_{u \rightarrow 1^-}\lambda_U(u), \quad \text{where} \quad \lambda_U(u) = P(X_1 \geq F_1^{-1}(u) | X_2
%\geq F_2^{-1}(u)). \label{defn:tail upper dependence}
%\end{equation}
$\textbf{Z}$ is said to have asymptotic lower  tail dependence if $\lambda_L$ exists and is positive. If $\lambda_L=0$, then $\textbf{Z}$ is said to be
asymptotically independent in the lower tail. 
%$\textbf{X}$ is said to have asymptotic lower (rest. upper) tail dependence if $\lambda_L$ (resp. $\lambda_U$) exists and is positive. If $\lambda_L=0$ (resp. $\lambda_U=0$), then $\textbf{X}$ is said to be
%asymptotically independent in the lower (resp. upper) tail. 
This quantity provides insight on
the tendency for the distribution to generate joint extreme event since it
measures the strength of dependence (or association) in the lower tails of a
bivariate distribution.  If the marginal distributions of these random variables are
continuous, then from (\ref{defn:tail dependence}), it follows that
$\lambda_L(u)$ can be expressed in terms of  the copula of
$\textbf{Z}$, $C(u_1,u_2)$, as
\begin{equation*}
\lambda_L(u) = \frac{P(Z_1\leq F_1^{-1}(u), Z_2\leq F_2^{-1}(u))}{P(Z_2 \leq F_2^{-1}(u))} =  \frac{C(u,u)}{u}.
\label{defn:tail dependence 3}
\end{equation*}

Thus the asymptotic behaviour  as $ u \to 0^+$ of the copula $C(u,u)$  is tantamount to that of $\lambda_L(u)$ 
through the relation:
\begin{equation*}
C(u,u)= u\lambda_L(u).  \label{copula} 
\end{equation*}

In this note our focus is  the bivariate skew normal distribution. Analysing $\lambda_L(u), $ we shall also obtain explicit forms for $\kappa,\,  \kappa >0  ,$ and the slowly varying function $\mathcal L(u)$ as $u \to 0^+,$ in
\begin{equation}
C(u,u)\sim u^{\kappa} \mathcal L(u)  \label{copula1}
\end{equation}
under completely general conditions. 

The bivariate skew normal distribution was introduced in \cite{azzaliniMultivariateSkewnormalDistribution1996} (which is discussed further in \cite{azzaliniStatisticalApplicationsMultivariate1999}); \cite{azzaliniSkewNormalRelatedFamilies2014}  contains a review. A random vector $\textbf{Z} = (Z_1, Z_2)^{\top}$ is said to have a bivariate skew normal distribution, denoted as $\textbf{Z} \sim SN_2(\bs{\alpha}, R)$, 
if the probability density of $\textbf{Z}$ is  
\begin{equation*}
f(\textbf{z}) = 2\phi_2(\textbf{z},R)\Phi(\bs{\alpha}^{\top}\textbf{z}), \label{pdf:SN_n}
\end{equation*}
where $\phi_2(\cdot,R)$ is density of a bivariate normal distribution with mean $\textbf{0}$ and correlation matrix $R$ and $\Phi(\cdot)$ is the cdf of a univariate standard normal distribution. The correlation matrix $R$ and skew vector $\bs{\alpha}$ are defined as $\left(\begin{smallmatrix} 1 & & \rho \\ \\ \rho & & 1 \end{smallmatrix}\right)$, with $-1<\rho<1$ and $\bs{\alpha} =(\alpha_1,\alpha_2)^{\top} \in \R^2$ respectively. Obviously, the (symmetric) bivariate normal is obtained as special case when $\bs{\alpha}=0$. 
%The results of \citet{lyse2009} and \citet{bort2010} show  that the skew normal distribution is tail independent, that is: $\lambda_L = \lim_{u \rightarrow 0^+}\lambda_L(u) =0.$ 
The marginal density for $Z_i$, $i=1$, $2$, is 
\begin{equation*}
f_{Z_i}(z_i) = 2 \phi(z_i) \Phi(\lambda_i z_i), \quad \text{for $z_i\in \R$}, \label{pdf:SN_1}
\end{equation*}
where $\phi(\cdot)$ and $\Phi(\cdot)$ are the pdf and cdf of the univariate standard normal and 
\begin{equation}
\lambda_i = \frac{\alpha_i+\rho \alpha_{3-i}}{\sqrt{1+\alpha_{3-i}^2(1-\rho^2)}}. \label{lambda}
\end{equation}
 This means that both $Z_1$ and $Z_2$ have a univariate skew normal distribution with skewness parameter $\lambda_i$ i.e. $Z_i \sim SN(\lambda_i)$, $i=1$, $2$.

 In our  paper \cite{fungTailAsymptoticsBivariate2016}, to which the present paper is a sequel and generalisation, following earlier authors cited there for terminology, we call $ \kappa $  the (lower) tail order of the copula (\ref{copula1}). The case  $\kappa = 2$ is that of asymptotic tail independence, the terminology deriving from the subcase when ${\mathcal L}(u) = const.$ The  tail order case $ 1 <\kappa <2 $  is considered as intermediate  tail dependence, as it corresponds to the copula having some level of positive dependence in the tail when $\lambda_L =0$, but the association is not as strong as when $\kappa =1,$ and $ \lambda_L(u) ={\mathcal L}(u) \to \lambda_L > 0, u \to 0^{+}.$ 

In \cite{fungTailAsymptoticsBivariate2016} the following was the  main result.

\begin{thmx}\label{thm:fungTailAsymptoticsBivariate2016}

Let $\textbf{X} \sim SN_2(\bs{\alpha}, R)$ with $\alpha_1 = \alpha_2 = \alpha$. As $u\to 0^+$, 
\begin{enumerate}
\item[(a)] if $\alpha>0$, {\adb
\begin{align}
%\notag & \lambda_L(u) \\
\lambda_L(u) &\sim u^{\beta^2}\frac{(2\pi\lambda)^{\beta^2}}{\sqrt{\pi}\beta(1+\beta^2)^2}[-\log u]^{\beta^2-\inv{2}} \nonumber,
%
%\frac{\alpha^3}{\pi \lambda^4\beta(1+\beta^2)^2} \sqrt{\frac{2}{\pi}} (2\pi \lambda)^{1+\beta^2}
%(\frac{1+\lambda^2}{2})^{\frac{3}{2}} \left[ -\log  u\right]^{\beta^2-\inv{2}}
%\\\frac{(2\pi\lambda)^{\beta^2}}{\sqrt{\pi}\beta(1+\beta^2)^2}[-\log u]^{\beta^2-\inv{2}}\times u^{\beta^2}
%\notag &\sim& u^{\frac{(1-\rho)(1+2\theta^2(1+\rho))}{1+\rho}} \frac{1}{\pi \beta(1+\beta^2)}\left(\frac{1+2\theta^2(1+\rho)}{1+\theta^2(1-\rho^2)}\right) \frac{(1+\rho)}{\alpha(2+4\theta^2(1+\rho))}\\
%&& \quad\times  \sqrt{\frac{2}{\pi}} (2\pi \lambda)^{\frac{2(1+\theta^2(1-\rho^2))}{1+\rho}} 
%(\frac{1+\lambda^2}{2})^{\frac{3}{2}} \left[ -\log u\right]^{\frac{2(1+\theta^2(1-\rho^2))}{1+\rho}-\frac{3}{2}},
\end{align}}
with $\lambda = \frac{\alpha(1+\rho)}{\sqrt{1+\alpha^2(1-\rho^2)}}$ and $\beta = \sqrt{\frac{(1-\rho)(1+2\alpha^2(1+\rho))}{1+\rho}}$;
% $\delta = \frac{\alpha(1+\rho)}{\sqrt{1+2\alpha^2(1+\rho)}}= \frac{\lambda}{\sqrt{1+\lambda^2}}$
\item[(b)] if $\alpha<0$,
\begin{equation}
\lambda_L(u) \sim u^{\frac{1-\rho}{1+\rho}} \times \frac{1+\rho}{2} \sqrt{\frac{1+\rho}{1-\rho}} (-\pi \log u)^{-\frac{\rho}{1+\rho}}.\nonumber
\end{equation}
\end{enumerate}
\end{thmx}

A  paper of \cite{berangerModelsExtremalDependence2017} enunciates  an upper-tail  version which does not contain the constraint $\alpha_1=\alpha_2= \alpha$ but requires the constraint $0 <\rho <1$ in particular.   The  proof, in their Appendix A.3,  is very condensed.  When  translated to our present  lower tail setting, we find that when  $\alpha_1=\alpha_2= \alpha$ the exponents of $u$ do agree with those of Theorem \ref{thm:fungTailAsymptoticsBivariate2016} above, but the slowly varying functions, always of asymptotic form $const (-\log)^{\tau_2}$, where $const$, $\tau_2$ depend on the signs of $\lambda_1, \lambda_2$ in particular,  are not asymptotically equivalent to those of Theorem \ref{thm:fungTailAsymptoticsBivariate2016}.  Our general approach encompasses the case $ -1 <\rho < 0$, and covers all possibilities.  

\section{Evaluation}
In the present  paper we proceed by noting that
\begin{equation}
C(u,u)= \int_0^u \frac{dC(x,x)}{dx}dx \label{defn:tail dependence13}
 \end{equation}
so that if ${\frac{dC(x,x)}{dx}}= x^{\theta}L(x),
\theta > 0 $ where $L(x)$ is a slowly varying function as $ x \to
0^+$, then by (applying with suitable transformation to regular
variation at $0$) a result of de Haan (see \cite{senetaRegularlyVaryingFunctions1976}, p. 87), we
obtain
\begin{equation}
\frac{C(u,u)}{u} = \inv{u}\int_0^u {\frac{dC(x,x)}{dx}}dx \sim  {\frac{u^{\theta}L(u)}{\theta + 1}}, \, \, u \to 0^+ .\label{defn:tail dependence133}\end{equation}
Therefore, it is sufficient for us to find a value of $\theta > 0$
which satisfies
\begin{equation*}
\frac{d C(u,u)}{du} = P(Z_1\leq F_1^{-1}(u)|Z_2 = F_2^{-1}(u)) + P(Z_2 \leq F_2^{-1}(u) | Z_1 = F_1^{-1}(u)) = u^{\theta} L(u)
\label{result}
\end{equation*}
for some slowly varying function $L(u)$, as $u\rightarrow 0^+$, so
that (\ref{defn:tail dependence133}) holds. Then (\ref{copula1}) holds with $ \kappa = \theta +1, \, \mathcal{L}(u) = L(u)/(\theta +1 ).$

We begin by considering, without loss of generality, 
\begin{eqnarray}
\notag && P(Z_1\leq F_1^{-1}(u) | Z_2 = F_2^{-1}(u)) \\
&=& \int^{F_{1}^{-1}(u)}_{-\infty} \inv{\sqrt{2\pi(1-\rho^2)}} e^{-\inv{2(1-\rho^2)} (z_1-\rho F_2^{-1}(u))^2} \frac{\Phi(\alpha_1 z_1+\alpha_2F_2^{-1}(u))}{\Phi(\lambda_2 F_2^{-1}(u))} \,dz_1, \label{original expression}
\end{eqnarray}
for $\bs{Z} = (Z_1, Z_2)^{\top} \sim SN_2(\bs{\alpha}, R).$ See\cite{azzaliniSkewNormalRelatedFamilies2014}'s (5.65) for the conditional 
quantity.

Notice that if $\alpha_1 = 0$, then $\lambda_2 = \frac{\alpha_2+\rho\alpha_1}{\sqrt{1+\alpha_1^2(1-\rho^2)}}=\alpha_2$ and (\ref{original expression}) becomes 
\begin{eqnarray*}
&& P(Z_1\leq F_1^{-1}(u) | Z_2 = F_2^{-1}(u)) \\
&=& \int^{F_1^{-1}(u)}_{-\infty} \inv{\sqrt{2\pi(1-\rho^2)}} e^{-\inv{2(1-\rho^2)}(z_1-\rho F_2^{-1}(u))^2}\frac{\Phi(\alpha_2 F_2^{-1}(u))}{\Phi(\lambda_2F_2^{-1}(u))}\,dz_1 \\
&=& \int^{F_1^{-1}(u)}_{-\infty} \inv{\sqrt{2\pi(1-\rho^2)}} e^{-\inv{2(1-\rho^2)}(z_1-\rho F_2^{-1}(u))^2}\,dz_1\\
&=& \Phi\left(\frac{F_1^{-1}(u)-\rho F_2^{-1}(u)}{\sqrt{1-\rho^2}}\right).
\end{eqnarray*}
So $\alpha_1=0$ is not a particularly interesting case and we shall exclude this from the subsequent theoretical development. 

%We will start with the assumption of 
%\begin{equation*}
%\alpha_1>0, \quad \alpha_2>0, \quad \lambda_1>0, \quad \lambda_2>0 \quad \& \quad \alpha_1F_1^{-1}(u) +\alpha_2F_2^{-1}(u) \to -\infty, \quad \text{as $u\to 0^+$}.
%\end{equation*}

It is obvious from (\ref{original expression}) that there are two different change of variables we can implement here. The first one is $x = \frac{z_1 - \rho F_2^{-1}(u)}{\sqrt{1-\rho^2}}$ and (\ref{original expression}) becomes {\adb
\begin{align}
\notag & \int^{\frac{F_1^{-1}(u)-\rho F_2^{-1}(u)}{\sqrt{1-\rho^2}}}_{-\infty} \inv{\sqrt{2\pi}}e^{-\frac{x^2}{2}} \frac{\Phi\left(\alpha_1\sqrt{1-\rho^2}x + (\alpha_2+\rho\alpha_1)F_2^{-1}(u)\right)}{\Phi(\lambda_2F_2^{-1}(u))}\,dx\\
\notag =& 
\int^{\frac{F_{1}^{-1}(u)-\rho F_2^{-1}(u)}{\sqrt{1-\rho^2}}}_{-\infty} \frac{\Phi(\alpha_1\sqrt{1-\rho^2} x+(\alpha_2+\rho\alpha_1)F_2^{-1}(u))}{\Phi(\lambda_2 F_2^{-1}(u))} \,d\Phi(x)\\
\notag =& \left[ \Phi(x) \times \frac{\Phi(\alpha_1\sqrt{1-\rho^2} x+(\alpha_2+\rho\alpha_1)F_2^{-1}(u))}{\Phi(\lambda_2 F_2^{-1}(u))}\right]^{\frac{F_{1}^{-1}(u)-\rho F_2^{-1}(u)}{\sqrt{1-\rho^2}}}_{-\infty} \\
& \quad - \int^{\frac{F_{1}^{-1}(u)-\rho F_2^{-1}(u)}{\sqrt{1-\rho^2}}}_{-\infty} \Phi(x) \frac{\inv{\sqrt{2\pi}} e^{-\inv{2}(\alpha_1\sqrt{1-\rho^2} x+(\alpha_2+\rho\alpha_1)F_2^{-1}(u))^2}}{\Phi(\lambda_2 F_2^{-1}(u))} \times \alpha_1\sqrt{1-\rho^2}\,dx.
\label{Route 1}
\end{align}}
Alternatively, we can  apply a change of variable of $x = \alpha_1z_1 + \alpha_2F_2^{-1}(u)$ when $\alpha_1\ne0$ to (\ref{original expression}) to get 
\begin{equation}
\int^{\alpha_1F_{1}^{-1}(u)+ \alpha_2F_2^{-1}(u)}_{-\infty} \inv{\sqrt{2\pi(1-\rho^2)}} e^{-\inv{2(1-\rho^2)} (\frac{x-\alpha_2F_2^{-1}(u)}{\alpha_1}- \rho F_2^{-1}(u))^2} \frac{\Phi(x)}{\Phi(\lambda_2 F_2^{-1}(u))}\inv{\alpha_1}\,dx,  \label{Route 2}
\end{equation}
if $\alpha_1>0$ or 
\begin{equation}
\int^{\infty}_{\alpha_1F_{1}^{-1}(u)+\alpha_2F_2^{-1}(u)}\inv{\sqrt{2\pi(1-\rho^2)}} e^{-\inv{2(1-\rho^2)} (\frac{x-\alpha_2F_2^{-1}(u)}{\alpha_1}- \rho F_2^{-1}(u))^2} \frac{\Phi(x)}{\Phi(\lambda_2 F_2^{-1}(u))}\inv{|\alpha_1|}\,dx, \label{Route 3}
\end{equation}
if $\alpha_1<0$.

%if $\alpha_1>0$ or 
%\begin{equation}
%\int^{\infty}_{\alpha_1F_{1}^{-1}(u)+\alpha_2F_2^{-1}(u)}\inv{\sqrt{2\pi(1-\rho^2)}} e^{-\inv{2(1-\rho^2)} (\frac{x-\alpha_2F_2^{-1}(u)}{\alpha_1}- \rho F_2^{-1}(u))^2} \frac{\Phi(x)}{\Phi(\lambda_2 F_2^{-1}(u))}\inv{|\alpha_1|}\,dx, \label{Route 3}
%\end{equation}
%if $\alpha_1<0$.
%
%The above integrals would be difficult to handle if the upper limit tends to $\infty$ as any expansion-based techniques work only when the integration is over the region when $x$ is close to $-\infty$. Our subsequent proof would show exactly that when the integrals are well behave. 
%
%%generally fails apart when $F_{1}^{-1}(u)-\rho F_2^{-1}(u) \to \infty$ in (\ref{Route 1}) or $F_{1}^{-1}(u)+\frac{\alpha_2}{\alpha_1}F_2^{-1}(u)\to \infty$ in (\ref{Route 2}). This is because we can no longer be sure that the integration is over the so most of  would not work anymore. 
%
%It is therefore a necessity to first show that it is impossible for both $F_{1}^{-1}(u)-\rho F_2^{-1}(u)$ and $F_{1}^{-1}(u)+\frac{\alpha_2}{\alpha_1}F_2^{-1}(u)$ tending to $\infty$ in the same setting. In other words, for any combinations of $\alpha_1$, $\alpha_2$, $\lambda_1$, $\lambda_2$ and $\rho$, at least one of $F_{1}^{-1}(u)-\rho F_2^{-1}(u)$ and $F_{1}^{-1}(u)+\frac{\alpha_2}{\alpha_1}F_2^{-1}(u)$ is strictly negative and is tending to $-\infty$ as $u\to 0^+$ so that we can funnel (\ref{original expression}) into a well-behaving integral in (\ref{Route 1}) or (\ref{Route 2}). 

We can see that the behaviour of these integrals depends on the asymptotic behaviour as $ u \to 0+$ of the quantities $ F_i(u), i=1,2$ and their expression through the quantities: 
\begin{equation}
A_i(u) = F_i^{-1}(u)-\rho F_{3-i}^{-1}(u) \quad \text{and}\quad  B(u) = \alpha_1F_1^{-1}(u)+\alpha_{2} F_{2}^{-1}(u). \label{Definition of A_1(u) and B(u)}
 \end{equation}
Analysis  of the asymptotic order of the remainder of a partial asymptotic expansion of the quantile function $F_i^{-1}(u)$ as $u\to 0^+$, by using regularly varying functions, was reported in \cite{fungQuantileFunctionExpansion2018} and we will simply initially  restate the results related to the skew normal quantile functions, and then state the consequences for the ratio of $F_i(u)$, $i = 1$, $2$.

\begin{Lemma}\label{Lemma: detailed F inverse}
Let $Z_i\sim SN(\lambda_i)$, with $i=1$, $2$, then {\adb
\begin{equation*}
F_i^{-1}(u) = K_{i,1}(-\sqrt{-2\log u})\left\{ 1+ \frac{K_{i,2}\log(-\log u)}{\log u} + \frac{K_{i,3}}{\log u} +O\left(\left(\frac{\log(-\log u)}{\log u}\right)\right)\right\}
\end{equation*}
as $u\to 0^+$} where {\adb
\begin{align*}
K_{i,1} &= \sqrt{1/(1+\lambda_i^2)},\, K_{i,2} = 1/2,\, K_{i,3} = \log(2\pi\lambda_i)/2,\, \text{when $\lambda_i>0$;}\\
K_{i,1} &= 1,\, K_{i,2} = 1/4,\, K_{i,3} = \log \pi/4,\, \text{when $\lambda_i<0$};\\
K_{i,1} &= 1,\, K_{i,2} = 1/4,\, K_{i,3} = \log 4\pi/4,\, \text{when $\lambda_i=0$}. 
\end{align*}
%{\adb
Hence:
\begin{align}
& \frac{F_i^{-1}(u)}{F_{3-i}^{-1}(u)} = \gamma_i\left\{1 + \frac{C_{i,1}\log|\log u|}{\log u} + \frac{C_{i,2}}{\log u} + O\left(\left(\frac{\log|\log u|}{\log u}\right)^2\right)\right\} \label{detail F ratio}
\end{align}
for $i=1,2,$ where 
\begin{align*}
\gamma_i & = \sqrt{\frac{1+\lambda_{3-i}^2}{1+\lambda_i^2}},\, C_{i,1} = 0,\, C_{i,2} = \frac{\log(\lambda_i/\lambda_{i-3})}{2}, \text{ when $\lambda_i>0$, $\lambda_{3-i}>0$;}\\
\gamma_i& = \sqrt{1+\lambda_{3-i}^2},\, C_{i,1} = -\inv{4},\, C_{i,2} = -\frac{\log(2\lambda_{3-i}\sqrt{\pi})}{2}, \text{ when $\lambda_i<0$, $\lambda_{3-i}>0$;}\\
\gamma_i&= \inv{\sqrt{1+\lambda_i^2}},\, C_{i,1} = \inv{4},\, C_{i,2} = \frac{\log(2\lambda_i\sqrt{\pi})}{2}, \text{ when $\lambda_i>0$, $\lambda_{3-i}<0$;}\\
\gamma_i&= \sqrt{1+\lambda^2_{3-i}},\, C_{i,1} = -\inv{4},\, C_{i,2} = - \frac{\log(\lambda_{3-i}\sqrt{\pi})}{2}, \text{ when $\lambda_i=0$, $\lambda_{3-i}>0$;}\\
\gamma_i&=\inv{\sqrt{1+\lambda_i^2}},\, C_{i,1} = \inv{4},\, C_{i,2} = \frac{\log(\lambda_i\sqrt{\pi})}{2},  \text{ when $\lambda_i>0$, $\lambda_{3-i}=0$;}\\
\gamma_i &= 1, \, C_{i,1} = C_{i,2} = 0,  \text{ when $\lambda_i<0$, $\lambda_{3-i}<0$};\\
\gamma_i &= 1,\, C_{i,1} = 0, C_{i,2} = - \frac{\log 2}{2}, \text{ when $\lambda_i<0$, $\lambda_{3-i}=0$}\\
\gamma_i &= 1,\, C_{i,1}=0, C_{i,2} = \frac{\log 2}{2}, \text{ when $\lambda_i=0$, $\lambda_{3-i}<0$.} 
\end{align*}} 
\end{Lemma} 
It turns out that  the two expressions $A_i(u)$ and $B(u)$  are related in a special way and the relationship is summarised in the following lemma which we need for the moment  only in the case $i=1,$ and the subsequent corollary. Notice from (\ref{Definition of A_1(u) and B(u)}), (\ref{detail F ratio}) that
\begin{equation} \label{limits}
\frac{A_1(u)}{F_2^{-1}(u)} \to \gamma_1 -\rho, \quad \text{and}\quad   \frac{B(u)}{F_2^{-1}(u)} \to \alpha_1\gamma_1+\alpha_2, \quad  \text{as $u \to 0^{+}$}. \end{equation}

\begin{Lemma}\label{Lemma: at least one condition is positive}
Suppose that $\alpha_1\ne 0$. Then at least one of  the limits $\gamma_1 - \rho$ and $\alpha_1 \gamma_1 + \alpha_2$ in (\ref{limits}) is strictly positive for any combination of $\alpha_1$, $\alpha_2$, $\lambda_1$, $\lambda_2$ and $\rho$. Whenever $\gamma_1 - \rho \leq 0$,  then $\alpha_1 >0, \,  \alpha_1 \gamma_1 + \alpha_2 >0.$
\end{Lemma}

The proof of Lemma \ref{Lemma: at least one condition is positive} can be found in the Appendix \ref{Appendix: Proof of Lemma 2}. Overall, Lemma \ref{Lemma: at least one condition is positive} shows that at least one of $\frac{A_1(u)}{F_2^{-1}(u)} = \frac{F_{1}^{-1}(u)}{F_2^{-1}(u)}-\rho$ and $\frac{B(u)}{F_2^{-1}(u)} = \alpha_1\frac{F_{1}^{-1}(u)}{F_2^{-1}(u)}+\alpha_2$ converges to a strictly positive constant for any combination of $\alpha_1$, $\alpha_2$, $\lambda_1$, $\lambda_2$ and $\rho$. 
In fact, since in the cases in (\ref{F ratio}) other than $ \lambda_1 >0, \lambda_2 >0$ we have shown that $\gamma_1 - \rho >0$, the statement:  whenever $\gamma_1-\rho\leq 0$, we have $\alpha_1\gamma_1+\alpha_2>0$ and $\alpha_1>0$ holds for all cases of  (\ref{F ratio}) .

The above lemmas imply the following corollary. 

\begin{Corollary}
At least one of $A_1(u)$ and $B(u)$ (assuming $\alpha_1 \ne 0$) tends to $-\infty$ as $u\to 0^+$ for any combination of $\alpha_1$, $\alpha_2$, $\lambda_1$, $\lambda_2$ and $\rho$.  Furthermore,  $A_1(u)\to -\infty$, if and only if  $\gamma_1-\rho>0$. If $\alpha_1\ne0$ then  $B(u) \to -\infty$, if and only if  $\alpha_1\gamma_1+\alpha_2>0$. Then, respectively, $A_1(u) \sim (\gamma_1-\rho)F_2^{-1}(u)$, $B(u) \sim (\alpha_1\gamma_1+\alpha_2)F_2^{-1}(u)$. 
\label{Corollary: at least one}
\end{Corollary}
\begin{proof}
It follows directly from Lemma \ref{Lemma: at least one condition is positive} that at least one of the limits 
$\frac{A_1(u)}{F_2^{-1}(u)} \to \gamma_1-\rho$ and $\frac{B(u)}{F_2^{-1}(u)} \to \alpha_1\gamma_1+\alpha_2$ is a positive constant, so  at least one of $A_1(u)$ and $B(u) \to -\infty$, as $u\to 0^+$.
% Furthermore, by using  (\ref{detail F ratio}) in Lemma \ref{Lemma: detailed F inverse}, we note that 
%\begin{equation*}
%F_2^{-1}(u) \left(\frac{F_1^{-1}(u)}{F_2^{-1}(u)}-\gamma_1\right) = O\left(\frac{\log(-\log u)}{\sqrt{-\log u}}\right).
%\end{equation*}
Now
\begin{align*}
A_1(u) & = F_2^{-1}(u)\left(\frac{F_1^{-1}(u)}{F_2^{-1}(u)}-\rho\right) \\
&= F_2^{-1}(u)\left(\frac{F_1^{-1}(u)}{F_2^{-1}(u)}-\gamma_1+\gamma_1-\rho\right)\\
&= O\left(\frac{\log|\log u|}{\sqrt{-\log u}}\right)+(\gamma_1-\rho)F_2^{-1}(u),
\end{align*}
from  (\ref{detail F ratio}). If $\gamma_1-\rho \leq 0$, this contradicts  $A_1(u) \to -\infty$. The other assertion follows similarly.
\end{proof}

The above Corollary \ref{Corollary: at least one} provides us with some structure on integrals (\ref{Route 1})--(\ref{Route 2}) but not for (\ref{Route 3}).  Fortunately, we shall never need to use (\ref{Route 3}) since  Lemma \ref{Lemma: at least one condition is positive}   indicates that  $A_1(u)\to -\infty$ and therefore (\ref{Route 1}) is applicable except in the case of $\lambda_1$, $\lambda_2>0$ in which setting  whenever $\gamma_1 - \rho \leq 0$ we have that  $\alpha_1 >0$, $B(u)\to -\infty$. Thus  when  we can not guarantee $A_1(u) \to -\infty$ we shall use (\ref{Route 2}) instead. As a result, dealing with integrals in (\ref{Route 1})--(\ref{Route 2}) inevitably forms the main part of our results. The following theorem summarises how to deal with such integrals in unified form.

\begin{Theorem}\label{Theorem: evaluate general integral}
Let $a(u)$ and $b(u)$ be two functions in $u$ such that $0> a(u) \to -\infty$ and $b(u)/a(u) \to k$, a constant, as $u\to 0^+$; $c$ be a fixed constant. Write $v(u)=\frac{a(u)(1+c^2) +cb(u)}{\sqrt{1+c^2}}.$  Then {\adb
\begin{eqnarray*}
&& \int^{a(u)}_{-\infty}\Phi(x) \inv{\sqrt{2\pi}} e^{-\inv{2}(cx+b(u))^2}\,dx \\
&\sim& 
\begin{cases}
\frac{e^{-\frac{b^2(u)}{2(1+c^2)}-\inv{2}v^2(u)}}{2\pi \sqrt{1+c^2} |a(u)| |v(u)|}, & \text{if} \, \,  v(u) \to - \infty;\\
\frac{\sqrt{1+c^2}e^{-\frac{b^2(u)}{2(1+c^2)}}}{2\sqrt{2\pi}|a(u)||ck|}, & \text{if $v(u) \to 0$};\\
\frac{\sqrt{1+c^2}e^{-\frac{b^2(u)}{2(1+c^2)}}}{\sqrt{2\pi}|a(u)||ck|}, & \text{if $v(u) \to \infty$},
\end{cases}
\end{eqnarray*}}
as $u\to 0^+$.
%where $C$ is a constant and $\lim_{u\to 0^+} e^{1+\frac{cb(u)}{(1+c^2)a(u)}}\leq C\leq 1$.

%Let $a(u)$, $\beta(u)$, and $w(u)$ be three functions in $u$ such that $\beta(u) >0$ and $a(u)$, $w(u) \downarrow 0$ as $u\to 0^+$. Then 
%\begin{equation}
% \int^{\infty}_{0} e^{-a(u)y^2 - \beta(u)y}\inv{1+ w(u)y}\,dy \sim \inv{\beta(u)}, \label{trick integral}
%\end{equation}
%as $u\to 0^+$.
\end{Theorem} 

%Notice that$v(u \to 0$ implies  $0 =(1+c^2 + ck),$ so that  $-ck = (1+c^2) = |ck|, $ so  an alternate asymptotic representation is possible in this case.  The expression in the statement of Theorem 1 is made to resemble that for when  $v(u) \to \infty.$

The detailed long  proof of the theorem  is deferred to Appendix \ref{Appendix: Proof of Theorem 1}, so as not to disrupt the flow of the mainstream development.

%The above theorem, immediately implies the following Corollaries 2 \& 3 which 
%involves a constant $\beta_{1.2} = \lim_{u\to 0^+} \beta_{1.2}(u)$ where 
%and an identity of: and we shall need these two quantities repeatedly for the rest of this note. 

\section{ Preliminary Lemmas.}
 We first  restate for convenience the definitions of (\ref{lambda}) and (\ref{Definition of A_1(u) and B(u)}) that for $i=1,2$:
\begin{align*}
\lambda_i &= \frac{\alpha_i+\rho \alpha_{3-i}}{\sqrt{1+\alpha_{3-i}^2(1-\rho^2)}} 
\tag*{(\ref{lambda})}\\
A_i(u)&= F_i^{-1}(u)-\rho F_{3-i}^{-1}(u) \quad \text{and}\quad  B(u) = \alpha_1F_1^{-1}(u)+\alpha_{2} F_{2}^{-1}(u), 
\tag*{(\ref{Definition of A_1(u) and B(u)})}
\end{align*} 
and introduce the definition:
\begin{align}
\beta_{i}(u) &= \left[ \frac{A_i(u)}{F_{3-i}^{-1}(u)(1-\rho^2)} + \frac{\alpha_i}{F_{3-i}^{-1}(u)}B(u)\right] \label{beta12u definition 3} \\
\notag & = \left[\frac{F_i^{-1}(u)-\rho F_{3-i}^{-1}(u)}{F_{3-i}^{-1}(u)(1-\rho^2)} + \frac{\alpha_i}{F_{3-i}^{-1}(u)}\left(\alpha_1F_1^{-1}(u)+\alpha_2F_2^{-1}(u)\right)\right] 
% \label{beta12u definition 1} 
\\
&= \frac{F_i^{-1}(u)}{F_{3-i}^{-1}(u)}\left(\alpha_i^2+\inv{1-\rho^2}\right)+\alpha_1\alpha_2-\frac{\rho}{1-\rho^2}. 
\label{beta12u definition 2} 
\end{align} so that 
\begin{align}
\notag \beta_{i} &= \lim_{u\to 0^+} \beta_{i}(u)  = \gamma_i\left(\alpha_i^2+\inv{1-\rho^2}\right)+\alpha_1\alpha_2-\frac{\rho}{1-\rho^2}\\
&=\frac{\gamma_i - \rho}{1 - \rho^2} + \alpha_i (\alpha_i \gamma_i + \alpha_{3-i}) \label{beta12 definition}
\end{align}
where $\gamma_i$ is defined by (\ref{detail F ratio}).

For the  two applications of Theorem \ref{Theorem: evaluate general integral}, expressed as Corollaries \ref{Corollary: specific integral with rho in the condition, second term} and \ref{Corollary: specific integral without rho in the condition} in the next section below, we shall  need the following Lemmas  \ref{Lemma: exponents equivalency} to \ref{Lemma: beta1>0} and their proofs can be found in Appendix \ref{Appendix: Proof of Lemmas}.
%We shall then need the following result. 
\begin{Lemma}\label{Lemma: exponents equivalency}
For $i=1,2$,
\begin{align}
& \lambda_{3-i}^2(F_{3-i}^{-1}(u))^2 +\left(\frac{1-\rho^2}{1+\alpha_i^2(1-\rho^2)}\right)\left(\beta_{i}(u)F_{3-i}^{-1}(u)\right)^2  =  \frac{A^2_i(u)}{1-\rho^2}+ B^2(u),\label{Lemma: eqn exponents equivalency}
%\notag =  & -\frac{(\alpha_2+\rho\alpha_1)^2(F_2^{-1}(u))^2}{2(1+\alpha_1^2(1-\rho^2))} -\inv{2}\left(\frac{1-\rho^2}{1+\alpha_1^2(1-\rho^2)}\right)\left(\frac{F_1^{-1}(u)-\rho F_2^{-1}(u)}{1-\rho^2} +\alpha_1(\alpha_1F_1^{-1}(u)+\alpha_2 F_2^{-1}(u))\right)^2\\ 
\end{align}

\end{Lemma}

We shall need to use various other groupings of the components of (\ref{Lemma: eqn exponents equivalency}) as exponents in the sequel, for example: {\adb
\begin{align} 
  -&\frac{A_1^2(u)}{1-\rho^2}+ \lambda_2^2(F_2^{-1}(u))^2  \nonumber\\
= & B^2(u) - \left(\frac{1-\rho^2}{1+\alpha_1^2(1-\rho^2)}\right)\left(\beta_{1}(u)F_2^{-1}(u)\right)^2  \nonumber\\
= &\left(B(u)+\frac{\sqrt{1-\rho^2}}{\sqrt{1+\alpha_1^2(1-\rho^2)}}\beta_{1}(u)F_2^{-1}(u)\right)\times \left(B(u) - \frac{\sqrt{1-\rho^2}}{\sqrt{1+\alpha_1^2(1-\rho^2)}}\beta_{1}(u)F_2^{-1}(u)\right) \nonumber \\
= &\left(B(u)+\frac{\sqrt{1-\rho^2}}{\sqrt{1+\alpha_1^2(1-\rho^2)}}\beta_{1}(u)F_2^{-1}(u)\right)\left[B(u) - \frac{\sqrt{1-\rho^2}}{\sqrt{1+\alpha_1^2(1-\rho^2)}}\left(\frac{A_1(u)}{1-\rho^2}+\alpha_1B(u)\right)\right] \nonumber\\
= &\left(B(u)+\frac{\sqrt{1-\rho^2}}{\sqrt{1+\alpha_1^2(1-\rho^2)}}\beta_{1}(u)F_2^{-1}(u)\right) \nonumber\\
& \quad \times \Biggl[\frac{\left(\sqrt{1+\alpha_1^2(1-\rho^2)}-\alpha_1\sqrt{1-\rho^2}\right)B(u)}{\sqrt{1+\alpha_1^2(1-\rho^2)}} - \frac{\sqrt{1-\rho^2}}{\sqrt{1+\alpha_1^2(1-\rho^2)}}\left(\frac{A_1(u)}{1-\rho^2}\right)\Biggr]\label{special}.   
\end{align} }
 using (\ref{beta12u definition 3}). 

\begin{Lemma} \label{Lemma: beta limit}
\begin{equation*}
\left(\beta_{i}(u) -\beta_{i}\right)\times F_{3-i}^{-1}(u) = O\left(\frac{\log(-\log u)}{\sqrt{-\log u}}\right) \to 0, \quad \text{as $u\to 0^+$.}
\end{equation*}
\end{Lemma}

In particular, when $\lim_{u\to 0^+} \beta_{i}(u) = \beta_{i} = 0$
we have 
\begin{equation}
\beta_{i}(u)F_{3-i}^{-1}(u) = O\left(\frac{\log(-\log u)}{\sqrt{-\log u}}\right) \to 0 \label{eqn beta12uF2inv goes to 0}
\end{equation}
as $u\to 0^+$. 
\begin{Lemma}\label{Lemma: beta1>0}
If $\lambda_2\geq0$ or $B(u) \to -\infty$, or $B(u) \to 0$ then $\beta_{1}>0$. Similarly, we have $\beta_2>0$ when $\lambda_1\geq 0$ or $B(u) \to -\infty$ or $B(u) \to 0$.
\end{Lemma}

Although we shall need Lemma  \ref{Lemma: beta1>0} only when  $B(u) \to -\infty$ and when $B(u) \to 0$, the proof, by its nature, would not be simplified from that for the present more general statement. Next, the present statement is more general, since it is possible that $\lambda_2 \geq 0$ but neither of $B(u) \to -\infty$, $B(u) \to 0$ holds.  For example, when $\lambda_1 <0$, $\lambda_2 >0$, $-1 < \rho < 0$, the combination  $\alpha_1 <0$, $\alpha_2 <0$ can occur.  Finally, some conditions, as in Lemma \ref{Lemma: beta1>0}, are needed  to ensure that $\beta_1 >0$, for in other circumstances  $\beta_1 =  0$, $\beta_1 < 0$ sometimes do occur, and we use a classification according to the sign of $\beta_1$ in the sequel.

\section{Two Applications of Theorem \ref{Theorem: evaluate general integral} and the Main Result}

Now we can finally discuss  applications of Theorem \ref{Theorem: evaluate general integral}.

In its Corollaries \ref{Corollary: specific integral with rho in the condition, second term} and \ref{Corollary: specific integral without rho in the condition} we shall see that even though $a(u)$, $b(u)$ and $c$ are going to be defined differently they share the  expressions:  
\begin{equation}
- \frac{b^2(u)}{2(1+c^2)} = -\inv{2}\lambda_2^2(F_2^{-1}(u))^2, \label{eqn relating b(u) and c^2 with lambda2 and F2inv}
\end{equation}
and 
\begin{equation}
\frac{a(u)(1+c^2)+cb(u)}{\sqrt{1+c^2}} = v(u)= \frac{\sqrt{1-\rho^2}}{\sqrt{1+\alpha_1^2(1-\rho^2)}}F_2^{-1}(u)\beta_{1}(u) \label{eqn relating a(u), b(u) with beta12u}
\end{equation}
where  $\lambda_2 = \frac{\alpha_2+\rho\alpha_1}{\sqrt{1+\alpha_1^2(1-\rho^2)}}$. 
From (\ref{eqn relating a(u), b(u) with beta12u}) and (\ref{eqn beta12uF2inv goes to 0}), this means that 
\begin{equation}
\beta_{1} = \lim_{u\to 0^+}\beta_{1}(u) >0,\, =0,\,<0 \Rightarrow   v(u) \to - \infty,\, v(u) \to 0,\, v(u) \to \infty  \label{eqn stating a(u), b(u) combo has 3 cases}
\end{equation}
 respectively. 
Further,  by using (\ref{eqn relating b(u) and c^2 with lambda2 and F2inv}), (\ref{eqn relating a(u), b(u) with beta12u}) and Lemma \ref{Lemma: exponents equivalency} we  get 
\begin{align} \label{special1}
 -\frac{b^2(u)}{2(1+c^2)}-\inv{2}v^2(u)
%=& -\inv{2}\frac{(\alpha_2+\rho\alpha_1)^2(F_2^{-1}(u))^2}{\alpha_1^2(1-\rho^2)} \times \frac{\alpha_1^2(1-\rho^2)}{1+\alpha_1^2(1-\rho^2)} - \inv{2}\left(\frac{1-\rho^2}{1+\alpha_1^2(1-\rho^2)}\right)\left(\beta_{1.2}(u)F_2^{-1}(u)\right)^2 \\
%= & -\inv{2}\left(\lambda_2F_2^{-1}(u)\right)^2 - \inv{2}\left(\frac{1-\rho^2}{1+\alpha_1^2(1-\rho^2)}\right)\left(\beta_{1.2}(u)F_2^{-1}(u)\right)^2 \\
= &  -\frac{A_1^2(u)}{2(1-\rho^2)} - \frac{B^2(u)}{2}.
\end{align}

\begin{Corollary}\label{Corollary: specific integral with rho in the condition, second term}
Consider a combination of $\alpha_1$, $\alpha_2$, $\lambda_1$, $\lambda_2$ and $\rho$ that guarantees $A_1(u) = F_1^{-1}(u)-\rho F_2^{-1}(u) \to -\infty$ as $u\to 0^+$.  The behaviour of the integral in (\ref{Route 1}) under such combination of parameters  is then {\adb
\begin{align*}
&  \frac{\alpha_1\sqrt{1-\rho^2}}{\Phi(\lambda_2F_2^{-1}(u)) )}\int^{\frac{A_1(u)}{\sqrt{1-\rho^2}}}_{-\infty} \Phi(x) \inv{\sqrt{2\pi}} e^{-\inv{2}(\alpha_1\sqrt{1-\rho^2} x+(\alpha_2+\rho\alpha_1)F_2^{-1}(u))^2}\,dx\\
\sim & \begin{cases}
\frac{\alpha_1\sqrt{1-\rho^2}e^{-\inv{2}\frac{A^2_1(u)}{1-\rho^2} - \inv{2}B^2(u)}}{2\pi\Phi(\lambda_2F_2^{-1}(u))|F_2^{-1}(u)| |A_1(u)|\beta_{1}}, & \text{if $\beta_{1} >0$;}\\
\frac{\alpha_1(1-\rho^2) e^{-\inv{2}\lambda_2^2 (F_2^{-1}(u))^2}}{2 \sqrt{2\pi}\Phi(\lambda_2F_2^{-1}(u))|A_1(u)| \sqrt{1+\alpha_1^2(1-\rho^2)}}, & \text{if $\beta_{1} =0$};\\
 \frac{\alpha_1e^{-\inv{2}\lambda_2^2 (F_2^{-1}(u))^2}}{\sqrt{2\pi}\Phi(\lambda_2F_2^{-1}(u))|\alpha_1\lambda_2F_2^{-1}(u)|}  , & \text{if $\beta_{1} <0$},
\end{cases}
\end{align*}}
as $u\to 0^+$.
%where $\beta_{1.2}$ was defined in Theorem 1.
\end{Corollary}
\begin{proof}
It is obvious that the integral satisfies Theorem \ref{Theorem: evaluate general integral} with $a(u) = \frac{A_1(u)}{\sqrt{1-\rho^2}}$, $c = \alpha_1\sqrt{1-\rho^2}$, $b(u) = (\alpha_2+\rho\alpha_1)F_2^{-1}(u)$ such that $a(u) \to -\infty$ and 
\begin{equation*}
\frac{b(u)}{a(u)}= \frac{(\alpha_2+\rho\alpha_1)\sqrt{1-\rho^2}}{A_1(u)/F_2^{-1}(u)} = \frac{(\alpha_2+\rho\alpha_1)\sqrt{1-\rho^2}}{\frac{F_1^{-1}(u)}{F_2^{-1}(u)}-\rho}\to \frac{(\alpha_2+\rho\alpha_1)\sqrt{1-\rho^2}}{\gamma_1-\rho} =k,
\end{equation*} 
from (\ref{F ratio}). 
So follows from Theorem \ref{Theorem: evaluate general integral}, (\ref{Lemma: eqn exponents equivalency}), and (\ref{eqn relating b(u) and c^2 with lambda2 and F2inv})--(\ref{eqn stating a(u), b(u) combo has 3 cases}), we have {\adb
\begin{align*}
&  \frac{\alpha_1\sqrt{1-\rho^2}}{\Phi(\lambda_2F_2^{-1}(u))}\int^{\frac{A_1(u)}{\sqrt{1-\rho^2}}}_{-\infty} \Phi(x) \inv{\sqrt{2\pi}} e^{-\inv{2}(\alpha_1\sqrt{1-\rho^2} x+(\alpha_2+\rho\alpha_1)F_2^{-1}(u))^2}\,dx\\
\sim & \begin{cases}
\frac{\alpha_1\sqrt{1-\rho^2}e^{-\inv{2}\frac{A_1^2(u)}{1-\rho^2} - \inv{2}B^2(u)}}{2\pi \Phi(\lambda_2F_2^{-1}(u))|F_2^{-1}(u)| |A_1(u)|\beta_{1}}, & \text{if  $\beta_{1} >0$;}\\
\frac{\alpha_1(1-\rho^2) e^{-\inv{2}\lambda_2^2 (F_2^{-1}(u))^2}}{2 \sqrt{2\pi}\Phi(\lambda_2F_2^{-1}(u))|A_1(u)| \sqrt{1+\alpha_1^2(1-\rho^2)}}, & \text{if $\beta_{1} =0$}\\
\frac{\alpha_1\sqrt{1-\rho^2}\sqrt{[1+\alpha_1^2(1-\rho^2)](1-\rho^2)} e^{-\inv{2}\lambda_2^2 (F_2^{-1}(u))^2}}{\sqrt{2\pi}\Phi(\lambda_2F_2^{-1}(u))|A_1(u)|\left|\alpha_1\sqrt{1-\rho^2}\left(\frac{(\alpha_2+\rho\alpha_1)\sqrt{1-\rho^2}}{\gamma_1-\rho}\right)\right|} , & \text{if $\beta_{1} <0$}
\end{cases}
\end{align*}}
as $u\to 0^+$. We can further simplify the expression in the case of  $\beta_{1}<0$. From Corollary \ref{Corollary: at least one}, 
\begin{equation*}
\frac{A_1(u)}{\gamma_1-\rho} \sim F_2^{-1}(u)
\end{equation*}
so when $\beta_{1}<0$, the result becomes {\adb
\begin{align*}
 & \frac{\alpha_1\sqrt{1-\rho^2}\sqrt{[1+\alpha_1^2(1-\rho^2)](1-\rho^2)} e^{-\inv{2}\lambda_2^2 (F_2^{-1}(u))^2}}{\sqrt{2\pi}\Phi(\lambda_2F_2^{-1}(u))|A_1(u)|\left|\alpha_1\sqrt{1-\rho^2}\left(\frac{(\alpha_2+\rho\alpha_1)\sqrt{1-\rho^2}}{\gamma_1-\rho}\right)\right|} \\
=& \frac{\alpha_1e^{-\inv{2}\lambda_2^2 (F_2^{-1}(u))^2}}{\sqrt{2\pi}\Phi(\lambda_2F_2^{-1}(u))|F_2^{-1}(u)|\left|\alpha_1\left(\frac{\alpha_2+\rho\alpha_1}{\sqrt{1+\alpha_1^2(1-\rho^2)}}\right)\right|} \\
= & \frac{\alpha_1e^{-\inv{2}\lambda_2^2 (F_2^{-1}(u))^2}}{\sqrt{2\pi}\Phi(\lambda_2F_2^{-1}(u))|\alpha_1\lambda_2F_2^{-1}(u)|},
\end{align*}}
as $\lambda_2 = \frac{\alpha_2+\rho\alpha_1}{\sqrt{1+\alpha_1^2(1-\rho^2)}}$ and this completes the proof. 
\end{proof}

Assuming we have $A_1(u) \to -\infty$ as $u\to 0^+$, recall from (\ref{Route 1}) that {\adb
\begin{align}
\notag & P(Z_1\leq F_1^{-1}(u) | Z_2 = F_2^{-1}(u)) \\
%=& \int^{F_{1}^{-1}(u)}_{-\infty} \inv{\sqrt{2\pi(1-\rho^2)}} e^{-\inv{2(1-\rho^2)} (z_1-\rho F_2^{-1}(u))^2} \frac{\Phi(\alpha_1 z_1+\alpha_2F_2^{-1}(u))}{\Phi(\lambda_2 F_2^{-1}(u))} \,dz_1\\
%=& \int^{\frac{F_1^{-1}(u)-\rho F_2^{-1}(u)}{\sqrt{1-\rho^2}}}_{-\infty} \inv{\sqrt{2\pi}}e^{-\frac{x^2}{2}} \frac{\Phi\left(\alpha_1\sqrt{1-\rho^2}x + (\alpha_2+\rho\alpha_1)F_2^{-1}(u)\right)}{\Phi(\lambda_2F_2^{-1}(u))}\,dx,
%\end{align*}
%by applying a change of variable of $x = \frac{z_1-\rho F_2^{-1}(u)}{\sqrt{1-\rho^2}}$;
%\begin{align}
%\notag =& \int^{\frac{F_{1}^{-1}(u)-\rho F_2^{-1}(u)}{\sqrt{1-\rho^2}}}_{-\infty} \frac{\Phi(\alpha_1\sqrt{1-\rho^2} x+(\alpha_2+\rho\alpha_1)F_2^{-1}(u))}{\Phi(\lambda_2 F_2^{-1}(u))} \,d\Phi(x)\\
\notag =& \left[ \Phi(x) \times \frac{\Phi(\alpha_1\sqrt{1-\rho^2} x+(\alpha_2+\rho\alpha_1)F_2^{-1}(u))}{\Phi(\lambda_2 F_2^{-1}(u))}\right]^{\frac{A_1(u)}{\sqrt{1-\rho^2}}}_{-\infty} \\
\notag & \quad - \int^{\frac{A_1(u)}{\sqrt{1-\rho^2}}}_{-\infty} \Phi(x) \frac{\inv{\sqrt{2\pi}} e^{-\inv{2}(\alpha_1\sqrt{1-\rho^2} x+(\alpha_2+\rho\alpha_1)F_2^{-1}(u))^2}}{\Phi(\lambda_2 F_2^{-1}(u))} \times \alpha_1\sqrt{1-\rho^2}\,dx\\
\notag =& \Phi\left(\frac{A_1(u)}{\sqrt{1-\rho^2}}\right)\times \frac{\Phi(B(u))}{\Phi(\lambda_2F_2^{-1}(u))}\\
& - \frac{\alpha_1\sqrt{1-\rho^2}}{\Phi(\lambda_2F_2^{-1}(u))} \int^{\frac{A_1(u)}{\sqrt{1-\rho^2}}}_{-\infty} \Phi(x) \inv{\sqrt{2\pi}} e^{-\inv{2}(\alpha_1\sqrt{1-\rho^2} x+(\alpha_2+\rho\alpha_1)F_2^{-1}(u))^2} \,dx \label{Theorem1:intergral A}
\end{align}}

Corollary \ref{Corollary: specific integral with rho in the condition, second term} has treated the asymptotic  behaviour of  the integral part  in the preceding expression.  We need also to treat the first part. We label this  as Corollary \ref{Corollary: specific integral with rho in the condition, first term}, for convenience rather than due to  its nature.

\begin{Corollary}\label{Corollary: specific integral with rho in the condition, first term}Suppose $A_1(u) = F_1^{-1}(u)-\rho F_2^{-1}(u) \to -\infty$ as $u\to 0^+$.Then  {\adb
%We will need to consider the two terms in (\ref{Theorem1:intergral A}) in separately.
\begin{align*}
&\Phi\left(\frac{A_1(u)}{\sqrt{1-\rho^2}}\right) \frac{\Phi(B(u))}{\Phi(\lambda_2F_2^{-1}(u))}
\sim  \frac{\sqrt{1-\rho^2}e^{-\frac{A_1^2(u)}{2(1-\rho^2)}}
  \Phi(B(u))}{\sqrt{2\pi}|A_1(u)|\Phi(\lambda_2F_2^{-1}(u))}\\
 \sim & 
\begin{cases}
 \frac{\sqrt{1-\rho^2}e^{-\frac{A_1^2(u)}{2(1-\rho^2)}-\inv{2}B^2(u)}}{2\pi|A_1(u)||B(u)|\Phi(\lambda_2F_2^{-1}(u))}, & \text{if $B(u) \to -\infty$};\\
 \frac{\sqrt{1-\rho^2}e^{-\frac{A_1^2(u)}{2(1-\rho^2)}}
 }{2\sqrt{2\pi}|A_1(u)|\Phi(\lambda_2F_2^{-1}(u))}, & \text{if $B(u) \to 0$};\\
  \frac{\sqrt{1-\rho^2}e^{-\frac{A_1^2(u)}{2(1-\rho^2)}}
 }{\sqrt{2\pi}|A_1(u)|\Phi(\lambda_2F_2^{-1}(u))}, & \text{if $B(u) \to \infty$};
\end{cases}
\end{align*}}
as $u\to 0^+$
\end{Corollary}
\begin{proof}
The proof begins by using the expansion of $\Phi(\cdot)$ from \cite{fellerIntroductionProbabilityTheory1968} Chapter VII Lemma 2: 
%(the following will be incorporated into a new Corollary 1 if the updated Theorem 1 holds): when $\lambda_1=0$ we have 
\begin{equation}
\Phi(x) = \inv{\sqrt{2\pi}|x|}e^{-\frac{x^2}{2}}\left(1+O\left(\inv{x^2}\right)\right), \quad \text{for $x<1$}, \label{normal cdf expansion1}
%\quad \text{and}\quad y(u) = -\sqrt{-2\log u\sqrt{-4\pi \log u}},
\end{equation}
Then
\begin{equation*}
\Phi\left(\frac{A_1(u)}{\sqrt{1-\rho^2}}\right) \sim \frac{\sqrt{1-\rho^2}}{\sqrt{2\pi}|A_1(u)|}e^{-\inv{2}\left(\frac{A_1(u)}{\sqrt{1-\rho^2}}\right)^2}
\end{equation*}
when $A_1(u)\to -\infty$ as $u\to 0^+$ and
\begin{equation*}
\Phi(B(u)) = \begin{cases}
\frac{e^{-\inv{2}B^2(u)}}{\sqrt{2\pi |B(u)|}}, & \text{if $B(u) \to -\infty$}; \\
\inv{2}, & \text{if $B(u) \to 0$}; \\
1, & \text{if $B(u) \to \infty$}, 
\end{cases}
\end{equation*}
as $u\to 0^+$. The result follows by combining these expressions.
\end{proof}

%where $\beta_{1.2}$ was defined in Theorem 1.

The purpose of the Corollaries \ref{Corollary: specific integral with rho in the condition, second term} and \ref{Corollary: specific integral with rho in the condition, first term} is  to evaluate (\ref{Route 1}) asymptotically.  By the  argument presented in Lemma \ref{Lemma: at least one condition is positive}, we know that most of the cases would are covered by Corollaries \ref{Corollary: specific integral with rho in the condition, second term} and \ref{Corollary: specific integral with rho in the condition, first term}. Corollary \ref{Corollary: specific integral without rho in the condition}  below will simply cover the remaining possibility, when $\gamma_1 - \rho \leq 0$.

%The corollary below handles integral in the form of (\ref{Route 2}) which requires $\alpha_1>0$. But that's not a problem as In the circumstances that we need the guarantee of $B(u) \to -\infty$ as $u\to 0^+$ as part of Lemma \ref{Lemma: at least one condition is positive}

\begin{Corollary}\label{Corollary: specific integral without rho in the condition}
Suppose that $\alpha_1> 0$ and $B(u) = \alpha_1F_1^{-1}(u)+\alpha_2F_2^{-1}(u) \to -\infty$ as $u\to 0^+$. 
%and $\alpha_1$, $\beta_{1.2}>0$, where $\beta_{1.2}$ is defined in (\ref{beta12 definition}), 
Then (\ref{Route 2}) can be expressed as {\adb
\begin{align*}
 & \inv{\alpha_1\sqrt{1-\rho^2}\Phi(\lambda_2F_2^{-1}(u))}\int^{B(u)}_{-\infty} \inv{\sqrt{2\pi}} e^{-\inv{2}\left(\frac{x}{\alpha_1\sqrt{1-\rho^2}}-\frac{\alpha_2+\rho\alpha_1}{\alpha_1\sqrt{1-\rho^2}}F_2^{-1}(u)\right)^2}\Phi(x)\,dx\\
\sim&  \frac{e^{-\frac{A_1^2(u)}{2(1-\rho^2)} - \inv{2}B^2(u)}}{2\pi \sqrt{1-\rho^2}\Phi(\lambda_2F_2^{-1}(u))|B(u)| |F_2^{-1}(u)|\beta_{1}(u)},
\end{align*}}
as $u\to 0^+$.
\end{Corollary}
\begin{proof}
It is obvious that the integral satisfies the conditions in Theorem \ref{Theorem: evaluate general integral} with $a(u) = B(u)$, $c = \inv{\alpha_1\sqrt{1-\rho^2}}$, $b(u) = - \frac{(\alpha_2+\rho\alpha_1)}{\alpha_1\sqrt{1-\rho^2}}F_2^{-1}(u)$ such that $a(u) \to -\infty$ and 
\begin{equation*}
\frac{b(u)}{a(u)} = \frac{-(\alpha_2+\rho\alpha_1)}{\alpha_1\sqrt{1-\rho^2}\left(\alpha_1\frac{F_1^{-1}(u)}{F_2^{-1}(u)}+\alpha_2\right)}\to \frac{-(\alpha_2+\rho\alpha_1)}{\alpha_1\sqrt{1-\rho^2}\left(\alpha_1\gamma_1+\alpha_2\right)} = k,
\end{equation*} which is a constant as $ u \to 0^{+}$. By Lemma \ref{Lemma: beta1>0}, 
we have $\beta_{1}>0$ so by (\ref{eqn stating a(u), b(u) combo has 3 cases}),  from Theorem \ref{Theorem: evaluate general integral}, we have 
\begin{align*}
& \inv{\alpha_1\sqrt{1-\rho^2}\Phi(\lambda_2F_2^{-1}(u)}\int^{B(u)}_{-\infty} \inv{\sqrt{2\pi}} e^{-\inv{2}\left(\frac{x}{\alpha_1\sqrt{1-\rho^2}}-\frac{\alpha_2+\rho\alpha_1}{\alpha_1\sqrt{1-\rho^2}}F_2^{-1}(u)\right)^2}\Phi(x)\,dx\\
\sim & \frac{e^{-\frac{A^2_1(u)}{2(1-\rho^2)} - \inv{2}B^2(u)}}{2\pi \sqrt{1-\rho^2}\Phi(\lambda_2F_2^{-1}(u))|B(u)| |\beta_{1}F_2^{-1}(u)|},
%= & \frac{e^{-\frac{(F_1^{-1}(u)-\rho F_2^{-1}(u))^2}{2(1-\rho^2)} - \inv{2}(\alpha_1F_1^{-1}(u) + \alpha_2 F_2^{-1}(u))^2}}{2\pi \sqrt{1-\rho^2}\Phi(\lambda_2F_2^{-1}(u))|\alpha_1F_1^{-1}(u)+\alpha_2F_2^{-1}(u)| |F_2^{-1}(u)|\beta_{1}(u)},
\end{align*}
as $u\to 0^+$, from (\ref{special1}).
\end{proof}

The following Theorem \ref{Thm: main result} summarises the main result of this note.  Recall that at least one of $A_1(u) = F_1^{-1}(u)-\rho F_2^{-1}(u) \to-\infty$ and  $B(u) = \alpha_1F_1^{-1}(u) + \alpha_2F_2^{-1}(u) \to -\infty$ as $ u \to 0^{+},$ according to Corollary \ref{Corollary: at least one}. 
%different to what we had in Lemma \ref{Lemma: at least one condition is positive}, which is depends on the behaviour of $F_1^{-1}(u)-\rho F_2^{-1}(u)$ and $F_1^{-1}(u) + \frac{\alpha_2}{\alpha_1}F_2^{-1}(u)$ instead.

%Although 
%
%Next, we shall consider the behaviour of the integrals in (\ref{Route 1}) and (\ref{Route 2}) and the results are summarised into the following theorem. 

\begin{Theorem} \label{Thm: main result} Suppose that $\beta_{1} = \lim_{u\to 0^+} \beta_{1}(u)$ is as defined as in (\ref{beta12 definition}) and $\alpha_1\ne 0$. Consider in the sequel $ u \to 0^+$.

\begin{enumerate}[label=(\alph*)]
\item Assuming $B(u) \to -\infty$ \underline{and} $A_1(u) \to -\infty$, $0$ or $\infty$, then $\beta_1>0$ and 
\begin{equation*}
P(Z_1\leq F_1^{-1}(u) | Z_2 = F_2^{-1}(u)) \sim \frac{e^{-\frac{A^2_1(u)}{2(1-\rho^2)}- \inv{2}B^2(u)}}{2\pi \sqrt{1-\rho^2} \beta_{1}|F_2^{-1}(u)|\Phi(\lambda_2F_2^{-1}(u))|B(u)|};
\end{equation*}

\item Assuming $B(u) \to 0$, then $A_1(u) \to -\infty$, and 
\begin{align*}
P(Z_1\leq F_1^{-1}(u) | Z_2 = F_2^{-1}(u)) \sim 
\frac{\sqrt{1-\rho^2}e^{-\frac{A_1^2(u)}{2(1-\rho^2)}}}{2\sqrt{2\pi}|A_1(u)|\Phi(\lambda_2F_2^{-1}(u))};
\end{align*}

\item Assuming $B(u) \to \infty$,  then   $A_1(u)  \to -\infty$, and:
\begin{align*}
P(Z_1\leq F_1^{-1}(u) | Z_2 = F_2^{-1}(u)) \sim 
\frac{\sqrt{1-\rho^2}e^{-\frac{A_1^2(u)}{2(1-\rho^2)}}}{\sqrt{2\pi}|A_1(u)|\Phi(\lambda_2F_2^{-1}(u))}.
\end{align*}
\end{enumerate}
\end{Theorem} 

The proof can be found in Appendix \ref{Appendix: Proof of Theorem 2}. Notice how the structure of the asymptotic results in the statement  of Theorem \ref{Thm: main result} reflects that of Theorem \ref{Theorem: evaluate general integral}.

\section{Asymptotic forms for components in Theorem \ref{Thm: main result}}

To proceed, we need to express the asymptotic  forms on the right of Theorem \ref{Thm: main result} (a) (b) (c)  in terms of the ``$\log$"- type functions  for each of the various 
sign combinations of $\lambda_1, \lambda_2$, which means we need 
the  asymptotic expressions of each of 
$$  |\lambda_2 F_2^{-1}(u)|\Phi\left(\lambda_2F_2^{-1}(u)\right),\quad  e^{-\frac{A^2_1(u)}{2(1-\rho^2)}},\quad  e^{-\frac{B^2(u)}{2}},  $$
as $ u \to 0^+$.  By noting when $\lambda_2>0$,
$$|\lambda_2 F_2^{-1}(u)|\Phi\left(\lambda_2F_2^{-1}(u)\right) \sim \inv{\sqrt{2\pi}} e^{-\inv{2}\lambda_2^2(F_2^{-1}(u))^2},$$
using (\ref{normal cdf expansion1}), it is possible to provide a unified treatment on all three expressions based on the common factor which has a general form $e^{-\frac{1}{2}\left(G_1F_1^{-1}(u) - G_2F_2^{-1}(u)\right)^2}$ where $G_1$ and $G_2$ are two constants so that  
\begin{equation*}
G_1F_1^{-1}(u) - G_2F_2^{-1}(u) = 
\begin{cases}
\lambda_2F_2^{-1}(u), & \text{if  $G_1 = 0$, $G_2= -\lambda_2$;}\\
\frac{F_1^{-1}(u) - \rho F_2^{1}(u)}{\sqrt{1-\rho^2}} = \frac{A_1(u)}{\sqrt{1-\rho^2}}, & \text{if $G_1 = \inv{\sqrt{1-\rho^2}}$, $G_2= \frac{\rho} {\sqrt{1-\rho^2}}$;}\\
\alpha_1F_1^{-1}(u) + \alpha_2F_2^{-1}(u) = B(u), & \text{if $G_1 = \alpha_1$, $G_2 = -\alpha_2$.}
\end{cases}
\end{equation*}

%We are considering  $e^{-\inv{2}\lambda_2^2(F_2^{-1}(u))^2}$ as 
%(the following will be incorporated into a new Corollary 1 if the updated Theorem 1 holds): when $\lambda_1=0$ we have 
%\begin{equation}
%\Phi(x) = \inv{\sqrt{2\pi}|x|}e^{-\frac{x^2}{2}}\left(1+O\left(\inv{x^2}\right)\right), \quad \text{for $x<1$}, \label{normal cdf expansion}
%\quad \text{and}\quad y(u) = -\sqrt{-2\log u\sqrt{-4\pi \log u}},
%\end{equation}
%We can give intially a unified treatment by noting that 
%respectively,  so we shall proceed with the general form $G_1F_1^{-1}(u) - GF_2^{-1}(u)$.

This unified treatment has been summarised into the following Lemma.
\begin{Lemma}\label{lemma: general form of G1F1-GF2}
For any constants $G_1$ and $G_2$, we have 
\begin{align*}
 \left|G_1F_1^{-1}(u)-G_2F_2^{-1}(u)\right| \sim& \left|K_{2,1}(G_1\gamma_1-G_2)\right| \sqrt{-2\log u}, \quad \text{as $u\to 0^+$;} \\
e^{-\inv{2}\left(G_1F_1^{-1}(u)-G_2 F_2^{-1}(u)\right)^2}
\notag =& \tau_1 u^{\theta}(-\log u)^{\tau_2}\left(1+O\left(\frac{\left[\log|\log u|\right]^2}{\log u}\right)\right),
%u^{G_1^2K_{2,2}^2(\gamma_1-G_2)^2}(-\log u)^{\text{constant}_1}\times \text{constant}_2
\end{align*}
where 
\begin{align*}
\theta =& K_{2,1}^2(G_1\gamma_1-G_2)^2;\\
\tau_1 = & e^{2\left[K_{2,1}^2K_{2,3}(G_1\gamma_1-G_2)^2+K_{2,1}^2C_{1,2}G_1\gamma_1(G_1\gamma_1-G_2)\right]};\\
\tau_2 =& 2\left[K_{2,1}^2K_{2,2}(G_1\gamma_1-G_2)^2+K_{2,1}^2C_{1,1}G_1\gamma_1(G_1\gamma_1-G_2)\right],
\end{align*}
and $\gamma_1$, $K_{2.1}$, $K_{2,2}$, $K_{2,3}$, $C_{1,1}$ $C_{1,2}$ are as defined in Lemma \ref{Lemma: detailed F inverse}.
\end{Lemma}

The proof of Lemma \ref{lemma: general form of G1F1-GF2} can be found in Appendix \ref{Appendix: Proof of Lemma: general form of G1F1-GF2}.

%Hence, using the asymptotic behaviour of the numerator from (\ref{generic nominator}), we have 
%\begin{align*}
%& \frac{e^{-\inv{2}G_1^2(F_2^{-1}(u))^2\left(\frac{F_1^{-1}(u)}{F_2^{-1}(u)}-G_2\right)^2}}{\left|G_1F_2^{-1}(u)\left(\frac{F_1^{-1}(u)}{F_2^{-1}(u)}-G_2\right)\right|}\\
%\sim & u^{G_1^2K_2^2(\gamma(1,2)-G_2)^2}(-\log u)^{\text{constant}_1-\inv{2}}\times \frac{\text{constant}_2}{\left| \sqrt{2}G_1K_2(\gamma(1,2)-G_2)\right|}
%\end{align*}

Using the results from Lemma \ref{lemma: general form of G1F1-GF2}, we can now present the  asymptotic result on $|\lambda_2 F_2^{-1}(u)|\Phi(\lambda_2F_2^{-1}(u))$ for $\lambda_2>0$; $e^{-\frac{A^2_1(u)}{2(1-\rho^2)}} = e^{-\inv{2}\left(\frac{F_1^{-1}(u)-\rho F_2^{-1}(u)}{\sqrt{1-\rho^2}}\right)^2}$; \& $e^{-\frac{B^2(u)}{2}}= e^{-\inv{2}(\alpha_1F_1^{-1}(u)+\alpha_2F_2^{-1}(u))^2}$ as $u\to 0^+$  in the following three Lemmas and their proofs can be found in Appendix \ref{Appendix: Proof of Lemma: general form of G1F1-GF2} and its subsections.

Notice that for $e^{-\frac{A^2_1(u)}{2(1-\rho^2)}}$ and $e^{-\frac{B^2(u)}{2}}$, there are 8 distinct cases to be treated for these two expressions:$ [1]\, \lambda_1 >0, \lambda_2 >0; [2]\, \lambda_1 <0, \lambda_2 >0; [3]\, \lambda_1>0, \lambda_2 <0;[4]\,  \lambda_1=0, \lambda_2 >0; [5]\, \lambda_1 >0, \lambda_2=0;[6]\,  \lambda_1 <0, \lambda_2 <0; [7]\, \lambda_1<0, \lambda_2 =0; [8]\, \lambda_1=0, \lambda_2 <0.$
 The expression  for the ``constant" multiplier  $ \tau_1$ is very similar in the pair of cases [2] and [4] while $\theta, \tau_2$ are the same. A similar situation holds for the pair of cases [3] and [5],  and for the triple of cases [6], [7], and [8]. So in  the sequel we will  specify $\theta, \tau_2$ only, for the four cases (the last three amalgamated): $$ \lambda_1 >0, \lambda_2 >0; \lambda _1 \leq 0, \lambda_2 >0; \lambda_1 >0, \lambda_2 \leq 0; \lambda_1 \leq 0, \lambda_2 \leq 0.$$
\begin{Lemma} \label{Lemma: Phi in RV}
For $\lambda_2>0$, 
\begin{equation*}
|\lambda_2 F_2^{-1}(u)|\Phi(\lambda_2F_2^{-1}(u)) \sim  
 \inv{\sqrt{2\pi}} u^{\frac{\lambda_2^2}{1+\lambda_2^2}}|\log u|^{\frac{\lambda_2^2}{1+\lambda_2^2}}(2\pi\lambda_2)^{\frac{\lambda_2^2}{1+\lambda_2^2}}, \quad \text{as $u\to 0^+$.}
%\frac{\sqrt{1+\lambda_2^2}(2\pi\lambda_2)^{\frac{\lambda_2^2}{1+\lambda_2^2}}}{2\sqrt{\pi}\lambda_2}u^{\frac{\lambda_2^2}{1+\lambda_2^2}}|\log u|^{\frac{\lambda_2^2}{1+\lambda_2^2}-\inv{2}}
\end{equation*}
\end{Lemma} 

%{\bf  Thomas, please edit to end of this Section 5, using hand-written notes, pp.3,4,5, dated 01.10.19.}

\begin{Lemma} \label{Lemma: exp A in RV}
As $u\to 0^+$, we have  

\begin{align}
\notag \frac{|A_1(u)|}{\sqrt{1-\rho^2}} &= \frac{|F_1^{-1}(u)- \rho F_2^{-1}(u)|}{\sqrt{1-\rho^2}} \\
&\sim 
\begin{cases}
\inv{\sqrt{1-\rho^2}}\left|\inv{\sqrt{1+\lambda_1^2}}-\frac{\rho}{\sqrt{1+\lambda_2^2}}\right|\sqrt{-2\log u}, & \text{if $\lambda_1>0$, $\lambda_2>0$};\\
\inv{\sqrt{1-\rho^2}}\left|1-\frac{\rho}{\sqrt{1+\lambda_2^2}}\right|\sqrt{-2\log u}, & \text{if $\lambda_1\leq 0$, $\lambda_2>0$};\\
\inv{\sqrt{1-\rho^2}}\left|\inv{\sqrt{1+\lambda_1^2}}-\rho\right|\sqrt{-2\log u}, 
&\text{if $\lambda_1>0$, $\lambda_2\leq 0$};\\
\sqrt{\frac{1-\rho}{1+\rho}}\sqrt{-2\log u},& \text{if $\lambda_1\leq 0$, $\lambda_2\leq 0$ (};
\end{cases} \label{eqn: detailed A in RV}
\end{align}
and 
\begin{align*}
& e^{-\frac{A^2_1(u)}{2(1-\rho^2)}} =  e^{-\inv{2}\left(\frac{F_1^{-1}(u)-\rho F_2^{-1}(u)}{\sqrt{1-\rho^2}}\right)^2}\\
\sim& 
\begin{cases} 
(2\pi\lambda_2)^{\inv{1-\rho^2}\left(\inv{\sqrt{1+\lambda_1^2}}-\frac{\rho}{\sqrt{1+\lambda_2^2}}\right)^2}\times \left(\frac{\lambda_1}{\lambda_2}\right)^{\inv{(1-\rho^2)\sqrt{1+\lambda_1^2}}\left(\inv{\sqrt{1+\lambda_1^2}}-\frac{\rho}{\sqrt{1+\lambda_2^2}}\right)}  \\
\quad \times u^{\inv{1-\rho^2}\left(\inv{\sqrt{1+\lambda_1^2}}-\frac{\rho}{\sqrt{1+\lambda_2^2}}\right)^2}(-\log u)^{\inv{1-\rho^2}\left(\inv{\sqrt{1+\lambda_1^2}}-\frac{\rho}{\sqrt{1+\lambda_2^2}}\right)^2}, & \text{if $\lambda_1>0$, $\lambda_2>0$};\\
(2\pi\lambda_2)^{\inv{1-\rho^2}\left(1-\frac{\rho}{\sqrt{1+\lambda_2^2}}\right)^2}\left(2\lambda_2\sqrt{\pi}\right)^{-\inv{1-\rho^2}\left(1-\frac{\rho}{\sqrt{1+\lambda_2^2}}\right)}\\
\quad \times u^{\inv{1-\rho^2}\left(1-\frac{\rho}{\sqrt{1+\lambda_2^2}}\right)^2}(-\log u)^{\inv{1-\rho^2}\left(1-\frac{\rho}{\sqrt{1+\lambda_2^2}}\right)^2-\inv{2(1-\rho^2)}\left(1-\frac{\rho}{\sqrt{1+\lambda_2^2}}\right)}, & \text{if $\lambda_1<0$, $\lambda_2>0$};\\
\pi^{\inv{2(1-\rho^2)}\left(\inv{\sqrt{1+\lambda_1^2}}-\rho\right)^2}(2\lambda_1\sqrt{\pi})^{\inv{2(1-\rho^2)\sqrt{1+\lambda_1^2}}\left(\inv{\sqrt{1+\lambda_1^2}}-\rho\right)} \\
\quad \times  u^{\inv{1-\rho^2}\left(\inv{\sqrt{1+\lambda_1^2}}-\rho\right)^2}(-\log u)^{\inv{2(1-\rho^2)}\left(\inv{\sqrt{1+\lambda_1^2}}-\rho\right)^2+\inv{2(1-\rho^2)\sqrt{1+\lambda_1^2}}\left(\inv{\sqrt{1+\lambda_1^2}}-\rho\right)}, & \text{if $\lambda_1>0$, $\lambda_2<0$}; \\
(2\pi\lambda_2)^{\inv{1-\rho^2}\left(1-\frac{\rho}{\sqrt{1+\lambda_2^2}}\right)^2}(\lambda_2\sqrt{\pi})^{-\inv{1-\rho^2}\left(1-\frac{\rho}{\sqrt{1+\lambda_2^2}}\right)}\\
\quad \times u^{\inv{1-\rho^2}\left(1-\frac{\rho}{\sqrt{1+\lambda_2^2}}\right)^2}(-\log u)^{\inv{1-\rho^2}\left(1-\frac{\rho}{\sqrt{1+\lambda_2^2}}\right)^2 - \inv{2(1-\rho^2)}\left(1-\frac{\rho}{\sqrt{1+\lambda_2^2}}\right)}, & \text{if $\lambda_1=0$, $\lambda_2>0$};\\
(4\pi)^{\inv{2(1-\rho^2)}\left(\inv{\sqrt{1+\lambda_1^2}}-\rho\right)^2}(\lambda_1\sqrt{\pi})^{\inv{(1-\rho^2)\sqrt{1+\lambda_1^2}}\left(\inv{\sqrt{1+\lambda_1^2}}-\rho\right)}\\
\quad \times u^{\inv{1-\rho^2}\left(\inv{\sqrt{1+\lambda_1^2}}-\rho\right)^2}(-\log u)^{ \inv{2(1-\rho^2)}\left(\inv{\sqrt{1+\lambda_1^2}}-\rho\right)^2+\inv{2(1-\rho^2)\sqrt{1+\lambda_1^2}}\left(\inv{\sqrt{1+\lambda_1^2}}-\rho\right)}, & \text{if $\lambda_1>0$, $\lambda_2=0$};\\
\pi^{\inv{2}\left(\frac{1-\rho}{1+\rho}\right)} u^{\frac{1-\rho}{1+\rho}}(-\log u)^{\inv{2}\left(\frac{1-\rho}{1+\rho}\right)}, &\text{if $\lambda_1<0$, $\lambda_2<0$};\\
(4\pi)^{\inv{2}\left(\frac{1-\rho}{1+\rho}\right)}2^{-\inv{1+\rho}}u^{\frac{1-\rho}{1+\rho}}(-\log u)^{\inv{2}\left(\frac{1-\rho}{1+\rho}\right)}, &\text{$\lambda_1<0$, $\lambda_2=0$};\\
\pi^{\inv{2}\left(\frac{1-\rho}{1+\rho}\right)}2^{\inv{1+\rho}} u^{\frac{1-\rho}{1+\rho}}(-\log u)^{\inv{2}\left(\frac{1-\rho}{1+\rho}\right)}, &\text{$\lambda_1=0$, $\lambda_2<0$.}
\end{cases}
\end{align*}
\end{Lemma}

\begin{Lemma}
 \label{Lemma: exp B in RV}
 As $u\to 0^+$, we have 
 \begin{align}
\notag  |B(u)| & = |\alpha_1F_1^{-1}(u)+\alpha_2F_2^{-1}(u)| \\
& \sim 
\begin{cases}
\left|\frac{\alpha_1}{\sqrt{1+\lambda_1^2}}+\frac{\alpha_2}{\sqrt{1+\lambda_2^2}}\right|\sqrt{-2\log u}, 
&\text{if $\lambda_1>0$, $\lambda_2>0$};\\
\left|\alpha_1+\frac{\alpha_2}{\sqrt{1+\lambda_2^2}}\right|\sqrt{-2\log u}, 
& \text{if $\lambda_1\leq 0$, $\lambda_2>0$};\\
\left|\frac{\alpha_1}{\sqrt{1+\lambda_1^2}}+\alpha_2\right|\sqrt{-2\log u}, 
&\text{if $\lambda_1>0$, $\lambda_2\leq 0$};\\
|\alpha_1+\alpha_2|\sqrt{-2\log u}, &\text{if $\lambda_1\leq 0$, $\lambda_2\leq 0$.}
\end{cases} \label{eqn: detailed B in RV}
\end{align}

\begin{align*}
& e^{-\frac{B^2(u)}{2}} =  e^{-\inv{2}\left(\alpha_1F_1^{-1}(u)+\alpha_2F_2^{-1}(u)\right)^2}\\
=& 
\begin{cases}
(2\pi\lambda_2)^{\left(\frac{\alpha_1}{\sqrt{1+\lambda_1^2}}+\frac{\alpha_2}{\sqrt{1+\lambda_2^2}}\right)^2}\left(\frac{\lambda_1}{\lambda_2}\right)^{\frac{\alpha_1}{\sqrt{1+\lambda_1^2}}\left(\frac{\alpha_1}{\sqrt{1+\lambda_1^2}}+\frac{\alpha_2}{\sqrt{1+\lambda_2^2}}\right)}\\
\quad \times u^{\left(\frac{\alpha_1}{\sqrt{1+\lambda_1^2}}+\frac{\alpha_2}{\sqrt{1+\lambda_2^2}}\right)^2}
(-\log u)^{\left(\frac{\alpha_1}{\sqrt{1+\lambda_1^2}}+\frac{\alpha_2}{\sqrt{1+\lambda_2^2}}\right)^2}, 
&\text{if $\lambda_1>0$, $\lambda_2>0$};\\
(2\pi\lambda_2)^{\left(\alpha_1+\frac{\alpha_2}{\sqrt{1+\lambda_2^2}}\right)^2}\left(2\lambda_2\sqrt{\pi}\right)^{-\alpha_1\left(\alpha_1+\frac{\alpha_2}{\sqrt{1+\lambda_2^2}}\right)}\\
\quad \times u^{\left(\alpha_1+\frac{\alpha_2}{\sqrt{1+\lambda_2^2}}\right)^2}(-\log u)^{\left(\alpha_1+\frac{\alpha_2}{\sqrt{1+\lambda_2^2}}\right)^2-\inv{2}\alpha_1\left(\alpha_1+\frac{\alpha_2}{\sqrt{1+\lambda_2^2}}\right)}, 
&\text{if $\lambda_1<0$, $\lambda_2>0$};\\
\pi^{\inv{2}\left(\frac{\alpha_1}{\sqrt{1+\lambda_1^2}}+\alpha_2\right)^2}\times (2\lambda_1\sqrt{\pi})^{\frac{\alpha_1}{\sqrt{1+\lambda_1^2}}\left(\frac{\alpha_1}{\sqrt{1+\lambda_1^2}}+\alpha_2\right)}\\
\quad \times u^{\left(\frac{\alpha_1}{\sqrt{1+\lambda_1^2}}+\alpha_2\right)^2}(-\log u)^{\inv{2}\left(\frac{\alpha_1}{\sqrt{1+\lambda_1^2}}+\alpha_2\right)^2+\inv{2}\frac{\alpha_1}{\sqrt{1+\lambda_1^2}}\left(\frac{\alpha_1}{\sqrt{1+\lambda_1^2}}+\alpha_2\right)}, 
&\text{if $\lambda_1>0$, $\lambda_2<0$};\\
(2\pi\lambda_2)^{\left(\alpha_1+\frac{\alpha_2}{\sqrt{1+\lambda_2^2}}\right)^2}(\lambda_2\sqrt{\pi})^{-\alpha_1\left(\alpha_1+\frac{\alpha_2}{\sqrt{1+\lambda_2^2}}\right)}\\
\quad \times u^{\left(\alpha_1+\frac{\alpha_2}{\sqrt{1+\lambda_2^2}}\right)^2}(-\log u)^{\left(\alpha_1+\frac{\alpha_2}{\sqrt{1+\lambda_2^2}}\right)^2-\frac{\alpha_1}{2}\left(\alpha_1+\frac{\alpha_2}{\sqrt{1+\lambda_2^2}}\right)}, 
& \text{$\lambda_1=0$, $\lambda_2>0$};\\
(4\pi)^{\inv{2}\left(\frac{\alpha_1}{\sqrt{1+\lambda_1^2}}+\alpha_2\right)^2}(\lambda_1\sqrt{\pi})^{\frac{\alpha_1}{\sqrt{1+\lambda_1^2}}\left(\frac{\alpha_1}{\sqrt{1+\lambda_1^2}}+\alpha_2\right)}\\
\quad \times u^{\left(\frac{\alpha_1}{\sqrt{1+\lambda_1^2}}+\alpha_2\right)^2}(-\log u)^{\inv{2}\left(\frac{\alpha_1}{\sqrt{1+\lambda_1^2}}+\alpha_2\right)^2+\inv{2}\frac{\alpha_1}{\sqrt{1+\lambda_1^2}}\left(\frac{\alpha_1}{\sqrt{1+\lambda_1^2}}+\alpha_2\right)},
& \text{if $\lambda_1>0$, $\lambda_2=0$};\\
\pi^{\frac{(\alpha_1+\alpha_2)^2}{2}}u^{(\alpha_1+\alpha_2)^2}(-\log u)^{\inv{2}(\alpha_1+\alpha_2)^2}, 
&\text{if $\lambda_1<0$, $\lambda_2<0$};\\
(4\pi)^{\frac{(\alpha_1+\alpha_2)^2}{2}}2^{-\alpha_1(\alpha_1+\alpha_2)}u^{(\alpha_1+\alpha_2)^2}(-\log u)^{\inv{2}(\alpha_1+\alpha_2)^2}, &\text{if $\lambda_1<0$, $\lambda_2=0$};\\
\pi^{\frac{(\alpha_1+\alpha_2)^2}{2}}2^{\alpha_1(\alpha_1+\alpha_2)}u^{(\alpha_1+\alpha_2)^2}(-\log u)^{\inv{2}(\alpha_1+\alpha_2)^2}, &\text{if $\lambda_1=0$, $\lambda_2<0$}.
\end{cases}
\end{align*}
\end{Lemma}      .

\section{Summation.Theorem 3}

Recall, from (\ref{defn:tail dependence13}),  (\ref{defn:tail dependence133})
that  it is sufficient for us to find a value of $\theta > 0$
which satisfies
\begin{equation*}
\frac{d C(u,u)}{du} = P(Z_1\leq F_1^{-1}(u)|Z_2 = F_2^{-1}(u)) + P(Z_2 \leq F_2^{-1}(u) | Z_1 = F_1^{-1}(u)) = u^{\theta} L(u)
\label{result}
\end{equation*} for some slowly varying function $L(u).$ In the preceding sections we have focussed without loss of generality on the asymptotic behaviour of the first of these two summands, and in effect found such an expression for it for all combinations of $\lambda_1, \lambda_2$, with $u^{\theta} L(u) $ of the form  $ \tau_1 u^{\theta}(-\log u)^{\tau_2}\left(1+O\left(\frac{\left[\log|\log u|\right]^2}{\log u}\right)\right)$.  By interchange of subscripts the general 
form will be the same for the second summand. although even for a given combination of  $\lambda_1, \lambda_2$, 
the parameter values $\tau_1, \theta, \tau_2 $ will differ for each summand.  But in any case it is clear that the result of the summation for any combination of $\lambda_1, \lambda_2$, will   still be of form:
$$u^{\theta}L(u)= \tau_1 u^{\theta}(-\log u)^{\tau_2}\left(1+O\left(\frac{\left[\log|\log u|\right]^2}{\log u}\right)\right),$$ and by abuse of notation will continue to use this notation $\tau_1, \theta, \tau_2$, to express the asymptotic result of this sum of two terms.

\begin{Theorem} \label{Theorem: Summation}
Let $Z\sim SN_2(\bs{\alpha}, R)$ with $R = \left(\begin{smallmatrix} 1 & & \rho \\ \\ \rho & & 1 \end{smallmatrix}\right)$. Then 
\begin{equation*}
\frac{dC(u,u)}{du} = u^{\theta}L(u)
\end{equation*}
as $u\to 0^+$, where 
\begin{enumerate}
\item when $\lambda_1$, $\lambda_2<0$ or equivalently $\alpha_1+\rho\alpha_2<0$ \& $\alpha_2+\rho\alpha_1<0$, then\footnote{Notice that the conditions here cover the condition $\alpha_1=\alpha_2 = \alpha < 0$ of \cite{fungTailAsymptoticsBivariate2016}, re-expressed in our  Theorem \ref{thm:fungTailAsymptoticsBivariate2016}} 
\begin{equation*}
\theta = \frac{1-\rho}{1+\rho} \quad \left[ = \frac{2}{1+\rho}-1 \right], \quad L(u) \sim \sqrt{\frac{1+\rho}{1-\rho}}\pi^{-\frac{\rho}{1+\rho}}(-\log u)^{-\frac{\rho}{1+\rho}};
\end{equation*} 

\item when $\lambda_1=0$, $\lambda_2<0$ or equivalently $\alpha_1+\rho\alpha_2=0$ \& $\alpha_2+\rho\alpha_1<0$, then 
\begin{equation*}
\theta = \frac{1-\rho}{1+\rho}\quad \left[ = \frac{2}{1+\rho}-1 \right], \quad L(u) \sim \sqrt{\frac{1+\rho}{1-\rho}}\pi^{-\frac{\rho}{1+\rho}}2^{\inv{1+\rho}}(-\log u)^{-\frac{\rho}{1+\rho}};
\end{equation*}
\item when $\lambda_1<0$, $\lambda_2>0$ or equivalently $\alpha_1+\rho\alpha_2<0$ \& $\alpha_2+\rho\alpha_1>0$, then 
\begin{enumerate}
\item if $\alpha_1+\frac{\alpha_2}{\sqrt{1+\lambda_2^2}}>0$, then 
\begin{align*}
\theta = &\inv{1-\rho^2}\left( \inv{\sqrt{1+\lambda_2^2}}-\rho \right)^2+\left( \alpha_1+\frac{\alpha_2}{\sqrt{1+\lambda_2^2}} \right)^2\quad \\
& \left[ = \frac{(1-\rho^2)+\left(1+ \sqrt{1+\lambda_2^2} \right)^2}{(1-\rho^2)(1+\lambda_2^2)}+ \left( \alpha_1+\frac{\alpha_2}{\sqrt{1+\lambda_2^2}} \right)^2 - 1\right]\\
L(u) \sim &  \left(-\log u\right)^{\inv{2(1-\rho^2)}\left(\inv{\sqrt{1+\lambda_2^2}}-\rho\right)^2+\inv{2(1-\rho^2)\sqrt{1+\lambda_2^2}}\left(\inv{\sqrt{1+\lambda_2^2}}-\rho\right)+\left(\alpha_1+\frac{\alpha_2}{\sqrt{1+\lambda_2^2}}\right)^2-\inv{2}\alpha_1\left(\alpha_1+\frac{\alpha_2}{\sqrt{1+\lambda_2^2}}\right)-1}\\
& \quad \times \frac{(2\lambda_2)^{\left(\alpha_1+\frac{\alpha_2}{\sqrt{1+\lambda_2^2}} \right)^2 -\alpha_1\left(\alpha_1+\frac{\alpha_2}{\sqrt{1+\lambda_2^2}}\right)+\inv{2(1-\rho^2)\sqrt{1+\lambda_2^2}}\left( \inv{\sqrt{1+\lambda_2^2}}-\rho \right)}}{4\sqrt{1-\rho^2}\left| \alpha_1+\frac{\alpha_2}{\sqrt{1+\lambda_2^2}} \right|}\\
& \quad \times \pi^{\left(\alpha_1+\frac{\alpha_2}{\sqrt{1+\lambda_2^2}} \right)^2 -\frac{\alpha_1}{2}\left(\alpha_1+\frac{\alpha_2}{\sqrt{1+\lambda_2^2}}\right)+\inv{4(1-\rho^2)\sqrt{1+\lambda_2^2}}\left( \inv{\sqrt{1+\lambda_2^2}}-\rho \right)+\inv{2(1-\rho^2)}\left( \inv{\sqrt{1+\lambda_2^2}}-\rho \right)^2-1}\\
& \quad \times \left[ \inv{|\beta_2|}+\frac{(2\lambda_2\sqrt{\pi})^{\inv{2(1-\rho^2)\sqrt{1+\lambda_2^2}}\left( \inv{\sqrt{1+\lambda_2^2}}-\rho \right)}}{|\beta_1|} \right],
%
%\left(-\log u\right)^{\inv{2(1-\rho^2)}\left(\inv{\sqrt{1+\lambda_2^2}}-\rho\right)^2+\inv{2(1-\rho^2)\sqrt{1+\lambda_2^2}}\left(\inv{\sqrt{1+\lambda_2^2}}-\rho\right)+\left(\alpha_1+\frac{\alpha_2}{\sqrt{1+\lambda_2^2}}\right)^2-\inv{2}\alpha_1\left(\alpha_1+\frac{\alpha_2}{\sqrt{1+\lambda_2^2}}\right)-1}\\
%& \quad \times \Biggl[ \frac{(2\pi\lambda_2)^{\left(\alpha_1+\frac{\alpha_2}{\sqrt{1+\lambda_2^2}}\right)^2}\left(\lambda_2 \sqrt{\pi}\right)^{\inv{(1-\rho^2)\sqrt{1+\lambda_2^2}}\left(\inv{\sqrt{1+\lambda_2^2}}-\rho\right)-\alpha_1\left(\alpha_1+\frac{\alpha_2}{\sqrt{1+\lambda_2^2}}\right)}(4\pi)^{\inv{2(1-\rho^2)}\left(\inv{\sqrt{1+\lambda_2^2}}-\rho\right)^2}
%}{2\pi\sqrt{1-\rho^2}\left|\alpha_1+\frac{\alpha_2}{\sqrt{1+\lambda_2^2}}\right|\times \beta_2}\\
%& \quad + \frac{\lambda_2(2\pi\lambda_2)^{\inv{1-\rho^2}\left(1-\frac{\rho}{\sqrt{1+\lambda_2^2}}\right)^2+\left(\alpha_1+\frac{\alpha_2}{\sqrt{1+\lambda_2^2}}\right)^2-\frac{\lambda_2^2}{1+\lambda_2^2}}}{2\sqrt{\pi(1-\rho^2)}\left|\alpha_1+\frac{\alpha_2}{\sqrt{1+\lambda_2^2}}\right|\beta_1\times (\lambda_2\sqrt{\pi})^{\inv{1-\rho^2}\left(1-\frac{\rho}{\sqrt{1+\lambda_2^2}}\right)+\alpha_1\left(\alpha_1+\frac{\alpha_2}{\sqrt{1+\lambda_2^2}}\right)}}\Biggr] 
\end{align*}
%[The slowly varying function is quite different to that of Padoan.]

%\begin{enumerate}
%\item if $\beta_2>0$, then 
%\item if $\beta_2=0$, then 
%\begin{align*}
%\theta = \lambda_1^2, \quad L(u) = \frac{\pi^{\frac{\lambda_1^2}{2}-\inv{2}}|\alpha_2|(1-\rho^2)(-\log u)^{\frac{\lambda_1^2}{2}-\inv{2}}}{4\left| \inv{\sqrt{1+\lambda_2^2}}-\rho \right|\sqrt{1+\alpha_2^2(1-\rho^2)}};
%\end{align*}
%[This case is not covered by Padoan?]
%\item if $\beta_2<0$, then 
%{\adb
%\begin{align*}
%\theta = \lambda_1^2, \quad L(u) = \frac{\pi^{\frac{\lambda_1^2}{2}-\inv{2}}(-\log u)^{\frac{\lambda_1^2}{2}}}{2|\lambda_1|}
%\end{align*}
%}
%[This case is not covered by Padoan?]
%\end{enumerate}
\item if $\alpha_1+\frac{\alpha_2}{\sqrt{1+\lambda_2^2}} =0$, then 
{\adb
\begin{align*}
\theta & = \inv{1-\rho^2}\left( \inv{\sqrt{1+\lambda_2^2}}-\rho \right)^2, \\
L(u) &\sim \frac{\sqrt{1-\rho^2}}{4} (2\lambda_2)^{\inv{2(1-\rho^2)\sqrt{1+\lambda_2^2}}\left( \inv{\sqrt{1+\lambda_2^2}-\rho}\right)}\pi^{\inv{2(1-\rho^2)}\left( \inv{\sqrt{1+\lambda_2^2}}-\rho \right)^2-\inv{2}}\\
& \quad \times \left[ (2\lambda\pi)^{\inv{2(1-\rho^2)\sqrt{1+\lambda_2^2}}\left( \inv{\sqrt{1+\lambda_2^2}}-\rho \right)}+ \pi^{-\inv{4(1-\rho^2)\sqrt{1+\lambda_2^2}}\left( \inv{\sqrt{1+\lambda_2^2}}-\rho \right)}
 \right] \\
\notag & \quad \times (-\log u)^{\inv{2(1-\rho^2)}\left( \inv{\sqrt{1+\lambda_2^2}}-\rho \right)^2+\inv{2(1-\rho^2)\sqrt{1+\lambda_2^2}}\left( \inv{\sqrt{1+\lambda_2^2}}-\rho \right)-\inv{2}}
%
%\Biggl[\frac{|\lambda_2|(2\pi\lambda_2)^{\inv{1-\rho^2}\left( 1-\frac{\rho}{\sqrt{1+\lambda_2^2}} \right)^2-\frac{\lambda_2^2}{1+\lambda_2^2}}(2\lambda_2\sqrt{\pi})^{-\inv{1-\rho^2}\left( 1-\frac{\rho}{\sqrt{1+\lambda_2^2}} \right)}}{2\left| \sqrt{1+\lambda_2^2}-\rho \right|}\\
%& \quad +  \frac{\pi^{\inv{2(1-\rho^2)}\left( \inv{\sqrt{1+\lambda_2^2}}-\rho \right)^2-\inv{2}}(2\lambda_2\sqrt{\pi})^{\inv{2(1-\rho^2)\sqrt{1+\lambda_2^2}}\left( \inv{\sqrt{1+\lambda_2^2}}-\rho \right)}}{4\left| \inv{\sqrt{1+\lambda_2^2}}-\rho \right|} \Biggr]\\
%& \quad \times (-\log u)^{\inv{2(1-\rho^2)}\left( \inv{\sqrt{1+\lambda_2^2}}-\rho \right)^2+\inv{2(1-\rho^2)\sqrt{1+\lambda_2^2}}\left( \inv{\sqrt{1+\lambda_2^2}}-\rho \right)-\inv{2}}
\end{align*}
}
%Looks like the limit of the $\beta_2>0$ case by letting $\alpha_1+\frac{\alpha_2}{\sqrt{1+\lambda_2^2}} \to 0$ but that still won't account for one has $-\inv{2}$ and another one has $-1$ in the $(-\log u)$ term? [This case is definitely not covered by Padoan.] 
\item if $\alpha_1+\frac{\alpha_2}{\sqrt{1+\lambda_2^2}} <0$, then 
{\adb
\begin{align*}
\theta &= \inv{1-\rho^2}\left( \inv{\sqrt{1+\lambda_2^2}}-\rho \right)^2 \quad \left[ = \frac{(1-\rho^2)+\left(1+\sqrt{1+\lambda_2^2}\right)^2}{(1-\rho^2)(1+\lambda_2^2)} - 1\right],\\
L(u) &\sim \frac{\sqrt{1-\rho^2}}{2} (2\lambda_2)^{\inv{2(1-\rho^2)\sqrt{1+\lambda_2^2}}\left( \inv{\sqrt{1+\lambda_2^2}-\rho}\right)}\pi^{\inv{2(1-\rho^2)}\left( \inv{\sqrt{1+\lambda_2^2}-\rho}-\rho \right)^2-\inv{2}}\\
& \quad \times \left[ (2\lambda\pi)^{\inv{2(1-\rho^2)\sqrt{1+\lambda_2^2}}\left( \inv{\sqrt{1+\lambda_2^2}}-\rho \right)}+ \pi^{-\inv{4(1-\rho^2)\sqrt{1+\lambda_2^2}}\left( \inv{\sqrt{1+\lambda_2^2}}-\rho \right)}
 \right] \\
\notag & \quad \times (-\log u)^{\inv{2(1-\rho^2)}\left( \inv{\sqrt{1+\lambda_2^2}}-\rho \right)^2+\inv{2(1-\rho^2)\sqrt{1+\lambda_2^2}}\left( \inv{\sqrt{1+\lambda_2^2}}-\rho \right)-\inv{2}}
\end{align*}
}
\end{enumerate}
\item if $\lambda_1=0$, $\lambda_2>0$ or equivalently $\alpha_1+\rho\alpha_2=0$ \& $\alpha_2+\rho\alpha_1>0$, then 
\begin{enumerate}
\item if $\alpha_1+\frac{\alpha_2}{\sqrt{1+\lambda_2^2}} > 0$, then 
\end{enumerate}
{\adb
\begin{align*}
\theta &= \inv{1-\rho^2}\left( \inv{\sqrt{1+\lambda_2^2}}-\rho \right)^2+\left( \alpha_1+\frac{\alpha_2}{\sqrt{1+\lambda_2^2}} \right)^2 \\
& \left[ = \frac{(1-\rho^2)+\left(1+ \sqrt{1+\lambda_2^2} \right)^2}{(1-\rho^2)(1+\lambda_2^2)}+ \left( \alpha_1+\frac{\alpha_2}{\sqrt{1+\lambda_2^2}} \right)^2 - 1\right]\\
L(u) &\sim  \left(-\log u\right)^{\inv{2(1-\rho^2)}\left(\inv{\sqrt{1+\lambda_2^2}}-\rho\right)^2+\inv{2(1-\rho^2)\sqrt{1+\lambda_2^2}}\left(\inv{\sqrt{1+\lambda_2^2}}-\rho\right)+\left(\alpha_1+\frac{\alpha_2}{\sqrt{1+\lambda_2^2}}\right)^2-\inv{2}\alpha_1\left(\alpha_1+\frac{\alpha_2}{\sqrt{1+\lambda_2^2}}\right)-1}
\\
& \quad \times \frac{(2\pi\lambda_2)^{\left(\alpha_1+\frac{\alpha_2}{\sqrt{1+\lambda_2^2}}\right)^2}(\lambda_2\sqrt{\pi})^{\inv{1-\rho^2)\sqrt{1+\lambda_2^2}}\left( \inv{\sqrt{1+\lambda_2^2}}-\rho \right)-\alpha_1\left(\alpha_1+\frac{\alpha_2}{\sqrt{1+\lambda_2^2}}\right)}\left( 2\sqrt{\pi}\right)^{\inv{1-\rho^2}\left( \inv{\sqrt{1+\lambda_2^2}}-\rho \right)^2}}{2\sqrt{\pi(1-\rho^2)}\left| \alpha_1+\frac{\alpha_2}{\sqrt{1+\lambda_2^2}} \right|}\\
&\quad \times \left[ \inv{\sqrt{\pi}\beta_2}+\inv{\beta_1} \right]
\end{align*}
}
\begin{enumerate}
\item[(b)] if $\alpha_1+\frac{\alpha_2}{\sqrt{1+\lambda_2^2}}=0$, then 
\end{enumerate}
{\adb
\begin{align*}
\theta & =  \inv{1-\rho^2}\left(\inv{\sqrt{1+\lambda_2^2}}-\rho \right)^2 \quad \left[ = \frac{(1-\rho^2)+\left( 1+\sqrt{1+\lambda_2^2} \right)^2}{(1-\rho^2)(1+\lambda_2^2)}-1 \right]\\
L(u) &\sim \frac{\sqrt{1-\rho^2}}{2\sqrt{\pi}}\left(\lambda_2\sqrt{\pi}\right)^{\inv{(1-\rho^2)\sqrt{1+\lambda_2^2}}\left( \inv{\sqrt{1+\lambda_2^2}}-\rho \right)}(2\sqrt{\pi})^{\inv{1-\rho^2}\left( \inv{\sqrt{1+\lambda_2^2}}-\rho \right)^2}\\
& \quad \times  \left[ \inv{\sqrt{1+\lambda_2^2}-\rho}+\inv{\left| \inv{\sqrt{1+\lambda_2^2}}-\rho \right|}
 \right](-\log u)^{\inv{2(1-\rho^2)}\left( \inv{\sqrt{1+\lambda_2^2}}-\rho \right)^2+\inv{2(1-\rho^2)\sqrt{1+\lambda_2^2}}\left( \inv{\sqrt{1+\lambda_2^2}}-\rho \right)-\inv{2}};
\end{align*}
}
\begin{enumerate}
\item[(c)] if $\alpha_1+\frac{\alpha_2}{\sqrt{1+\lambda_2^2}}<0$, then 
\end{enumerate}
{\adb
\begin{align*}
\theta &= \inv{1-\rho^2}\left( \inv{\sqrt{1+\lambda_2^2}}-\rho \right)^2 \quad \left[ = \frac{(1-\rho^2)+\left(1+\sqrt{1+\lambda_2^2}\right)^2}{(1-\rho^2)(1+\lambda_2^2)} - 1\right],\\
L(u) &\sim \frac{\sqrt{1-\rho^2}}{\sqrt{\pi}}\left(\lambda_2\sqrt{\pi}\right)^{\inv{(1-\rho^2)\sqrt{1+\lambda_2^2}}\left( \inv{\sqrt{1+\lambda_2^2}}-\rho \right)}(2\sqrt{\pi})^{\inv{1-\rho^2}\left( \inv{\sqrt{1+\lambda_2^2}}-\rho \right)^2}\\
& \quad \times  \left[ \inv{\sqrt{1+\lambda_2^2}-\rho}+\inv{\left| \inv{\sqrt{1+\lambda_2^2}}-\rho \right|}
 \right](-\log u)^{\inv{2(1-\rho^2)}\left( \inv{\sqrt{1+\lambda_2^2}}-\rho \right)^2+\inv{2(1-\rho^2)\sqrt{1+\lambda_2^2}}\left( \inv{\sqrt{1+\lambda_2^2}}-\rho \right)-\inv{2}};
\end{align*}
}
\item if $\lambda_1>0$, $\lambda_2>0$ or equivalently $\alpha_1+\rho\alpha_2>0$ \& $\alpha_2+\rho\alpha_1>0$, then \footnote{Notice that the conditions here cover the condition $\alpha_1=\alpha_2 = \alpha > 0$ of \cite{fungTailAsymptoticsBivariate2016}, re-expressed in our  Theorem \ref{thm:fungTailAsymptoticsBivariate2016}.}
{\adb
\begin{align*}
\theta &= \inv{1-\rho^2}\left( \inv{\sqrt{1+\lambda_2^2}}-\frac{\rho}{\sqrt{1+\lambda_1^2}} \right)^2+\left( \frac{\alpha_1}{\sqrt{1+\lambda_1^2}}+\frac{\alpha_2}{\sqrt{1+\lambda_2^2}} \right)^2-\frac{\lambda_1^2}{1+\lambda_1^2}\\
& \quad \left[ =\inv{1-\rho^2}\left( \frac{1+\alpha_1^2(1-\rho^2)}{1+\lambda_1^2}+\frac{1+\alpha_2^2(1-\rho^2)}{1+\lambda_2^2}+\frac{2\left(\alpha_1\alpha_2(1-\rho^2)-\rho\right)}{\sqrt{(1+\lambda_1^2)(1+\lambda_2^2)}} \right) -1\right]\\
L(u) &\sim \left[ \inv{\beta_1}+\inv{\beta_2}\right] \times \frac{(2\pi)^{\inv{1-\rho^2}\left( \inv{\sqrt{1+\lambda_2^2}}-\frac{\rho}{\sqrt{1+\lambda_1^2}} \right)^2+\left( \frac{\alpha_1}{\sqrt{1+\lambda_1^2}}+\frac{\alpha_2}{\sqrt{1+\lambda_2^2}} \right)^2-\frac{\lambda_1^2}{1+\lambda_1^2}-\inv{2}}}{\sqrt{2(1-\rho^2)}\left| \frac{\alpha_1}{\sqrt{1+\lambda_1^2}}+\frac{\alpha_2}{\sqrt{1+\lambda_2^2}} \right|}\\
& \quad \times 
\lambda_1^{\inv{(1-\rho^2)\sqrt{1+\lambda_1^2}}\left( \inv{\sqrt{1+\lambda_1^2}}-\frac{\rho}{\sqrt{1+\lambda_2^2}} \right)+\frac{\alpha_1}{\sqrt{1+\lambda_1^2}}\left( \frac{\alpha_1}{\sqrt{1+\lambda_1^2}}+\frac{\alpha_2}{\sqrt{1+\lambda_2^2}}\right)}\\
& \quad \times 
\lambda_2^{\inv{(1-\rho^2)\sqrt{1+\lambda_2^2}}\left( \inv{\sqrt{1+\lambda_2^2}}-\frac{\rho}{\sqrt{1+\lambda_1^2}} \right)+\frac{\alpha_2}{\sqrt{1+\lambda_2^2}}\left( \frac{\alpha_1}{\sqrt{1+\lambda_1^2}}+\frac{\alpha_2}{\sqrt{1+\lambda_2^2}}\right)}\\
& \quad \times (-\log u)^{\inv{1-\rho^2}\left( \inv{\sqrt{1+\lambda_2^2}}-\frac{\rho}{\sqrt{1+\lambda_1^2}} \right)^2+\left( \frac{\alpha_1}{\sqrt{1+\lambda_1^2}}+\frac{\alpha_2}{\sqrt{1+\lambda_2^2}} \right)^2-\frac{\lambda_1^2}{1+\lambda_1^2}-\inv{2}}\\
\end{align*}
}
\end{enumerate} 
\end{Theorem} 

\appendix

\section{Proof of Lemma \ref{Lemma: at least one condition is positive}} \label{Appendix: Proof of Lemma 2}

\begin{proof}
%Based on the behaviour $F_i^{-1}(u)$ stated above, it is equivalent to show that at least one of $\frac{F_{1}^{-1}(u)}{F_2^{-1}(u)}-\rho$ and $\frac{F_{1}^{-1}(u)}{F_2^{-1}(u)}+\frac{\alpha_2}{\alpha_1}$ is asymptotically strictly positive for any combination of $\alpha_1$, $\alpha_2$, $\lambda_1$, $\lambda_2$ and $\rho$. 

Since from  Lemma \ref{Lemma: detailed F inverse} the behaviour of $F_i^{-1}(u)$ depends solely on the sign of $\lambda_i$, we split the proof according to  the different combinations of $\lambda_1$ and $\lambda_2$. Using  (\ref{detail F ratio}) we have 
\begin{align}
%=&  \begin{cases}
%\frac{\sqrt{1+\lambda_2^2}}{\sqrt{1+\lambda_1^2}}\left(1+\frac{\log(\lambda_1/\lambda_2)}{2\log u}+O\left(\left(\frac{\log(-\log u)}{\log u}\right)^2\right)\right), & \text{when $\lambda_1$, $\lambda_2>0$};\\
% \sqrt{1+\lambda_2^2}\left(1-\frac{\log\left(2\lambda_2\sqrt{-\pi\log u}\right)}{2\log u}+O\left(\left(\frac{\log(-\log u)}{\log u}\right)^2\right)\right), &\text{when $\lambda_1<0$, $\lambda_2>0$;}\\
% \inv{\sqrt{1+\lambda_1^2}}\left(1+\frac{\log\left(2\lambda_1\sqrt{-\pi\log u}\right)}{2\log u}+O\left(\left(\frac{\log(-\log u)}{\log u}\right)^2\right)\right), & \text{when $\lambda_1>0$, $\lambda_2<0$};\\
% \sqrt{1+\lambda_2^2}\left(1-\frac{\log\left(\lambda_2\sqrt{-\pi\log u}\right)}{2\log u}+O\left(\left(\frac{\log(-\log u)}{\log u}\right)^2\right)\right), &\text{when $\lambda_1=0$, $\lambda_2>0$};\\
% \inv{\sqrt{1+\lambda_1^2}}\left(1+\frac{\log\left(\lambda_1\sqrt{-\pi\log u}\right)}{2\log u}+O\left(\left(\frac{\log(-\log u)}{\log u}\right)^2\right)\right), &\text{when $\lambda_1>0$, $\lambda_2=0$;}\\
% 1 + O\left(\left(\frac{\log(-\log u)}{\log u}\right)^2\right), & \text{when $\lambda_1$, $\lambda_2<0$:}\\
% 1 - \frac{\log 2}{2\log u} + O\left(\left(\frac{\log(-\log u)}{\log u}\right)^2\right), & \text{when $\lambda_1<0$, $\lambda_2=0$;}\\
% 1 + \frac{\log 2}{2\log u} + O\left(\left(\frac{\log(-\log u)}{\log u}\right)^2\right), &\text{when $\lambda_1=0$, $\lambda_2<0$,}
%\end{cases} \\
\frac{F_1^{-1}(u)}{F_2^{-1}(u)} \sim & \gamma_1 = \begin{cases}
\sqrt{\frac{1+\lambda_2^2}{1+\lambda_1^2}}, & \text{when } \lambda_1, \lambda_2 > 0; \\
\sqrt{1+\lambda_2^2} & \text{when } \lambda_1 \leq 0 \text{ \& } \lambda_2 > 0; \\
\inv{\sqrt{1+\lambda_1^2}}& \text{when } \lambda_1 > 0 \text{ \& } \lambda_2 \leq 0; \\
1 & \text{when } \lambda_1, \lambda_2 \leq 0. 
\end{cases} \label{F ratio}
\end{align}
as $u\to 0^+.$ Therefore it is obvious whenever $-1<\rho \leq 0$, or in the cases  ($\lambda_1\leq 0$ \& $\lambda_2>0$) or  ($\lambda_1$, $\lambda_2\leq 0$), that  $\gamma_1 -\rho>0$ . 

Next we  show for the case of $\lambda_1>0$, $\lambda_2\leq 0$ and $0<\rho<1$ that  $\gamma_1 -\rho =  \inv{\sqrt{1+\lambda_1^2}}-\rho >0$.  In this case we also gain $\alpha_2< 0$. Ths  is a crucial condition to the rest of this proof. It is easy to see that   when $\alpha_2=0, 0 < \rho <1 $ it is impossible to  have $\lambda_1 = \alpha_1> 0$, $\lambda_2 = \frac{\rho\alpha_1}{\sqrt{1+\alpha_1^2(1-\rho^2)}} \leq 0$. In order for us to show that $\alpha_2<0$, suppose the opposite is true with $\alpha_2 > 0.$ Then 
\begin{equation*}
\lambda_2\leq 0 \quad \Rightarrow \quad \alpha_2+\rho\alpha_1 \leq 0 \quad \Rightarrow \quad 
0< \alpha_2 \leq \rho(-\alpha_1) \quad \Rightarrow \quad 0<\alpha_2 \leq \rho|\alpha_1|  \quad \Rightarrow \quad
0<\inv{\rho} \leq |\alpha_1|/\alpha_2.
\end{equation*}
Similarly,
\begin{equation*}
\lambda_1> 0 \quad \Rightarrow \quad \alpha_1+\rho\alpha_2 >  0 \quad \Rightarrow \quad 
\rho \alpha_2 > |\alpha_1| \quad \Rightarrow \quad 1 > \rho >|\alpha_1|/\alpha_2\geq 0. 
\end{equation*}
By combining the two inequalities we get 
\begin{equation*}
1>\rho> |\alpha_1|/\alpha_2 \geq 1/\rho >0 \quad \Rightarrow \quad 1 > \rho > 1/\rho >0
\end{equation*}
which is a contradiction. This means that when $\lambda_1>0$, $\lambda_2\leq 0$ and $0<\rho<1$, we have  $\alpha_2 < 0$. Notice that this further implies $\alpha_1>0$ as $\lambda_1>0 \Rightarrow \alpha_1+\rho\alpha_2> 0$ but $\rho\alpha_2\leq 0$.

Now
\begin{eqnarray*}
0< \gamma_1-\rho = \inv{\sqrt{1+\lambda_1^2}}-\rho = \frac{\sqrt{1+\alpha_2^2(1-\rho^2)}-\rho\sqrt{1+\alpha_1^2+2\rho\alpha_1\alpha_2+\alpha_2^2}}{\sqrt{1+\alpha_1^2+2\rho\alpha_1\alpha_2+\alpha_2^2}}.
\end{eqnarray*}
%[There was a mistake here in the previous version of this note. As $\lambda_1 = \frac{\alpha_1+\rho\alpha_2}{\sqrt{1+\alpha_2^2(1-\rho^2)}}$, so 
%\begin{equation*}
%1+\lambda_1^2 = 1+\frac{\alpha_1^2+2\rho\alpha_1\alpha_2+\rho^2\alpha_2^2}{1+\alpha_2^2(1-\rho^2)} 
%= \frac{1+\alpha_1^2+2\rho\alpha_1\alpha_2+\alpha_2^2}{1+\alpha_2^2(1-\rho^2)}
%\end{equation*}
%and 
%\begin{equation*}
%\inv{\sqrt{1+\lambda_1^2}} = \frac{\sqrt{1+\alpha_2^2(1-\rho^2)}}{\sqrt{1+\alpha_1^2+2\rho\alpha_1\alpha_2+\alpha_2^2}} \text{ instead of } \frac{\sqrt{1+\alpha_1^2(1-\rho^2)}}{\sqrt{1+\alpha_1^2+2\rho\alpha_1\alpha_2+\alpha_2^2}}]
%\end{equation*}
since {\adb
\begin{eqnarray*}
&& \left(\sqrt{1+\alpha_2^2(1-\rho^2)}-\rho\sqrt{1+\alpha_1^2+2\rho\alpha_1\alpha_2+\alpha_2^2}\right)\times \left(\sqrt{1+\alpha_2^2(1-\rho^2)}+\rho\sqrt{1+\alpha_1^2+2\rho\alpha_1\alpha_2+\alpha_2^2}\right)\\
&=& 1+\alpha_2^2(1-\rho^2)-\rho^2(1+\alpha_1^2+2\rho\alpha_1\alpha_2+\alpha_2^2)\\
&=& 1-\rho^2+\alpha_2^2-\rho^2\alpha_1^2-2\rho^2\alpha_2^2-2\rho^3\alpha_1\alpha_2\\
&=& 1-\rho^2+(\alpha_2+\rho\alpha_1)(\alpha_2(1-\rho^2)-\rho(\alpha_1+\rho\alpha_2))>0
\end{eqnarray*}}
and  we know that $\alpha_2<0$, $\alpha_1+\rho\alpha_2>0$, $\alpha_2+\rho\alpha_1\leq 0$ and $0<\rho<1$. 

In the remaining the case  $\lambda_1= \frac{\alpha_1+\rho\alpha_2}{\sqrt{1+\alpha_2^2(1-\rho^2)}}>0$ \& $\lambda_2=\frac{\alpha_2+\rho\alpha_1}{\sqrt{1+\alpha_1^2(1-\rho^2)}}>0$  so that $\alpha_1+\rho\alpha_2>0$ and $\alpha_2+\rho\alpha_1>0$. Adding these two together we have 
\begin{equation}
\alpha_1+\alpha_2+\rho(\alpha_1+\alpha_2) = (1+\rho)(\alpha_1+\alpha_2)>0 \quad \Rightarrow \quad \alpha_1+\alpha_2>0. \label{eqn: lambdas positive implies sum of alphas positive}
\end{equation} 
Now 
\begin{equation*}
\gamma_1 = \sqrt{\frac{1+\lambda_2^2}{1+\lambda_1^2}} = \sqrt{\frac{1+\frac{(\alpha_2+\rho\alpha_1)^2}{1+\alpha_1^2(1-\rho^2)}}{1+\frac{(\alpha_1+\rho\alpha_2)^2}{1+\alpha_2^2(1-\rho^2)}}}\quad \Rightarrow 
\quad \gamma_1-\rho = \sqrt{\frac{1+\alpha_2^2(1-\rho^2)}{1+\alpha_1^2(1-\rho^2)}}-\rho.
\end{equation*}
Thus $\gamma_1-\rho>0$ if $-1<\rho\leq 0$, or if $|\alpha_2| \geq |\alpha_1|$, since then $\gamma_1\geq 1$. If $|\alpha_2|< |\alpha_1|$, and $0<\rho<1$, then 
{\adb
\begin{align*}
\gamma_1-\rho & >0,  \quad \text{if $\alpha_2^2>\alpha_1^2\rho^2-1$};\\
& =0, \quad \text{if $\alpha_2^2=\alpha_1^2\rho^2-1$};\\
& <0, \quad \text{if $\alpha_2^2<\alpha_1^2\rho^2-1$}.
\end{align*}
}
Now if $\alpha_1\leq 0$, $|\alpha_2|<|\alpha_1|$ would contradict (\ref{eqn: lambdas positive implies sum of alphas positive}), so we must have $|\alpha_2|<|\alpha_1|$ and $\alpha_1>0$. Thus when 
$\lambda_1>0$, $\lambda_2>0$, whenever $\gamma_1-\rho\leq 0$, we have  $\alpha_1>0$.

Finally, when $\lambda_1>0$, $\lambda_2>0$, 
%This further implies that at least one of $\alpha_1$ or $\alpha_2$ is positive here. Notice that we would show primarily $\frac{B(u)}{F_2^{-1}(u)}\to \alpha_1\gamma_1+\alpha_2$ is a positive constant here. Focusing on $\alpha_1\gamma_1+\alpha_2 $,
 we have  {\adb
\begin{equation}
\alpha_1\gamma_1+\alpha_2  = \alpha_1\sqrt{\frac{1+\alpha_2^2(1-\rho^2)}{1+\alpha_1^2(1-\rho^2)}}+\alpha_2 = \frac{\alpha_1\sqrt{1+\alpha_2^2(1-\rho^2)}+\alpha_2\sqrt{1+\alpha_1^2(1-\rho^2)}}{\sqrt{1+\alpha_1^2(1-\rho^2)}}, \label{F inverse ratio 1}
 \end{equation}}and we show that this is  always strictly positive when $\lambda_1>0$, $\lambda_2>0$.

From (\ref{eqn: lambdas positive implies sum of alphas positive}) at  least one of $\alpha_1, \alpha_2$ is strictly positive.
Assume without loss of generality, that $\alpha_1 >0.$  If $\alpha_2 \geq 0$ (\ref{F inverse ratio 1}) is clearly strictly positive.  If  $\alpha_2 < 0$ then from  (\ref{eqn: lambdas positive implies sum of alphas positive}) $\alpha_1 > |\alpha_2|,$ and  (\ref{F inverse ratio 1}) is still strictly positive, since its  numerator is ${\sqrt{\alpha_1^2+\alpha_1^2\alpha_2^2(1-\rho^2)}-\sqrt{\alpha_2^2+\alpha_2^2\alpha_1^2(1-\rho^2)}}.$ 

\end{proof}

\section{Proof of Theorem 1.} \label{Appendix: Proof of Theorem 1}
%as $u\to 0^+$.

\begin{proof}

%[The following trick was first developed in [2014-07-18]]. Using integration-by-parts [with $e^{-\beta(u)y}$] on the integral in (\ref{trick integral}): 
The proof begins by using the expansion of $\Phi(\cdot)$ from Feller (1968) Chapter VII Lemma 2: 
%(the following will be incorporated into a new Corollary 1 if the updated Theorem 1 holds): when $\lambda_1=0$ we have 
\begin{equation}
\Phi(x) = \inv{\sqrt{2\pi}|x|}e^{-\frac{x^2}{2}}\left(1+O\left(\inv{x^2}\right)\right), \quad \text{for $x<1$}, \label{normal cdf expansion}
%\quad \text{and}\quad y(u) = -\sqrt{-2\log u\sqrt{-4\pi \log u}},
\end{equation}
and the integral becomes
{\adb
\begin{align}
& \int^{a(u)}_{-\infty}\Phi(x)\inv{\sqrt{2\pi}} e^{-\inv{2}(cx +b(u))^2}\,dx \label{Theorem: evaluate general integral, initial}\\
\notag \sim & \int^{a(u)}_{-\infty} \inv{\sqrt{2\pi}|x|}e^{-\frac{x^2}{2}} \inv{\sqrt{2\pi}} e^{-\inv{2}(cx + b(u))^2}\,dx\\
%, \quad \text{using the expansion of $\Phi(\cdot)$;}\\
\notag = & \inv{2\pi} \int^{a(u)}_{-\infty} \inv{|x|} e^{-\frac{x^2}{2}} e^{-\inv{2}\left(c^2x^2+2cb(u)+b^2(u)\right)}\,dx \\
\notag =& \frac{e^{-\inv{2}b^2(u)}}{2\pi} \int^{a(u)}_{-\infty} \inv{|x|} e^{-\inv{2}\left((1+c^2)x^2+2cb(u)\right)}\,dx \\
\notag =& \frac{e^{-\inv{2}b^2(u) + \frac{c^2b^2(u)}{2(1+c^2)}}}{2\pi}\int^{a(u)}_{-\infty} \inv{|x|}e^{-\inv{2}(1+c^2)\left(x+\frac{cb(u)}{1+c^2}\right)^2}\,dx.
\end{align}}
Apply a change of variable $z = \sqrt{1+c^2}(x-a(u))$ and the expression becomes:
\begin{align}
 =& \frac{e^{-\frac{b^2(u)}{2(1+c^2)}}}{2\pi\sqrt{1+c^2}}\int^{0}_{-\infty} \inv{\left|\frac{z}{\sqrt{1+c^2}}+a(u)\right|}e^{-\inv{2}\left(z+a(u)\sqrt{1+c^2}+ \frac{cb(u)}{\sqrt{1+c^2}}\right)^2}\,dz  \label{Theorem: evaluate general integral, case3}\\
=& \frac{e^{-\frac{b^2(u)}{2(1+c^2)}}}{2\pi\sqrt{1+c^2}|a(u)|}J(u), \label{Theorem: evaluate general integral, integral1}
\end{align}
where 
\begin{equation}
J(u) = \int^{0}_{-\infty} \inv{\left|\frac{z}{a(u)\sqrt{1+c^2}}+1\right|}e^{-\inv{2}\left(z+v(u)\right)^2}\,dz = \int^{0}_{-\infty} \inv{\left(\frac{z}{a(u)\sqrt{1+c^2}}+1\right)}e^{-\inv{2}\left(z+v(u)\right)^2}\,dz 
\label{Theorem: evaluate general integral, J term}
\end{equation}
as $z$, $a(u)<0 \Rightarrow z/a(u)>0$, and $v(u) = \frac{a(u)(1+c^2)+cb(u)}{\sqrt{1+c^2}}$. As $-\infty<z<0$ and $a(u)<0$ for all $u$, $J(u)$, as defined in (\ref{Theorem: evaluate general integral, J term}), is bounded above by
\begin{align}
J(u) &\leq \int^{0}_{-\infty}e^{-\inv{2}\left(z+v(u)\right)^2}\,dz = \sqrt{2\pi}\Phi\left(v(u)\right).\label{Theorem: evaluate general integral, J term ub}
%= \sqrt{2\pi}\Phi\left(\frac{a(u)(1+c^2)+cb(u)}{\sqrt{1+c^2}}\right). 
\end{align}
To get the lower bound of (\ref{Theorem: evaluate general integral, J term}), we use the fact that for $x>0$,
\begin{equation*}
e^x = 1+x+\frac{x^2}{2}+\ldots\quad \Rightarrow \quad \frac{e^x}{1+x} \geq \frac{1+x}{1+x} = 1 \quad \Rightarrow \quad \inv{1+x} \geq e^{-x},
\end{equation*}
and $J(u)$,  as defined in (\ref{Theorem: evaluate general integral, J term}), becomes {\adb
\begin{align}
\notag J(u) 
%\frac{e^{-\frac{b^2(u)}{2(1+c^2)}}}{2\pi\sqrt{1+c^2}}\int^{0}_{-\infty} \inv{\left|\frac{z}{\sqrt{1+c^2}}+a(u)\right|}e^{-\inv{2}\left(z+\frac{a(u)(1+c^2)+ cb(u)}{\sqrt{1+c^2}}\right)^2}\,dz\\
% = & \int^{0}_{-\infty} \inv{\left(1+\frac{z}{a(u)\sqrt{1+c^2}}\right)}e^{-\inv{2}\left(z+v(u)\right)^2}\,dz, \quad \text{as $z$, $a(u)<0 \Rightarrow z/a(u)>0$;}\\
\notag \geq & \int^{0}_{-\infty} e^{-\inv{2}\left(z+v(u)\right)^2- \frac{z}{a(u)\sqrt{1+c^2}}}\,dz\\
\notag  = & \int^{0}_{-\infty} e^{-\inv{2}\left(z^2+2z\left(v(u)+\inv{a(u)\sqrt{1+c^2}}\right)+ \left(v(u)\right)^2\right)}\,dz \\
\notag  =& e^{\inv{2}\left(v(u)+\inv{a(u)\sqrt{1+c^2}}\right)^2-\inv{2}\left(v(u)\right)^2}\int^{0}_{-\infty} e^{-\inv{2}\left(z+v(u)+\inv{a(u)\sqrt{1+c^2}}\right)^2}\,dz\\
=& \sqrt{2\pi} e^{\inv{2}\left(v(u)+\inv{a(u)\sqrt{1+c^2}}\right)^2-\inv{2}\left(v(u)\right)^2}\Phi\left(v(u)+\inv{a(u)\sqrt{1+c^2}}\right).
%=& \frac{e^{-\frac{b^2(u)}{2(1+c^2)}+\inv{2}\left(\frac{a(u)(1+c^2)+cb(u)}{\sqrt{1+c^2}}+\inv{a(u)\sqrt{1+c^2}}\right)^2-\inv{2}\left(\frac{a(u)(1+c^2)+cb(u)}{\sqrt{1+c^2}}\right)^2}}{\sqrt{2\pi}\sqrt{1+c^2}|a(u)|}\Phi\left(\frac{a(u)(1+c^2)+cb(u)}{\sqrt{1+c^2}}+\inv{a(u)\sqrt{1+c^2}}\right). 
\label{Theorem: evaluate general integral, integral1lb}
\end{align}}
This means that 
\begin{equation*}
\sqrt{2\pi} e^{\inv{2}\left(v(u)+\inv{a(u)\sqrt{1+c^2}}\right)^2-\inv{2}\left(v(u)\right)^2}\Phi\left(v(u)+\inv{a(u)\sqrt{1+c^2}}\right) \leq J(u) \leq \sqrt{2\pi}\Phi(v(u)).
\end{equation*}
Thus, the lower and upper bounds for (\ref{Theorem: evaluate general integral, integral1}) are 
%As a results, (\ref{Lemma2:integral1}) is bounded by: 
\begin{align}
\frac{e^{-\frac{b^2(u)}{2(1+c^2)}+\inv{2}\left(v(u)+\inv{a(u)\sqrt{1+c^2}}\right)^2-\inv{2}\left(v(u)\right)^2}\Phi\left(v(u)+\inv{a(u)\sqrt{1+c^2}}\right)}{\sqrt{2\pi}\sqrt{1+c^2}|a(u)|} \quad \& \quad 
%&< \frac{e^{-\frac{b^2(u)}{2(1+c^2)}}}{2\pi\sqrt{1+c^2}|a(u)|}J(u) 
\frac{e^{-\frac{b^2(u)}{2(1+c^2)}}\Phi\left(v(u)\right)}{\sqrt{2\pi}\sqrt{1+c^2}|a(u)|},
\label{Theorem: evaluate general integral, integral1bounds}
\end{align}
respectively. 
So the behaviour of the bounds depends on the behaviour of $v(u).$ 

If $v(u)\to -\infty$ as $u\to 0^+$, then the upper bound of (\ref{Theorem: evaluate general integral, integral1}) (i.e. the second expression in 
(\ref{Theorem: evaluate general integral, integral1bounds})) is asymptotically equivalent to 
\begin{equation*}
 \frac{e^{-\frac{b^2(u)}{2(1+c^2)}}}{\sqrt{2\pi}\sqrt{1+c^2}|a(u)|}\Phi\left(v(u)\right) \sim  \frac{e^{-\frac{b^2(u)}{2(1+c^2)}}}{\sqrt{2\pi}\sqrt{1+c^2}|a(u)|} \inv{\sqrt{2\pi}|v(u)|} e^{-\inv{2}\left(v(u)\right)^2} = \frac{e^{-\frac{b^2(u)}{2(1+c^2)}-\inv{2}\left(v(u)\right)^2}}{2\pi\sqrt{1+c^2}|a(u)v(u)|}, 
\end{equation*}
by using (\ref{normal cdf expansion}). Similarly, the lower bound of (\ref{Theorem: evaluate general integral, integral1}) (i.e. the first expression in (\ref{Theorem: evaluate general integral, integral1bounds})) is asymptotically equivalent to 
\begin{align*}
& \frac{e^{-\frac{b^2(u)}{2(1+c^2)}+\inv{2}\left(v(u)+\inv{a(u)\sqrt{1+c^2}}\right)^2-\inv{2}\left(v(u)\right)^2}}{\sqrt{2\pi}\sqrt{1+c^2}|a(u)|}\Phi\left(v(u)+\inv{a(u)\sqrt{1+c^2}}\right)\\
\sim& \frac{e^{-\frac{b^2(u)}{2(1+c^2)}+\inv{2}\left(v(u)+\inv{a(u)\sqrt{1+c^2}}\right)^2-\inv{2}\left(v(u)\right)^2}}{\sqrt{2\pi}\sqrt{1+c^2}|a(u)|}\times \inv{\sqrt{2\pi}\left|v(u)+\inv{a(u)\sqrt{1+c^2}}\right|} e^{-\inv{2}\left(v(u)+\inv{a(u)\sqrt{1+c^2}}\right)^2}\\
\sim & \frac{e^{-\frac{b^2(u)}{2(1+c^2)}-\inv{2}\left(v(u)\right)^2}}{2\pi\sqrt{1+c^2}|a(u)|\left|v(u)+\inv{a(u)\sqrt{1+c^2}}\right|}, 
\end{align*}
by using (\ref{normal cdf expansion}). As $1/a(u)\sqrt{1+c^2}\to 0$ as $u\to 0^+$, so both lower and upper bounds of (\ref{Theorem: evaluate general integral, integral1}) (i.e. the two expression in (\ref{Theorem: evaluate general integral, integral1bounds})) coincide asymptotically in this case and we have 
\begin{equation*}
\int^{a(u)}_{-\infty}\Phi(x) \inv{\sqrt{2\pi}}e^{-\inv{2}(cx +b(u))^2}\,dx \sim \frac{e^{-\frac{b^2(u)}{2(1+c^2)}-\inv{2}\left(v(u)\right)^2}}{2\pi\sqrt{1+c^2}|a(u)v(u)|} 
\end{equation*}
if $v(u)\to -\infty$ as $u\to 0^+$. 

Next, if $v(u) \to 0$ then the upper bound of (\ref{Theorem: evaluate general integral, integral1}) i.e. (the second expression in (\ref{Theorem: evaluate general integral, integral1bounds})) is asymptotically equivalent to 
\begin{align*}
& \frac{e^{-\frac{b^2(u)}{2(1+c^2)}}}{\sqrt{2\pi}\sqrt{1+c^2}|a(u)|}\Phi\left(v(u)\right)  \sim \frac{e^{-\frac{b^2(u)}{2(1+c^2)}}}{\sqrt{2\pi}\sqrt{1+c^2}|a(u)|} \times \inv{2}.
\end{align*}
On the other hand, the lower bound of (\ref{Theorem: evaluate general integral, integral1}) (i.e. the first expression in (\ref{Theorem: evaluate general integral, integral1bounds})), when $v(u) \to 0$ as $u\to 0^+$, behaves asymptotically as 
\begin{align*}
& \frac{e^{-\frac{b^2(u)}{2(1+c^2)}+\inv{2}\left(v(u)+\inv{a(u)\sqrt{1+c^2}}\right)^2-\inv{2}\left(v(u)\right)^2}}{\sqrt{2\pi}\sqrt{1+c^2}|a(u)|}\Phi\left(v(u)+\inv{a(u)\sqrt{1+c^2}}\right)\\
\sim & \frac{e^{-\frac{b^2(u)}{2(1+c^2)}+\inv{2}\left(\inv{a(u)\sqrt{1+c^2}}\right)^2}}{\sqrt{2\pi}\sqrt{1+c^2}|a(u)|}\Phi\left(\inv{a(u)\sqrt{1+c^2}}\right)\\
\sim & \frac{e^{-\frac{b^2(u)}{2(1+c^2)}}}{\sqrt{2\pi}\sqrt{1+c^2}|a(u)|} \times \inv{2},
\end{align*}
as $a(u)\to -\infty$ when $u\to 0^+$, so both lower and upper bounds of (\ref{Theorem: evaluate general integral, integral1}) (i.e. the two expressions in 
(\ref{Theorem: evaluate general integral, integral1bounds})) again coincide asymptotically in this case and we have 
\begin{equation*}
\int^{a(u)}_{-\infty}\Phi(x) \inv{\sqrt{2\pi}}e^{-\inv{2}(cx +b(u))^2}\,dx \sim \frac{e^{-\frac{b^2(u)}{2(1+c^2)}}}{2\sqrt{2\pi}\sqrt{1+c^2}|a(u)|},
\end{equation*}
if $v(u) \to 0$.

Finally, if $v(u) \to \infty$ as $u\to 0^+$, it follows that 
\begin{equation}
1+c^2+\frac{cb(u)}{a(u)} \to 1+c^2+ck \leq 0. \label{Theorem: evaluate general integral, polynomial} 
\end{equation}
%(In the situation $a(u)(1+c^2)+cb(u) \to -\infty$, the same argument yields $1+c^2+ck >0$.) 
In this case, we put $w = \frac{z}{a(u)}$ in $J(u)$, which is defined in (\ref{Theorem: evaluate general integral, J term}), to get:
%Then from 
%Look at the ratio between (\ref{Lemma2:initial}) and $\frac{e^{-\frac{b^2(u)}{2(1+c^2)}}}{\sqrt{2\pi}|a(u)|\sqrt{1+c^2}}$ as $u\to 0^+$,  then follows from (\ref{Lemma2:case3}) we have 
\begin{eqnarray}
%\notag && \frac{|a(u)|}{\sqrt{2\pi}}\int^0_{-\infty} \inv{\left|\frac{z}{\sqrt{1+c^2}}+a(u)\right|} e^{-\inv{2}\left(z+\frac{a(u)(1+c^2)+cb(u)}{\sqrt{1+c^2}}\right)^2}\,dz \\
\frac{1}{\sqrt{2\pi}}J(u) &=& \frac{|a(u)|}{\sqrt{2\pi}} \int^{\infty}_{0} \inv{\left(1+\frac{w}{\sqrt{1+c^2}}\right)}e^{-\frac{a^2(u)}{2}\left(w+ \frac{1+c^2+\frac{cb(u)}{a(u)}}{\sqrt{1+c^2}}\right)^2}\,dw. \label{Lemma2:case3}
%\quad \text{putting $w= \frac{z}{a(u)}$}. 
\end{eqnarray}
Now suppose $g(w) = \inv{1+\frac{|w|}{\sqrt{1+c^2}}}, \, w \geq 0; = 0, w <0,$  and consider $Eg(W_u)$ where 
$W_u \sim N\left(-\frac{1+c^2+\frac{cb(u)}{a(u)}}{\sqrt{1+c^2}}, \inv{a^2(u)}\right)$. This means as $u\to 0^+$, 
\begin{equation}
W_u \overset{\mathcal{P}}{\to} -\frac{1+c^2+ck}{\sqrt{1+c^2}}\quad \text{and}\quad Eg(W_u) \to \inv{1+\frac{|1+c^2+ck|}{1+c^2}} = \inv{\frac{1+c^2-(1+c^2+ck)}{1+c^2}}  = -\frac{1+c^2}{ck} \quad (>0),
\end{equation}
where $\overset{\mathcal{P}}{\to}$ denotes the convergence in probability and $1+c^2+ck \leq 0 \Rightarrow ck \leq 0$ by (\ref{Theorem: evaluate general integral, polynomial}). 
%Now the difference between $Eg(W_u)$ and $\frac{1}{\sqrt{2\pi}}J(u)$ is 
%\begin{eqnarray*}
%&&  \frac{|a(u)|}{\sqrt{2\pi}} \int^{0}_{-\infty} \inv{\left(1+\frac{|w|}{\sqrt{1+c^2}}\right)}e^{-\frac{a^2(u)}{2}\left(w+ \frac{1+c^2+\frac{cb(u)}{a(u)}}{\sqrt{1+c^2}}\right)^2}\,dw\\
%&\leq& \frac{|a(u)|}{\sqrt{2\pi}} \int^{0}_{-\infty} e^{-\frac{a^2(u)}{2}\left(w+ \frac{1+c^2+\frac{cb(u)}{a(u)}}{\sqrt{1+c^2}}\right)^2}\,dw \\
%&=& P(W_u <0) \to 0, 
%\end{eqnarray*}
%as $W_u \overset{\mathcal{P}}{\to} -\frac{1+c^2+ck}{\sqrt{1+c^2}} >0$. 
As a result, it follows from (\ref{Theorem: evaluate general integral, integral1}) that
\begin{eqnarray*}
&&  \frac{\int^{a(u)}_{-\infty} \Phi(x) \inv{\sqrt{2\pi}} e^{-\inv{2}(cx + b(u))^2}\,dx}{\frac{e^{-\frac{b^2(u)}{2(1+c^2)}}}{\sqrt{2\pi}|a(u)|\sqrt{1+c^2}}}\sim \frac{1}{\sqrt{2\pi}}J(u) \to \lim_{u\to 0^+} Eg(W_u)  = -\frac{1+c^2}{ck} = \frac{1+c^2}{|ck|}.
\end{eqnarray*}
This completes the proof of Theorem \ref{Theorem: evaluate general integral}.

\end{proof}

Notice that $v(u) \to 0$ implies  $0 =(1+c^2 + ck),$ so that  $-ck = (1+c^2) = |ck|, $ so  an alternate asymptotic representation is possible in this case.  The expression in the statement of Theorem 1 is made to resemble that for when  $v(u) \to \infty.$

\section{Proof of Lemmas \ref{Lemma: exponents equivalency} to \ref{Lemma: beta1>0}} \label{Appendix: Proof of Lemmas}

\subsection{Proof of Lemma \ref{Lemma: exponents equivalency}}

\begin{proof}
As $B(u)$ can be expressed as $\alpha_iF_i^{-1}(u) + \alpha_{3-i}F_{3-i}^{-1}(u)$, combining $\alpha_iA_i(u)$ and $B(u)$ to solve for $F_{3-i}^{-1}(u)$, we obtain
\begin{equation*}
F_{3-i}^{-1}(u) = \frac{B(u)-\alpha_iA_i(u)}{\alpha_{3-i}+\rho\alpha_i} \quad \Rightarrow\quad \lambda_{3-i}^2(F_{3-i}^{-1}(u))^2 = \frac{(B(u)-\alpha_iA_i(u))^2}{1+\alpha_{i}^2(1-\rho^2)}, 
\end{equation*}
so the LHS of (\ref{Lemma: eqn exponents equivalency}) is 
\begin{align*}
& \inv{1+\alpha_i^2(1-\rho^2)} \left\{ (B(u)-\alpha_iA_i(u))^2+(1-\rho^2)\left(\frac{A_i(u)}{1-\rho^2}+\alpha_iB(u)\right)^2\right\}
\end{align*}
in view of (\ref{beta12u definition 3});{\adb
\begin{align*}
= & \inv{1+\alpha_i^2(1-\rho^2)} \Biggl\{ B^2(u)+ \alpha_i^2A_i^2(u)-2\alpha_iA_i(u)B(u) \\
& \quad + (1-\rho^2)\left[ \frac{A^2_i(u)}{(1-\rho^2)^2}+\alpha_i^2B^2(u)+\frac{2\alpha_iA_i(u)B(u)}{(1-\rho^2)}\right]\Biggr\}\\
=& \inv{1+\alpha_i^2(1-\rho^2)} \Biggl\{ B^2(u)\left(1+\alpha_i(1-\rho^2)\right)+ A_i^2(u)\left(\alpha_i^2+\inv{1-\rho^2}\right)\Biggr\}\\
= & \frac{A^2_i(u)}{1-\rho^2}+B^2(u) = \text{ RHS of (\ref{Lemma: eqn exponents equivalency})}.
\end{align*}}
\end{proof}

\subsection{Proof of Lemma \ref{Lemma: beta limit}}

\begin{proof}
From (\ref{detail F ratio}) in Lemma \ref{Lemma: detailed F inverse}, we have 
\begin{equation}
\left(\frac{F_i^{-1}(u)}{F_{3-i}^{-1}(u)} - \gamma_i\right) \times F_{3-i}^{-1}(u) = O\left(\frac{\log(-\log u)}{\sqrt{-
\log u}}\right) \to 0 \label{F ratio minus dominating term goes to 0}
\end{equation}
as $u\to 0^+$. Consequently from the expression (\ref{beta12u definition 2}) for $\beta_{i}(u)$ we have 
\begin{align*}
& \left(\beta_{i}(u)-\beta_{i}\right)\times F_{3-i}^{-1}(u) \\
=& \left\{ \frac{F_i^{-1}(u)}{F_{3-i}^{-1}(u)}\left(\inv{1-\rho^2}+\alpha_i^2\right)-\frac{\rho}{1-\rho^2}+\alpha_1\alpha_2 -\beta_{i}\right\} \times F_{3-i}^{-1}(u) \\
=& \left(\frac{F_i^{-1}(u)}{F_{3-i}^{-1}(u)} -\gamma_i \right)\left(\inv{1-\rho^2}+\alpha_i^2\right) \times F_{3-i}^{-1}(u) 
\end{align*}
from (\ref{beta12 definition}), so from (\ref{F ratio minus dominating term goes to 0}) the Lemma is proved. 
\end{proof}

\subsection{Proof of Lemma \ref{Lemma: beta1>0}}

\begin{proof} Suppose $\lambda_2\geq 0$. From (\ref{beta12 definition}), we can see that $\beta_1$ is defined as 
\begin{equation*}
\beta_1 = \gamma_1\left( \frac{1+\alpha_1^2(1-\rho^2)}{1-\rho^2} \right)+\frac{\alpha_1\alpha_2(1-\rho^2)-\rho}{1-\rho^2}.
\end{equation*}
Now when $\lambda_2\geq 0$, from (\ref{F ratio}) we have $\gamma_1 \geq \sqrt{\frac{1+\lambda_2^2}{1+\lambda_1^2}}$, and
{\adb
\begin{align*}
\beta_1 & \geq \sqrt{\frac{1+\lambda_{2}^2}{1+\lambda_1^2}}\left(\frac{1+\alpha_1^2(1-\rho^2)}{1-\rho^2}\right) + \frac{\alpha_1\alpha_2(1-\rho^2)-\rho}{1-\rho^2} \\
&= \sqrt{\frac{1+\alpha_{2}^2(1-\rho^2)}{1+\alpha_1^2(1-\rho^2)}}\left(\frac{1+\alpha_1^2(1-\rho^2)}{1-\rho^2}\right)+ \frac{\alpha_1\alpha_2(1-\rho^2)-\rho}{1-\rho^2} \\
&= \frac{\sqrt{(1+\alpha_{2}^2(1-\rho^2))(1+\alpha_1^2(1-\rho^2))}+\alpha_1\alpha_2(1-\rho^2)-\rho}{1-\rho^2}.
\end{align*}
} To prove $\beta_1>0$ when $\lambda_2\geq 0$, consider the expression:
\begin{eqnarray*}& & \left(\sqrt{(1+\alpha_2^2(1-\rho^2))(1+\alpha_1^2(1-\rho^2))}+\alpha_1\alpha_2(1-\rho^2)-\rho\right) \\
 &\times& \left (\sqrt{(1+\alpha_2^2(1-\rho^2))(1+\alpha_1^2(1-\rho^2))}- \alpha_1\alpha_2(1-\rho^2)+\rho \right )\\
&=& (1+\alpha_2^2(1-\rho^2))(1+\alpha_1^2(1-\rho^2)) - (\alpha_1\alpha_2(1-\rho^2)- \rho)^2\\
&=&(1-\rho^2)(1+ \alpha_1^2 +\alpha_2^2 + 2 \alpha_1 \alpha_2 \rho)\\
&\geq&(1-\rho^2)(1+ \alpha_1^2 +\alpha_2^2  - 2 |\alpha_1|| \alpha_2 |)\\
&=& (1-\rho^2)( 1 + (|\alpha_1| - |\alpha_2|)^2) > 0,
\end{eqnarray*} so both factors in the expression must be positive.

Now suppose $B(u) \to -\infty$. If also $\lambda_2 \geq 0$, we have  from the above that $\beta_1 >0.$  

So now suppose that  $B(u) \to -\infty$ and also $\lambda_2<0$ i.e. $\alpha_2+\rho\alpha_1<0$.
Then $\lambda_1>0$, since  assuming the contrary that  $\lambda_1\leq 0$, from (\ref{F ratio}), $\lambda_1\leq 0$ and $\lambda_2<0$ imply that $\gamma_1=1$, so 
\begin{equation*}
\frac{B(u)}{F_2^{-1}(u)} \to  \alpha_1\gamma_1+\alpha_2 = \alpha_1+\alpha_2 >0, \quad \text{as $u\to 0^+$}.
\end{equation*}
Furthermore, $\lambda_1\leq 0$, $\lambda_2<0$ are equivalent to  {\adb
\begin{align*}
 \alpha_1+\rho\alpha_2\leq 0 \quad \text{and}\quad  \alpha_2+\rho\alpha_1<0,
\end{align*}
}and  by adding these we have $(\alpha_1+\alpha_2)(1+\rho)<0 \Rightarrow \alpha_1+\alpha_2<0$ as $1+\rho>0$ which is a contradiction to the assumption that $B(u) \to -\infty$. This means that $\lambda_1>0$ when $B(u)\to -\infty$ and $\lambda_2<0$. 
%At the moment, there is no restriction on $\lambda_1$, so let's consider $\lambda_1\leq 0$ and $\lambda_1>0$ separately.
%
%Assume $\lambda_1\leq 0$ (and $\lambda_2<0$). Then 
%This means that case (a) cannot occur under assumption (1), even through $\gamma_1-\rho>0$ when $\lambda_1\leq 0$, $\lambda_2<0$. (This is due to the fact that $\alpha_1\gamma_1+\alpha_2>0$ cannot occur).
%\item $\lambda_1>0$ (and $\lambda_2<0$), so 
%\begin{align}
%& \lambda_1>0 \iff \alpha_1+\rho\alpha_2>0 \\
%& \lambda_2<0 \iff \alpha_2+\rho\alpha_1<0 
%\end{align}
This also implies that $\gamma_1-\rho>0$ from Lemma \ref{Lemma: at least one condition is positive}, where it is shown that $ 0<\gamma_1 - \rho$ except possibly in the case $\lambda_1 >0, \lambda_2 >0.$ 
%And we are additionally assuming from (1) that 
%\begin{equation}
%\alpha_1\gamma_1+\alpha_2>0.
%\end{equation}
Using $\lambda_1>0$ and $-\lambda_2>0$, we have
{\adb
\begin{align*}
\alpha_1+\rho\alpha_2-(\alpha_2+\rho\alpha_1)>0 \quad \Rightarrow & \quad \alpha_1(1-\rho)+\alpha_2(\rho-1)>0\\
\Rightarrow & \quad (\alpha_1-\alpha_2)(1-\rho)>0 \Rightarrow \, \,   \alpha_1>\alpha_2.
\end{align*}}
If $\alpha_1<0$, then $\alpha_2<0$ as well. This means that $\alpha_1\gamma_1+\alpha_2<0$ as $\gamma_1>0$ always, which is a contradiction to $B(u) \to -\infty$. As a result, $\alpha_1$ must be positive, since by global assumption $\alpha_1 \neq 0.$ So now we have  $\alpha_1> 0$, $\gamma_1-\rho>0$ and $\alpha_1\gamma_1+\alpha_2>0$, so  $\beta_1>0$ from (\ref{beta12 definition}) 

Finally consider $B(u) \to 0$, so $\alpha_1 \gamma_1 + \alpha_2 =0 $. Then  from (\ref{beta12 definition}) we have $\beta_1= \frac{\gamma_1 - \rho}{1 - \rho^2}.$ But from Lemma 2 whenever  $\gamma_1 - \rho \leq 0$ we  have  $\alpha_1 \gamma_1 + \alpha_2 > 0$, a contradiction, so $\beta_1 >0.$

%Now, if $\alpha_1\geq 0$, from (1) and (2), we have $\beta_1>0$. 
%
%On the other hand, if $\alpha_1<0$ from (5), we have $\alpha_2<0$. Then since $\gamma_1>0$, $\alpha_1\gamma_1+\alpha_2<0$ is contradicting to (5). So $\alpha_2<0$ cannot occur. (This is due to the fact that when $\lambda_1>0$ and $\lambda_2<0$ and $\alpha_2<0$, $\alpha_1\gamma_1+\alpha_2>0$ cannot occur). (From the proof of Lemma 2, we know that $\lambda_1>0$, $\lambda_2<0$, $0<\rho<1$ that $\alpha_1>0$.)
%
%Thus if $\lambda_1<0$, $\lambda<0$ and (1) holds, $\alpha_1\geq 0$ and so $\beta_1>0$.
%So under (1), $\beta_1>0$.

%Overall, 
%\begin{equation}
%\beta_{1} = \lim_{u\to 0^+} \beta_{1}(u) >0. \label{beta12>0}
%\end{equation}
%Notice that $\frac{B(u)}{F_1^{-1}(u)} \to \alpha_2\gamma_2+\alpha_1 = \frac{\alpha_2}{\gamma_1} +\alpha_1 = \inv{\gamma_1}\left( \alpha_1\gamma_1+\alpha_2 \right)$.  

By ``symmetry'', when $\lambda_1\geq 0$ or $B(u) \to -\infty$ or  $B(u) \to -\infty$ or $B(u) \to 0$  , we have $\beta_2>0$.
\end{proof}

\section{Proof of Theorem \ref{Thm: main result}} \label{Appendix: Proof of Theorem 2}

\begin{proof}

\begin{enumerate}[label=(\alph*)]
 \item Assuming $A(u) \to -\infty$ and $B(u) \to -\infty$,  the behaviour of the two terms in (\ref{Theorem1:intergral A}) was stated in Corollaries \ref{Corollary: specific integral with rho in the condition, first term} and \ref{Corollary: specific integral with rho in the condition, second term} respectively.  Now we have to compute the difference to see which of the two terms will be dominating. 
% By comparing the structure of the terms in (\ref{eqn:comparing 2 terms}), we can see that the behaviour of $\alpha_1F_1^{-1}(u)+\alpha_1F_2^{-1}(u)$ will decide which of the two terms will be dominating, so let's consider them separately.  

As $B(u) \to -\infty$ as $u\to 0^+$,  we have $\beta_{1}>0$ from Lemma \ref{Lemma: beta1>0} and  (\ref{Theorem1:intergral A}) is  equivalent to {\adb
 \begin{align}
\notag &  \frac{\sqrt{1-\rho^2}e^{-\frac{A_1^2(u)}{2(1-\rho^2)}-\inv{2}B^2(u)}}{2\pi|A_1(u)||B(u)|\Phi(\lambda_2F_2^{-1}(u))}(1+o_1(1))\\
\notag & \quad   - \frac{\alpha_1\sqrt{1-\rho^2}e^{-\frac{A_1^2(u)}{2(1-\rho^2)} - \inv{2}B^2(u)}}{2\pi|F_2^{-1}(u)| \Phi(\lambda_2F_2^{-1}(u)) |A_1(u)|\beta_{1}(u)}(1+o_2(1))\\
\notag =& \frac{\sqrt{1-\rho^2}e^{-\frac{A_1^2(u)}{2(1-\rho^2)}-\inv{2}B^2(u)}}{2\pi|A_1(u)|\Phi(\lambda_2F_2^{-1}(u))}\left\{
\left(\inv{|B(u)|}-\frac{\alpha_1}{\beta_{1}(u)|F_2^{-1}(u)|}\right)+ \frac{o_1(1)}{|B(u)|} + \frac{\alpha_1o_2(1)}{\beta_{1}(u)|F_2^{-1}(u)|}\right\}\\
\notag = &
\frac{\sqrt{1-\rho^2}e^{-\frac{A_1^2(u)}{2(1-\rho^2)}-\inv{2}B^2(u)}}{2\pi|A_1(u)|\Phi(\lambda_2F_2^{-1}(u))}
\Biggl\{
\frac{|A_1(u)|}{(1-\rho^2)|B(u)|\beta_{1}|F_2^{-1}(u)|}(1+o_3(1)) \\
& \quad +  \frac{o_1(1)}{|B(u)|} + \frac{\alpha_1o_2(1)}{\beta_{1}(u)|F_2^{-1}(u)|}\Biggr\}, 
\label{Theorem 2: part I i considering the two ratios}
%\sim& \frac{e^{-\frac{A_1^2(u)}{2(1-\rho^2)}-\inv{2}B^2(u)}}{2\pi\sqrt{1-\rho^2}|B(u)|\Phi(\lambda_2F_2^{-1}(u))\beta_{1.2}|F_2^{-1}(u)|} 
  \end{align}}
 as{\adb
 \begin{align*}
& \inv{|B(u)|}- \frac{\alpha_1}{\beta_{1}(u)|F_2^{-1}(u)|} \\
= & \frac{\beta_{1}(u)|F_2^{-1}(u)|-\alpha_1|B(u)|}{|B(u)|\beta_{1}(u)|F_2^{-1}(u)|}\\
=& \frac{- \left(\frac{A(u)}{1-\rho^2}+\alpha_1B(u)\right)+\alpha_1B(u)}{|B(u)|\beta_{1}(u)|F_2^{-1}(u)|}, \quad \text{by (\ref{beta12u definition 3})}\\
%= & \frac{|F_1^{-1}(u)-\rho F_2^{-1}(u)|}{(1-\rho^2)|\alpha_1F_1^{-1}(u)+\alpha_2F_2^{-1}(u)|\beta_{1}(u) |F_2^{-1}(u)|}\\
=& \frac{|A_1(u)|}{(1-\rho^2)|B(u)|\beta_{1}(u)|F_2^{-1}(u)|}\\
=& \frac{|A_1(u)|}{(1-\rho^2)|B(u)|\beta_{1} |F_2^{-1}(u)|}(1+o_3(u)). 
 \end{align*}}
Since both $\alpha_1\gamma_1+\alpha_2>0$ and $\gamma_1-\rho>0$ by Corollary \ref{Corollary: at least one} (since both $A_1(u)$ and $B(u)\to -\infty$ by assumption),  (\ref{Theorem 2: part I i considering the two ratios}) becomes
\begin{align*}
& \frac{\sqrt{1-\rho^2}e^{-\frac{A_1^2(u)}{2(1-\rho^2)}-\inv{2}B^2(u)}}{2\pi|A_1(u)|\Phi(\lambda_2F_2^{-1}(u))}\left\{ \frac{|A_1(u)|}{(1-\rho^2)|B(u)|\beta_{1}|F_2^{-1}(u)|}\left(1+o_3(1)+o_4(1)+o_5(1)\right)\right\}\\
%\sim & \frac{e^{-\frac{A_1^2(u)}{2(1-\rho^2)}-\inv{2}B^2(u)}}{2\pi\sqrt{1-\rho^2}|B(u)|\Phi(\lambda_2F_2^{-1}(u))\beta_{1.2}|F_2^{-1}(u)|}
\sim& \frac{e^{-\frac{A_1^2(u)}{2(1-\rho^2)}-\inv{2}B^2(u)}}{2\pi\sqrt{1-\rho^2}|B(u)|\Phi(\lambda_2F_2^{-1}(u))\beta_{1}|F_2^{-1}(u)|}.
\end{align*}

Now assuming $B(u) \to -\infty$ but $A_1(u) \not\to -\infty$ i.e. $A_1(u) \to 0$ or $\infty$,  from Lemma \ref{Lemma: beta1>0}, we have $\beta_{1}>0$, and from Lemma 1,  $\alpha_1 >0.$ Finally, from Corollary \ref{Corollary: specific integral without rho in the condition}:
\begin{align*}
& P(Z_1\leq F_1^{-1}(u) | Z_2 = F_2^{-1}(u)) \sim  \frac{e^{-\frac{A_1^2(u)}{2(1-\rho^2)} - \inv{2}B^2(u)}}{2\pi \sqrt{1-\rho^2}\Phi(\lambda_2F_2^{-1}(u))|B(u)| |F_2^{-1}(u)|\beta_{1}}.
\end{align*}

\item Assuming $B(u) \to 0$, so  $A_1(u) \to -\infty$,  we can see from Lemma \ref{Lemma: beta1>0} that $\beta_1>0$. As $A_1(u)\to -\infty$, the behaviour of $P(Z_1\leq F_1^{-1}(u)| Z_2= F_2^{-1}(u))$ can be expressed as (\ref{Theorem1:intergral A}) and the treatment of those two terms was explained in Corollaries \ref{Corollary: specific integral with rho in the condition, first term} and \ref{Corollary: specific integral with rho in the condition, second term} respectively. This means that 
{\adb
 \begin{align*}
&  P(Z_1\leq F_1^{-1}(u)| Z_2= F_2^{-1}(u)) \\
 \sim  &    \frac{\sqrt{1-\rho^2}e^{-\frac{A_1^2(u)}{2(1-\rho^2)}}
 }{2\sqrt{2\pi}|A_1(u)|\Phi(\lambda_2F_2^{-1}(u))}(1+o_1(1))\\
&  \quad - \frac{\alpha_1\sqrt{1-\rho^2}e^{-\inv{2}\frac{A_1^2(u)}{1-\rho^2} - \inv{2}B^2(u)}}{2\pi|F_2^{-1}(u)| \Phi(\lambda_2F_2^{-1}(u)) |A_1(u)|\beta_{1}}(1+o_2(1))\\
 \sim &  \frac{\sqrt{1-\rho^2}e^{-\frac{A_1^2(u)}{2(1-\rho^2)}}
 }{2\sqrt{2\pi}|A_1(u)|\Phi(\lambda_2F_2^{-1}(u))},
 \end{align*}}
 as the first term is the dominating one. 
\item Assuming $B(u) \to \infty$, so $A_1(u) \to -\infty$, we can subdivide the proof based on the behaviour of $\beta_1$.
\begin{enumerate}[label=(\roman*)]
\item If $\beta_1>0$, then (\ref{Theorem1:intergral A}) is  equivalent to 
\begin{align*}
&  P(Z_1\leq F_1^{-1}(u)| Z_2= F_2^{-1}(u)) \\
\sim   &   \frac{\sqrt{1-\rho^2}e^{-\frac{A_1^2(u)}{2(1-\rho^2)}}
 }{\sqrt{2\pi}|A_1(u)|\Phi(\lambda_2F_2^{-1}(u))}(1+o_1(1)) - \frac{\alpha_1\sqrt{1-\rho^2}e^{-\inv{2}\frac{A_1^2(u)}{1-\rho^2} - \inv{2}B^2(u)}}{2\pi|F_2^{-1}(u)| \Phi(\lambda_2F_2^{-1}(u)) |A_1(u)|\beta_{1}}(1 +o_2(1))\\
 \sim &  \frac{\sqrt{1-\rho^2}e^{-\frac{A_1^2(u)}{2(1-\rho^2)}}
 }{\sqrt{2\pi}|A_1(u)|\Phi(\lambda_2F_2^{-1}(u))}
\end{align*} 
 as the first term is again the dominating one.
\item If $\beta_1=0$, then (\ref{Theorem1:intergral A}) says  {\adb
\begin{align*}
&  P(Z_1\leq F_1^{-1}(u)| Z_2= F_2^{-1}(u)) \\
= &  \frac{\sqrt{1-\rho^2}e^{-\frac{A_1^2(u)}{2(1-\rho^2)}}}{\sqrt{2\pi}|A_1(u)|\Phi(\lambda_2F_2^{-1}(u))}\left( 1+o_1(1)) \right) \\
& \quad - \frac{\alpha_1(1-\rho^2)e^{-\inv{2}\lambda_2^2(F_2^{-1}(u))^2}}{2\sqrt{2\pi}\Phi(\lambda_2F_2^{-1}(u))|A_1(u)|\sqrt{1+\alpha_1^2(1-\rho^2)}}(1+o_2(1))\\ 
\sim & \frac{\sqrt{1-\rho^2}e^{-\frac{A_1^2(u)}{2(1-\rho^2)}}}{\sqrt{2\pi}|A_1(u)|\Phi(\lambda_2F_2^{-1}(u))}\left[1 - \frac{\alpha_1\sqrt{1-\rho^2}e^{-\inv{2}(\lambda_2F_2^{-1}(u))^2+\frac{A_1^2(u)}{2(1-\rho^2)}}}{2\sqrt{1+\alpha_1^2(1-\rho^2)}}\right] \\
= & \frac{\sqrt{1-\rho^2}e^{-\frac{A_1^2(u)}{2(1-\rho^2)}}}{\sqrt{2\pi}|A_1(u)|\Phi(\lambda_2F_2^{-1}(u))}
\left[ 1 - \frac{\alpha_1\sqrt{1-\rho^2} e^{-\inv{2}B^2(u)+ \inv{2}\left(\frac{1-\rho^2}{1+\alpha_1^2(1-\rho^2)}\right)\left(\beta_{1}(u)F_2^{-1}(u)\right)^2}}{2\sqrt{1+\alpha_1^2(1-\rho^2)}}\right], \, \, \text{by (\ref{Lemma: eqn exponents equivalency});}\\
  \sim & \frac{\sqrt{1-\rho^2}e^{-\frac{A_1^2(u)}{2(1-\rho^2)}}}{\sqrt{2\pi}|A_1(u)|\Phi(\lambda_2F_2^{-1}(u))}
\end{align*}}
as required since $B(u) \to \infty$ and $\beta_{1}(u)F_2^{-1}(u)\to 0$ from (\ref{eqn beta12uF2inv goes to 0}), as $u\to 0^+$. 
\item  If $\beta_1<0$ and  $B(u)\to \infty$ as well, we have,  using Corollaries \ref{Corollary: specific integral with rho in the condition, first term} and \ref{Corollary: specific integral with rho in the condition, second term}, that  (\ref{Theorem1:intergral A})  says
\begin{align*}
&  P(Z_1\leq F_1^{-1}(u)| Z_2= F_2^{-1}(u)) \\
= & \frac{\sqrt{1-\rho^2}e^{-\frac{A_1^2(u)}{2(1-\rho^2)}}}{\sqrt{2\pi}|A_1(u)|\Phi(\lambda_2F_2^{-1}(u))}(1+o_1(1)) - \frac{\alpha_1 e^{-\frac{\lambda_2^2(F_2^{-1}(u))^2}{2}}}{\sqrt{2\pi}|\lambda_2F_2^{-1}(u)| \Phi(\lambda_2F_2^{-1}(u))}(1+o_2(1))\\ 
 \sim & \frac{\sqrt{1-\rho^2}e^{-\frac{A_1^2(u)}{2(1-\rho^2)}}}{\sqrt{2\pi}|A_1(u)|\Phi(\lambda_2F_2^{-1}(u))},
\end{align*}
as from (\ref{special})  we have 
\begin{equation*}
 -\frac{A_1^2(u)}{1-\rho^2}+ \lambda_2^2(F_2^{-1}(u))^2 \to \infty 
\end{equation*}
since $\beta_{1}<0 \Rightarrow \beta_{1}(u)F_2^{-1}(u) \to \infty$, $B(u)\to \infty$ and so $A_1(u)\to -\infty$, so the first term will dominate here. 
\end{enumerate}
This means that the final expression in (c) is the same irrespective of the value  of $\beta_1$ and  
\begin{equation*}
 P(Z_1\leq F_1^{-1}(u)| Z_2= F_2^{-1}(u)) \sim   \frac{\sqrt{1-\rho^2}e^{-\frac{A_1^2(u)}{2(1-\rho^2)}}
 }{\sqrt{2\pi}|A_1(u)|\Phi(\lambda_2F_2^{-1}(u))}.
\end{equation*}
 \end{enumerate}

This completes the proof of Theorem \ref{Thm: main result}.
\end{proof}

\section{Proof of Lemmas \ref{lemma: general form of G1F1-GF2} to \ref{Lemma: exp B in RV}} \label{Appendix: Proof of Lemma: general form of G1F1-GF2}

\subsection{Proof of Lemmas \ref{lemma: general form of G1F1-GF2}}
\begin{proof}
From Lemma \ref{Lemma: detailed F inverse}, we have 
\begin{align*}
F_2^{-1}(u) = K_{2,1}(-\sqrt{-2\log u})\left\{ 1+ \frac{K_{2,2} \log(-\log u)}{\log u} + \frac{K_{2,3}}{\log u}+O\left(\left(\frac{\log(-\log u)}{\log u}\right)^2\right)\right\} % \label{generic F2inverse}
\end{align*}
 as $u\to 0^+$, where $K_{2,1}>0$, so
\begin{align*}
(F_2^{-1}(u))^2 = K_{2,1}^2(-2\log u)\left\{ 1+\frac{2K_{2,2}\log(-\log u)}{\log u} + \frac{2K_{2,3}}{\log u} + O\left(\left(\frac{\log(-\log u)}{\log u}\right)^2\right)\right\}. 
% \label{generic F2inverse sq}
\end{align*}
%(See Lemma 1 for specific $K_{2,2}$, $K_{2,3}$, $K_{2,4}$ in 3 cases of $\lambda_2=0$, $\lambda_2>0$ \& $\lambda_2<0$.)
Next, once again from Lemma \ref{Lemma: detailed F inverse} that, 
\begin{align*}
\frac{F_1^{-1}(u)}{F_2^{-1}(u)} - \gamma_1 = \gamma_1\left\{
\frac{C_{1,1}\log(-\log u)}{\log u} +\frac{C_{1,2}}{\log u} + O\left(\left(\frac{\log(-\log u)}{\log u}\right)^2\right)\right\} 
%\label{generic ratio of Finverses}
\end{align*}
where $C_{1,1}$ may be zero; so 
\begin{align*}
\left(\frac{F_1^{-1}(u)}{F_2^{-1}(u)}-\gamma_1\right)^2 = O\left(\left(\frac{\log(-\log u)}{\log u}\right)^2\right). 
%\label{generic ratio of Finverses sq}
\end{align*}
Since 
\begin{equation*}
G_1\frac{F_1^{-1}(u)}{F_2^{-1}(u)} - G_2 = G_1 \underbrace{\left\{\frac{F_1^{-1}(u)}{F_2^{-1}(u)}-\gamma_1
\right\}}_{\downarrow 0} + \underbrace{G_1\gamma_1-G_2}_{\text{constant}},
\end{equation*}
we have 
\begin{align*}
& \left(\frac{G_1F_1^{-1}(u)}{F_2^{-1}(u)} - G_2\right)^2 \\
=& \left(G_1\left\{\frac{F_1^{-1}(u)}{F_2^{-1}(u)} -\gamma_1\right\}+G_1\gamma_1- G_2\right)^2 \\
=& \left(G_1\left\{\frac{F_1^{-1}(u)}{F_2^{-1}(u)}-\gamma_1\right\}\right)^2+2G_1(G_1\gamma_1-G_2)\left(\frac{F_1^{-1}(u)}{F_2^{-1}(u)}-\gamma_1\right)+(G_1\gamma_1-G_2)^2\\
=& \left(G_1\gamma_1-G_2\right)^2 + 2G_1\left(G_1\gamma_1-G_2\right)\gamma_1\left(\frac{C_{1,1}\log|\log u|}{\log u}+\frac{C_{1,2}}{\log u}\right)+O\left(\left(\frac{\log(-\log u)}{\log u}\right)^2\right).
\end{align*}
Consider now {\adb
\begin{align}
& \notag (G_1F_1^{-1}(u) - G_2F_2^{-1}(u))^2 \\
\notag =& (F_2^{-1}(u))^2\left(G_1\frac{F_1^{-1}(u)}{F_2^{-1}(u)}-G_2\right)^2\\
\notag =& K_{2,1}^2(-2\log u)\left\{1+\frac{2K_{2,2}\log(-\log u)}{\log u}+\frac{2K_{2,3}}{\log u}+O\left(\left(\frac{\log(-\log u)}{\log u}\right)^2\right)\right\}\\
\notag & \quad \times \left\{ (G_1\gamma_1-G_2)^2+2(G_1\gamma_1-G_2)\gamma_1\left(\frac{C_{1,1}\log(-\log u)}{\log u}+\frac{C_{1,2}}{\log u}\right)+O\left(\left(\frac{\log(-\log u)}{\log u}\right)^2\right)\right\}\\
\notag =&  K_{2,1}^2(-2\log u)\Biggl\{ (G_1\gamma_1-G_2)^2+\frac{2K_{2,2}(G_1\gamma_1-G_2)^2\log(-\log u)}{\log u}+\frac{2K_{2,3}(G_1\gamma_1-G_2)^2}{\log u}\\
\notag & \quad +2(G_1\gamma_1-G_2)\gamma_1\left(\frac{C_{1,1}\log(-\log u)}{\log u}+\frac{C_{1,2}}{\log u}\right)+O\left(\left(\frac{\log(-\log u)}{\log u}\right)^2\right)\Biggr\}\\
%\notag =& \left\{G_1^2K_{2,1}^2(\gamma_1-G_2)^2\right\}(-2\log u)\\
%\notag & \quad -4\log(-\log u)\left\{G_1^2K_{2,1}^2K_{2,2}(\gamma_1-G_2)^2+G_1^2K_{2,1}^2\gamma_1(\gamma_1-G_2)C_{1,1}\right\}\\
%\notag  & \quad -4\left\{ K_{2,1}^2K_{2,3}(G_1\gamma_1-G_1G_2)^2+G_1^2K_{2,1}^2C_{1,2}\gamma_1(\gamma_1-G_2)\right\} + O\left(\frac{\left[\log(-\log u)\right]^2}{\log u}\right)\\
 \notag =& \left\{K_{2,1}^2(G_1\gamma_1-G_2)^2\right\}(-2\log u)\\
\notag & \quad -4\log(-\log u)\left\{K_{2,1}^2K_{2,2}(G_1\gamma_1-G_2)^2+K_{2,1}^2G_1\gamma_1(G_1\gamma_1-G_2)C_{1,1}\right\}\\
 & \quad -4\left\{K_{2,1}^2K_{2,3}(G_1\gamma_1-G_2)^2+K_{2,1}^2C_{1,2}G_1\gamma_1(G_1\gamma_1-G_2)\right\} + O\left(\frac{\left[\log(-\log u)\right]^2}{\log u}\right)
  \label{generic denominator}
\end{align}}
Thus 
\begin{align}
\notag & e^{-\inv{2}\left(G_1F_1^{-1}(u)-G_2F_2^{-1}(u)\right)^2}\\
\notag = & e^{-\inv{2}(F_2^{-1}(u))^2\left(G_1\frac{F_1^{-1}(u)}{F_2^{-1}(u)}-G_2\right)^2}\\
\notag =& e^{ \left\{K_{2,1}^2(G_1\gamma_1-G_2)^2\right\}(\log u) +2\log(-\log u)\left\{K_{2,1}^2K_{2,2}(G_1\gamma_1-G_2)^2+K_{2,1}^2C_{1,1}G_1\gamma_1(G_1\gamma_1-G_2)\right\}}\\
\notag & \quad \times e^{2\left\{ K_{2,1}^2K_{2,3}(G_1\gamma_1-G_2)^2+K_{2,1}^2C_{1,2}G_1\gamma_1(G_1\gamma_1-G_2)\right\} + O\left(\frac{\left[\log(-\log u)\right]^2}{\log u}\right)}\\
=& \tau_1 u^{\theta}(-\log u)^{\tau_2}\left(1+O\left(\frac{\left[\log|\log u|\right]^2}{\log u}\right)\right), \label{generic nominator}
%u^{G_1^2K_{2,2}^2(\gamma_1-G_2)^2}(-\log u)^{\text{constant}_1}\times \text{constant}_2
\end{align}
where 
\begin{align*}
\theta =& K_{2,1}^2(G_1\gamma_1-G_2)^2;\\
\tau_1 = & e^{2\left[K_{2,1}^2K_{2,3}(G_1\gamma_1-G_2)^2+K_{2,1}^2C_{1,2}G_1\gamma_1(G_1\gamma_1-G_2)\right]};\\
\tau_2 =& 2\left[K_{2,1}^2K_{2,2}(G_1\gamma_1-G_2)^2+K_{2,1}^2C_{1,1}G_1\gamma_1(G_1\gamma_1-G_2)\right].
\end{align*}

From (\ref{generic denominator}), we know that if $\left|G_1 F_1^{-1}(u)-G_1G_2F_2^{-1}(u)\right| \to \infty$ as $u\to 0^+$, then 
\begin{align*}
 (F_2^{-1}(u))^2\left(G_1\frac{F_1^{-1}(u)}{F_2^{-1}(u)}-G_2\right)^2
%=& \left\{G_1^2K_2^2(\gamma(1,2)-G_2)^2\right\}(-2\log u)\\
%& \quad -2\log(-\log u)\left\{G_1^2K_2^2K_3(\gamma(1,2)-G_2)^2+2G_1^2K_2^2C_3(\gamma(1,2)-G_2)\right\}\\
% & \quad -2\left\{ 2G_1^2K_2^2K_4(\gamma(1,2)-G_2)^2+2G_1^2K_2^2C_4(\gamma(1,2)-G_2)\right\} + O\left(\left(\frac{\log(-\log u)}{\log u}\right)^2\right)\\
% =& \left\{G_1^2K_2^2(\gamma(1,2)-G_2)^2\right\}(-2\log u)\\
% &\quad \Biggl\{ 1 + \frac{\log(-\log u)\left\{G_1^2K_2^2K_3(\gamma(1,2)-G_2)^2+2G_1^2K_2^2C_3(\gamma(1,2)-G_2)\right\}}{-\log u}\\
% &\quad \frac{-2\left\{ 2G_1^2K_2^2K_4(\gamma(1,2)-G_2)^2+2G_1^2K_2^2C_4(\gamma(1,2)-G_2)\right\} + O\left(\left(\frac{\log(-\log u)}{\log u}\right)^2\right)}{-\log u}\Biggr\}\\
 &\sim  \left\{K_{2,1}^2(G_1\gamma_1- G_2)^2\right\}(-2\log u)\\
 \Rightarrow\quad \left|F_2^{-1}(u)\left(G_1\frac{F_1^{-1}(u)}{F_2^{-1}(u)}-G_2\right)\right| &\sim \left|K_{2,1}(G_1\gamma_1-G_2)\right| \sqrt{-2\log u}  
\end{align*}
as $u\to 0^+$. 
\end{proof}

\subsection{Proof of Lemma \ref{Lemma: Phi in RV}}

\begin{proof}

Using Lemma \ref{lemma: general form of G1F1-GF2} with $G_1=0$ and $G_2 = -\lambda_2$, we have 
\begin{align}
\notag &  e^{-\inv{2}\lambda_2^2(F_2^{-1}(u))^2} \\
\notag =& e^{\lambda_2^2K_{2,1}^2\log u\left\{ 1+\frac{2K_{2,2}\log(-\log u)}{\log u} + \frac{2K_{2,3}}{\log u} + O\left(\left(\frac{\log(-\log u)}{\log u}\right)^2\right)\right\}}\\
\notag = & e^{\lambda_2^2K_{2,1}^2\log u + 2\lambda_2^2K_{2,1}^2K_{2,2}\log(-\log u)+2\lambda^2_2K_{2,1}^2K_{2,3}+O\left(\frac{[\log(-\log u)]^2}{\log u}\right)}\\
\notag = & u^{\lambda_2^2K_{2,1}^2}|\log u|^{2\lambda_2^2K_{2,1}^2K_{2,2}}e^{2\lambda_2^2K_{2,1}^2K_{2,3}}\left(1+O\left(\frac{[\log(-\log u)]^2}{\log u}\right)\right)\\
= &
\begin{cases}
u^{\frac{\lambda_2^2}{1+\lambda_2^2}}|\log u|^{\frac{\lambda_2^2}{1+\lambda_2^2}}(2\pi\lambda_2)^{\frac{\lambda_2^2}{1+\lambda_2^2}}\left(1+O\left(\frac{[\log(-\log u)]^2}{\log u}\right)\right), & \text{if $\lambda_2>0$};\\
u^{\lambda_2^2}|\log u|^{\frac{\lambda_2^2}{2}}\pi^{\frac{\lambda_2^2}{2}}\left(1+O\left(\frac{[\log(-\log u)]^2}{\log u}\right)\right), 
&\text{if $\lambda_2<0$}.
\end{cases} \label{eqn: exp F sq in RV}
\end{align}

When $\lambda_2>0$, combining (\ref{normal cdf expansion}) with (\ref{eqn: exp F sq in RV}), we have {\adb
\begin{align*}  
|\lambda_2 F_2^{-1}(u)|\Phi\left(\lambda_2F_2^{-1}(u)\right) \sim & \inv{\sqrt{2\pi}} e^{-\inv{2}\lambda_2^2(F_2^{-1}(u))^2} 
% =&   \inv{\sqrt{2\pi}}\frac{e^{-\inv{2}\lambda_2^2(F_2^{-1}(u))^2}}{\lambda_2K_{2,1}\sqrt{-2\log u}\left(1+\frac{K_{2,2}\log(-\log u)}{\log u}+\frac{K_{2,3}}{\log u}+O\left(\left(\frac{\log|\log u|}{\log u}\right)^2\right)\right)}\\
%  = & \frac{e^{-\inv{2}\lambda_2^2(F_2^{-1}(u))^2}}{2\sqrt{\pi}\lambda_2K_{2,1}\sqrt{-\log u}\left(1+\frac{K_{2,2}\log(-\log u)}{\log u}+\frac{K_{2,3}}{\log u}+O\left(\left(\frac{\log(-\log u)}{\log u}\right)^2\right)\right)}\\
  \sim  \inv{\sqrt{2\pi}} u^{\frac{\lambda_2^2}{1+\lambda_2^2}}|\log u|^{\frac{\lambda_2^2}{1+\lambda_2^2}}(2\pi\lambda_2)^{\frac{\lambda_2^2}{1+\lambda_2^2}},
%  \left(1+O\left(\frac{[\log(-\log u)]^2}{\log u}\right)\right)}{2\sqrt{\pi}\times \frac{\lambda_2}{\sqrt{1+\lambda_2^2}}\times \sqrt{-\log u}\left(1+O\left(\frac{\log(-\log u)}{\log u}\right)\right)}\\
%  =& \frac{\sqrt{1+\lambda_2^2}(2\pi\lambda_2)^{\frac{\lambda_2^2}{1+\lambda_2^2}}}{2\sqrt{\pi}\lambda_2}u^{\frac{\lambda_2^2}{1+\lambda_2^2}}|\log u|^{\frac{\lambda_2^2}{1+\lambda_2^2}-\inv{2}}\left(1+O\left(\frac{[\log(-\log u)]^2}{\log u}\right)\right)
%%  =& \frac{e^{2\lambda_2^2\times\inv{1+\lambda_2^2}\times \log(2\pi\lambda_2)/2}}{2\sqrt{\pi}\lambda_2/\sqrt{1+\lambda_2^2}} u^{\frac{\lambda_2^2}{1+\lambda_2^2}}|\log u|^{2\times \frac{\lambda_2^2}{1+\lambda_2^2}-\inv{2}}\\
%  =& \frac{\sqrt{1+\lambda_2^2}(2\pi\lambda_2)^{\frac{\lambda_2^2}{1+\lambda_2^2}}}{2\sqrt{\pi}\lambda_2}u^{\frac{\lambda_2^2}{1+\lambda_2^2}}|\log u|^{\frac{\lambda_2^2}{1+\lambda_2^2}-\inv{2}}, 
%  + O\left(\frac{\left(\log(-\log u)\right)^2}{\log u}\right)}}{\sqrt{2\pi} \lambda_2K_2\sqrt{-2\log u}}\\
% \sim & \frac{u^{K_2\lambda_2^2}(-\log u)^{K_2^2\lambda_2^2K_3}}{2\sqrt{\pi}\lambda_2K_2\sqrt{-\log u}}.
\end{align*}}
as $u\to 0^+$.

\end{proof}

\subsection{Proof of Lemma \ref{Lemma: exp A in RV}}

\begin{proof}
%On the other hand, let's consider $e^{-\inv{2}\left(\frac{F_1^{-1}(u)-\rho F_2^{-1}(u)}{\sqrt{1-\rho^2}}\right)^2} \sim \tau_1 u^{\theta}(-\log u)^{\tau_2}$. 
% 
The results follow directly from Lemma \ref{lemma: general form of G1F1-GF2} once we obtained $\gamma_1$, $\theta$, $\tau_1$ and $\tau_2$ for different combinations of $\lambda_1$ and $\lambda_2$. 

When $\lambda_1>0$, $\lambda_2>0$, we have $[G_1 = \inv{\sqrt{1-\rho^2}}$, $G_2 = \frac{\rho}{\sqrt{1-\rho^2}}$, $\gamma_1 = \sqrt{\frac{1+\lambda_2^2}{1+\lambda_1^2}}$, $K_{2,1} = \inv{\sqrt{1+\lambda_2^2}}$, $K_{2,2} = \inv{2}$, $K_{2,3} = \frac{\log(2\pi\lambda_2)}{2}$, $C_{1,1}=0$, $C_{1,2} = \frac{\log(\lambda_1/\lambda_2)}{2}]$, {\adb
\begin{align*}
\notag \frac{|A_1(u)|}{\sqrt{1-\rho^2}} &\sim |K_{2,1}(G_1\gamma_1-G_2)|\sqrt{-2\log u} \\
&= 
\inv{\sqrt{1-\rho^2}}\left|\inv{\sqrt{1+\lambda_1^2}}-\frac{\rho}{\sqrt{1+\lambda_2^2}}\right|\sqrt{-2\log u}, \quad \text{as $u\to \infty$}; \\ 
\theta = & K_{2,1}^2\left(G_1\gamma_1- G_2\right)^2  = \inv{1-\rho^2}\left(\inv{\sqrt{1+\lambda_1^2}}-\frac{\rho}{\sqrt{1+\lambda_2^2}}\right)^2;\\
\tau_1 =&  e^{2\left[K_{2,1}^2K_{2,3}(G_1\gamma_1- G_2)^2+K_{2,1}^2C_{1,2}G_1\gamma_1(G_1\gamma_1- G_2)\right]}\\
= & e^{2\left[\frac{\log(2\pi\lambda_2)}{2}\times \inv{1-\rho^2}\left(\inv{\sqrt{1+\lambda_1^2}}-\frac{\rho}{\sqrt{1+\lambda_2^2}}\right)^2 + \frac{\log(\lambda_1/\lambda_2)}{2}\inv{1-\rho^2}\inv{\sqrt{1+\lambda_1^2}}\left(\inv{\sqrt{1+\lambda_1^2}}-\frac{\rho}{\sqrt{1+\lambda_2^2}}\right)\right]}\\
= & e^{\log\left[ (2\pi\lambda_2)^{\inv{1-\rho^2}\left(\inv{\sqrt{1+\lambda_1^2}}-\frac{\rho}{\sqrt{1+\lambda_2^2}}\right)^2}\times \left(\frac{\lambda_1}{\lambda_2}\right)^{\inv{(1-\rho^2)\sqrt{1+\lambda_1^2}}\left(\inv{\sqrt{1+\lambda_1^2}}-\frac{\rho}{\sqrt{1+\lambda_2^2}}\right)}\right]}\\
= & (2\pi\lambda_2)^{\inv{1-\rho^2}\left(\inv{\sqrt{1+\lambda_1^2}}-\frac{\rho}{\sqrt{1+\lambda_2^2}}\right)^2}\times \left(\frac{\lambda_1}{\lambda_2}\right)^{\inv{(1-\rho^2)\sqrt{1+\lambda_1^2}}\left(\inv{\sqrt{1+\lambda_1^2}}-\frac{\rho}{\sqrt{1+\lambda_2^2}}\right)};\\
\tau_2 = & 2\left[K_{2,1}^2K_{2,2}(G_1\gamma_1-G_2)^2+K_{2,1}^2C_{1,1}G_1\gamma_1(G_1\gamma_1- G_2)\right]\\
= & 2\left[\inv{2}\times \inv{1-\rho^2}\left(\inv{\sqrt{1+\lambda_1^2}}-\frac{\rho}{\sqrt{1+\lambda_2^2}}\right)^2+0\right], \quad \text{as $C_{1,1}=0$};\\
= & \inv{1-\rho^2}\left(\inv{\sqrt{1+\lambda_1^2}}-\frac{\rho}{\sqrt{1+\lambda_2^2}}\right)^2.
\end{align*}}
When $\lambda_1<0$, $\lambda_2>0$, we have $[G_1 = \inv{\sqrt{1-\rho^2}}$, $G_2 = \frac{\rho}{\sqrt{1-\rho^2}}$, $\gamma_1 = \sqrt{1+\lambda_2^2}$, $K_{2,1} = \inv{\sqrt{1+\lambda_2^2}}$, $K_{2,2} = \inv{2}$, $K_{2,3} = \frac{\log(2\pi\lambda_2)}{2}$, $C_{1,1}=-\inv{4}$, $C_{1,2} = - \frac{\log(2\lambda_2\sqrt{\pi})}{2}]$, {\adb
\begin{align*}
\notag \frac{|A_1(u)|}{\sqrt{1-\rho^2}} &\sim  |K_{2,1}(G_1\gamma_1-G_2)|\sqrt{-2\log u}\\
&=
\inv{\sqrt{1-\rho^2}}\left|1-\frac{\rho}{\sqrt{1+\lambda_2^2}}\right|\sqrt{-2\log u}, \quad \text{as $u\to 0^+$};\\
\theta = & K_{2,1}^2\left(G_1\gamma_1-G_2\right)^2
=  \inv{1-\rho^2}\left(1-\frac{\rho}{\sqrt{1+\lambda_2^2}}\right)^2;\\
\tau_1 =&  e^{2\left[K_{2,1}^2K_{2,3}(G_1\gamma_1-G_2)^2+K_{2,1}^2C_{1,2}G_1\gamma_1(G_1\gamma_1-G_2)\right]}\\
= & e^{2\left[\frac{\log(2\pi\lambda_2)}{2}\inv{1-\rho^2}\left(1-\frac{\rho}{\sqrt{1+\lambda_2^2}}\right)^2- \frac{\log(2\lambda_2\sqrt{\pi})}{2}\inv{1-\rho^2}\left(1-\frac{\rho}{\sqrt{1+\lambda_2^2}}\right)\right]}\\
=& e^{\log\left[(2\pi\lambda_2)^{\inv{1-\rho^2}\left(1-\frac{\rho}{\sqrt{1+\lambda_2^2}}\right)^2}\left(2\lambda_2\sqrt{\pi}\right)^{-\inv{1-\rho^2}\left(1-\frac{\rho}{\sqrt{1+\lambda_2^2}}\right)}\right]}\\
= & (2\pi\lambda_2)^{\inv{1-\rho^2}\left(1-\frac{\rho}{\sqrt{1+\lambda_2^2}}\right)^2}\left(2\lambda_2\sqrt{\pi}\right)^{-\inv{1-\rho^2}\left(1-\frac{\rho}{\sqrt{1+\lambda_2^2}}\right)};\\
\tau_2 = & 2\left[K_{2,1}^2K_{2,2}(G_1\gamma_1-G_2)^2+K_{2,1}^2C_{1,1}G_1\gamma_1(G_1\gamma_1-G_2)\right]\\
= & 2\left[\inv{2}\times \inv{1-\rho^2}\left(1-\frac{\rho}{\sqrt{1+\lambda_2^2}}\right)^2-\inv{4}\inv{1-\rho^2}\left(1-\frac{\rho}{\sqrt{1+\lambda_2^2}}\right)\right]\\
= & \inv{1-\rho^2}\left(1-\frac{\rho}{\sqrt{1+\lambda_2^2}}\right)^2-\inv{2(1-\rho^2)}\left(1-\frac{\rho}{\sqrt{1+\lambda_2^2}}\right).
\end{align*}}
When $\lambda_1>0$, $\lambda_2<0$, we have $[G_1 = \inv{\sqrt{1-\rho^2}}$, $G_2 = \frac{\rho}{\sqrt{1-\rho^2}}$, $\gamma_1 = \frac{1}{\sqrt{1+\lambda_1^2}}$, $K_{2,1} = 1$, $K_{2,2} = \inv{4}$, $K_{2,3} = \frac{\log(\pi)}{4}$, $C_{1,1}=\inv{4}$, $C_{1,2} = \frac{\log(2\lambda_1\sqrt{\pi})}{2}]$, {\adb
\begin{align*}
\notag \frac{|A_1(u)|}{\sqrt{1-\rho^2}} &\sim |K_{2,1}(G_1\gamma_1-G_2)|\sqrt{-2\log u} \\
&= \inv{\sqrt{1-\rho^2}}\left|\inv{\sqrt{1+\lambda_1^2}}-\rho\right|\sqrt{-2\log u}, \quad \text{as $u\to 0^+$};\\
\theta = & K_{2,1}^2\left(G_1\gamma_1-G_2\right)^2
=  \inv{1-\rho^2}\left(\inv{\sqrt{1+\lambda_1^2}}-\rho\right)^2;\\
\tau_1 =&  e^{2\left[K_{2,1}^2K_{2,3}(G_1\gamma_1-G_2)^2+K_{2,1}^2C_{1,2}G_1\gamma_1(G_1\gamma_1-G_2)\right]}\\
= & e^{2\left[\frac{\log\pi}{4}\times \inv{1-\rho^2}\left(\inv{\sqrt{1+\lambda_1^2}}-\rho\right)^2 + \frac{\log(2\lambda_1\sqrt{\pi})}{4}\inv{1-\rho^2}\inv{\sqrt{1+\lambda_1^2}}\left(\inv{\sqrt{1+\lambda_1^2}}-\rho\right)\right]}\\
=& e^{\log\left[\pi^{\inv{2(1-\rho^2)}\left(\inv{\sqrt{1+\lambda_1^2}}-\rho\right)^2}(2\lambda_1\sqrt{\pi})^{\inv{2(1-\rho^2)\sqrt{1+\lambda_1^2}}\left(\inv{\sqrt{1+\lambda_1^2}}-\rho\right)}\right]}\\
=& \pi^{\inv{2(1-\rho^2)}\left(\inv{\sqrt{1+\lambda_1^2}}-\rho\right)^2}(2\lambda_1\sqrt{\pi})^{\inv{2(1-\rho^2)\sqrt{1+\lambda_1^2}}\left(\inv{\sqrt{1+\lambda_1^2}}-\rho\right)};\\
\tau_2 = & 2\left[K_{2,1}^2K_{2,2}(G_1\gamma_1-G_2)^2+K_{2,1}^2C_{1,1}\gamma_1(G_1\gamma_1-G_2)\right]\\
= & 2\left[\inv{4(1-\rho^2)}\left(\inv{\sqrt{1+\lambda_1^2}}-\rho\right)^2+\inv{4(1-\rho^2)\sqrt{1+\lambda_1^2}}\left(\inv{\sqrt{1+\lambda_1^2}}-\rho\right)\right]\\
= & \inv{2(1-\rho^2)}\left(\inv{\sqrt{1+\lambda_1^2}}-\rho\right)^2+\inv{2(1-\rho^2)\sqrt{1+\lambda_1^2}}\left(\inv{\sqrt{1+\lambda_1^2}}-\rho\right).
\end{align*}}
When $\lambda_1=0$, $\lambda_2>0$, we have $[G_1 = \inv{\sqrt{1-\rho^2}}$, $G_2 = \frac{\rho}{\sqrt{1-\rho^2}}$, 
$\gamma_1 = \sqrt{1+\lambda_2^2}$, $K_{2,1} = \inv{\sqrt{1+\lambda_2^2}}$, $K_{2,2} = \inv{2}$, $K_{2,3} = \frac{\log(2\pi\lambda_2)}{2}$, $C_{1,1}=-\inv{4}$, $C_{1,2} = -\frac{\log(\lambda_2\sqrt{\pi})}{2}]$, {\adb
\begin{align*}
\notag \frac{|A_1(u)|}{\sqrt{1-\rho^2}} &\sim |K_{2,1}(G_1\gamma_1-G_2)|\sqrt{-2\log u} \\
&=\inv{\sqrt{1-\rho^2}}\left|1-\frac{\rho}{\sqrt{1+\lambda_2^2}}\right|\sqrt{-2\log u}, \quad \text{as $u\to 0^+$}; \\
\theta = & K_{2,1}^2\left(G_1\gamma_1-G_2\right)^2
=  \inv{1-\rho^2}\left(1-\frac{\rho}{\sqrt{1+\lambda_2^2}}\right)^2;\\
\tau_1 =&  e^{2\left[K_{2,1}^2K_{2,3}(G_1\gamma_1-G_2)^2+K_{2,1}^2C_{1,2}G_1\gamma_1(G_1\gamma_1-G_2)\right]}\\
= & e^{2\left[\frac{\log(2\pi\lambda_2)}{2}\inv{1-\rho^2}\left(1-\frac{\rho}{\sqrt{1+\lambda_2^2}}\right)^2-\frac{\log(\lambda_2\sqrt{\pi})}{2}\inv{1-\rho^2}\left(1-\frac{\rho}{\sqrt{1+\lambda_2^2}}\right)\right]}\\
=& e^{\log\left[(2\pi\lambda_2)^{\inv{1-\rho^2}\left(1-\frac{\rho}{\sqrt{1+\lambda_2^2}}\right)^2}(\lambda_2\sqrt{\pi})^{-\inv{1-\rho^2}\left(1-\frac{\rho}{\sqrt{1+\lambda_2^2}}\right)}\right]}\\
= & (2\pi\lambda_2)^{\inv{1-\rho^2}\left(1-\frac{\rho}{\sqrt{1+\lambda_2^2}}\right)^2}(\lambda_2\sqrt{\pi})^{-\inv{1-\rho^2}\left(1-\frac{\rho}{\sqrt{1+\lambda_2^2}}\right)};\\
\tau_2 = & 2\left[K_{2,1}^2K_{2,2}(G_1\gamma_1-G_2)^2+K_{2,1}^2C_{1,1}G_1\gamma_1(G_1\gamma_1-G_2)\right]\\
= & 2\left[\inv{2(1-\rho^2)}\left(1-\frac{\rho}{\sqrt{1+\lambda_2^2}}\right)^2-\inv{4(1-\rho^2)}\left(1-\frac{\rho}{\sqrt{1+\lambda_2^2}}\right)\right]\\
= & \inv{1-\rho^2}\left(1-\frac{\rho}{\sqrt{1+\lambda_2^2}}\right)^2 - \inv{2(1-\rho^2)}\left(1-\frac{\rho}{\sqrt{1+\lambda_2^2}}\right).
\end{align*} 
}
When $\lambda_1>0$, $\lambda_2=0$, we have $[G_1 = \inv{\sqrt{1-\rho^2}}$, $G_2 = \frac{\rho}{\sqrt{1-\rho^2}}$, $\gamma_1 = \frac{1}{\sqrt{1+\lambda_1^2}}$, $K_{2,1} = 1$, $K_{2,2} = \inv{4}$, $K_{2,3} = \frac{\log(4\pi)}{4}$, $C_{1,1}=\inv{4}$, $C_{1,2} = \frac{\log(\lambda_1\sqrt{\pi})}{2}]$, {\adb
\begin{align*}
\notag \frac{|A_1(u)|}{\sqrt{1-\rho^2}} &\sim |K_{2,1}(G_1\gamma_1-G_2)|\sqrt{-2\log u} \\
&= \inv{\sqrt{1-\rho^2}}\left|\inv{\sqrt{1+\lambda_1^2}}-\rho\right|\sqrt{-2\log u}, \quad \text{as $u\to 0^+$}; \\
\theta = & K_{2,1}^2\left(G_1\gamma_1-G_2\right)^2 =  \inv{1-\rho^2}\left(\inv{\sqrt{1+\lambda_1^2}}-\rho\right)^2\\
\tau_1 =&  e^{2\left[K_{2,1}^2K_{2,3}(G_1\gamma_1- G_2)^2+K_{2,1}^2C_{1,2}G_1\gamma_1(G_1\gamma_1-G_2)\right]}\\
= & e^{2\left[\frac{\log4\pi}{4}\times \inv{1-\rho^2}\left(\inv{\sqrt{1+\lambda_1^2}}-\rho\right)^2+\frac{\log(\lambda_1\sqrt{\pi})}{2}\times\inv{1-\rho^2}\inv{\sqrt{1+\lambda_1^2}}\left(\inv{\sqrt{1+\lambda_1^2}}-\rho\right)\right]}\\
= & e^{\log\left[(4\pi)^{\inv{2(1-\rho^2)}\left(\inv{\sqrt{1+\lambda_1^2}}-\rho\right)^2}(\lambda_1\sqrt{\pi})^{\inv{(1-\rho^2)\sqrt{1+\lambda_1^2}}\left(\inv{\sqrt{1+\lambda_1^2}}-\rho\right)}\right]}\\
=& (4\pi)^{\inv{2(1-\rho^2)}\left(\inv{\sqrt{1+\lambda_1^2}}-\rho\right)^2}(\lambda_1\sqrt{\pi})^{\inv{(1-\rho^2)\sqrt{1+\lambda_1^2}}\left(\inv{\sqrt{1+\lambda_1^2}}-\rho\right)};\\
\tau_2 = & 2\left[K_{2,1}^2K_{2,2}(G_1\gamma_1- G_2)^2+K_{2,1}^2C_{1,1}G_1\gamma_1(G_1\gamma_1- G_2)\right]\\
= & 2\left[1\times \inv{4}\times \inv{1-\rho^2}\left(\inv{\sqrt{1+\lambda_1^2}}-\rho\right)^2+\inv{4}\times \inv{1-\rho^2}\times \inv{\sqrt{1+\lambda_1^2}}\left(\inv{\sqrt{1+\lambda_1^2}}-\rho\right)\right]\\
= & \inv{2(1-\rho^2)}\left(\inv{\sqrt{1+\lambda_1^2}}-\rho\right)^2+\inv{2(1-\rho^2)\sqrt{1+\lambda_1^2}}\left(\inv{\sqrt{1+\lambda_1^2}}-\rho\right).
\end{align*}}
When $\lambda_1<0$, $\lambda_2<0$, we have $[G_1 = \inv{\sqrt{1-\rho^2}}$, $G_2 = \frac{\rho}{\sqrt{1-\rho^2}}$, $\gamma_1 =1$, $K_{2,1} = 1$, $K_{2,2} = \inv{4}$, $K_{2,3} = \frac{\log\pi}{4}$, $C_{1,1}=  C_{1,2} = 0]$, {\adb
\begin{align*}
\notag \frac{|A_1(u)|}{\sqrt{1-\rho^2}} &\sim |K_{2,1}(G_1\gamma_1-G_2)|\sqrt{-2\log u}=
\sqrt{\frac{1-\rho}{1+\rho}}\sqrt{-2\log u}, \quad \text{as $u\to 0^+$}; \\
\theta = & K_{2,1}^2\left(G_1\gamma_1- G_2\right)^2 =   \frac{1-\rho}{1+\rho};\\
\tau_1 =&  e^{2\left[K_{2,1}^2K_{2,3}(G_1\gamma_1- G_2)^2+K_{2,1}^2C_{1,2}G_1\gamma_1(G_1\gamma_1- G_2)\right]}\\
= & e^{2\left[\frac{\log\pi}{4}\frac{1-\rho}{1+\rho} + 0\right]}, \quad \text{as $C_{1,2}=0$};\\
= & \pi^{\inv{2}\left(\frac{1-\rho}{1+\rho}\right)};\\
\tau_2 = & 2\left[K_{2,1}^2K_{2,2}(G_1\gamma_1- G_2)^2+K_{2,1}^2C_{1,1}G_1\gamma_1(G_1\gamma_1- G_2)\right]\\
= & 2\left[\inv{4}\times\frac{1-\rho}{1+\rho}+0\right], \quad \text{as $C_{1,1}=0$};\\
=& \inv{2}\left(\frac{1-\rho}{1+\rho}\right).
\end{align*}}
When $\lambda_1<0$, $\lambda_2=0$, we have $[G_1 = \inv{\sqrt{1-\rho^2}}$, $G_2 = \frac{\rho}{\sqrt{1-\rho^2}}$, $\gamma_1 = 1$, $K_{2,1} = 1$, $K_{2,2} = \inv{4}$, $K_{2,3} = \frac{\log(4\pi)}{4}$, $C_{1,1}=0$, $C_{1,2} = -\frac{\log 2}{2}]$, {\adb
\begin{align*}\notag \frac{|A_1(u)|}{\sqrt{1-\rho^2}} &\sim |K_{2,1}(G_1\gamma_1-G_2)|\sqrt{-2\log u}=
\sqrt{\frac{1-\rho}{1+\rho}}\sqrt{-2\log u}, \quad \text{as $u\to 0^+$}; \\
\theta = & K_{2,1}^2\left(G_1\gamma_1- G_2\right)^2 =   \frac{1-\rho}{1+\rho};\\
\tau_1 =&  e^{2\left[K_{2,1}^2K_{2,3}(G_1\gamma_1- G_2)^2+K_{2,1}^2C_{1,2}G_1\gamma_1(G_1\gamma_1- G_2)\right]}\\
=& e^{2\left[\frac{\log(4\pi)}{4}\times \frac{1-\rho}{1+\rho}+1^2\times (-\frac{\log 2}{2})\inv{\sqrt{1-\rho^2}}
\left(\frac{1-\rho}{\sqrt{1-\rho^2}}\right)\right]}\\
=& e^{\log\left[(4\pi)^{\inv{2}\left(\frac{1-\rho}{1+\rho}\right)}\times 2^{-\inv{1+\rho}}\right]}\\
=& (4\pi)^{\inv{2}\left(\frac{1-\rho}{1+\rho}\right)}\times 2^{-\inv{1+\rho}}\\
\tau_2 = & 2\left[K_{2,1}^2K_{2,2}(G_1\gamma_1- G_2)^2+K_{2,1}^2C_{1,1}G_1\gamma_1(G_1\gamma_1- G_2)\right]\\
=& 2\left[1^2\times \inv{4}\times\frac{1-\rho}{1+\rho}+0\right], \quad \text{as $C_{1,1}=0$};\\
=& \inv{2}\left(\frac{1-\rho}{1+\rho}\right).
\end{align*}}
When $\lambda_1=0$, $\lambda_2<0$, we have $[G_1 = \inv{\sqrt{1-\rho^2}}$, $G_2 = \frac{\rho}{\sqrt{1-\rho^2}}$, $\gamma_1 = 1$, $K_{2,1} =1$, $K_{2,2} = \inv{4}$, $K_{2,3} = \frac{\log(\pi)}{4}$, $C_{1,1}=0$, $C_{1,2} = \frac{\log 2}{2}]$, {\adb
\begin{align*}
\notag \frac{|A_1(u)|}{\sqrt{1-\rho^2}} &\sim |K_{2,1}(G_1\gamma_1-G_2)|\sqrt{-2\log u}=
\sqrt{\frac{1-\rho}{1+\rho}}\sqrt{-2\log u}, \quad \text{as $u\to 0^+$}; \\
\theta = & K_{2,1}^2\left(G_1\gamma_1- G_2\right)^2 =   \frac{1-\rho}{1+\rho};\\
\tau_1 =&  e^{2\left[K_{2,1}^2K_{2,3}(G_1\gamma_1- G_2)^2+K_{2,1}^2C_{1,2}G_1\gamma_1(G_1\gamma_1- G_2)\right]}\\
=& e^{2\left[\frac{\log(\pi)}{4}\times \frac{1-\rho}{1+\rho}+\frac{\log 2}{2}\times \inv{\sqrt{1-\rho^2}}\left(\frac{1-\rho}{\sqrt{1-\rho^2}}\right)\right]}\\
=& (\pi)^{\inv{2}\left(\frac{1-\rho}{1+\rho}\right)}\times 2^{\inv{1+\rho}}; \\
\tau_2 = & 2\left[K_{2,1}^2K_{2,2}(G_1\gamma_1- G_2)^2+K_{2,1}^2C_{1,1}G_1\gamma_1(G_1\gamma_1- G_2)\right]\\
=& 2\left[\inv{4}\times\frac{1-\rho}{1+\rho}+0\right], \quad \text{as $C_{1,1}=0$};\\
=& \inv{2}\left(\frac{1-\rho}{1+\rho}\right).
\end{align*}}
When $\lambda_1 = \lambda_2=0$, we have $F_1^{-1}(u) = F_2^{-1}(u) = \Phi^{-1}(u)$, and
\begin{align*}
 \frac{|A_1(u)|}{\sqrt{1-\rho^2}} & = \frac{1-\rho}{\sqrt{1-\rho^2}}\left|\Phi^{-1}(u)\right| \sim \sqrt{\frac{1-\rho}{1+\rho}}\sqrt{-2\log u}, \quad \text{as $u\to 0^+$}.
\end{align*}

\end{proof}

\subsection{Proof of Lemma \ref{Lemma: exp B in RV}}

\begin{proof}
The results follow directly from Lemma \ref{lemma: general form of G1F1-GF2} once we obtained $\gamma_1$, $\theta$, $\tau_1$ and $\tau_2$ for different combinations of $\lambda_1$ and $\lambda_2$. 

When $\lambda_1$, $\lambda_2>0$, we have $[G_1 = \alpha_1$, $G_2 = -\alpha_2$, $\gamma_1 = \sqrt{\frac{1+\lambda_2^2}{1+\lambda_1^2}}$, $K_{2,1} = \inv{\sqrt{1+\lambda_2^2}}$, $K_{2,2} = \inv{2}$, $K_{2,3} = \frac{\log(2\pi\lambda_2)}{2}$, $C_{1,1}=0$, $C_{1,2} = \frac{\log(\lambda_1/\lambda_2)}{2}]$, {\adb
\begin{align*}
\notag  |B(u)| & \sim |K_{2,1}(G_1\gamma_1-G_2)|\sqrt{-2\log u} = 
\left|\frac{\alpha_1}{\sqrt{1+\lambda_1^2}}+\frac{\alpha_2}{\sqrt{1+\lambda_2^2}}\right|\sqrt{-2\log u}, \quad \text{as $u\to 0^+$}; \\
\theta =& K_{2,1}^2(G_1\gamma_1-G_2)^2  = \left(\frac{\alpha_1}{\sqrt{1+\lambda_1^2}}+\frac{\alpha_2}{\sqrt{1+\lambda_2^2}}\right)^2\\
\tau_1 =& e^{2\left[K_{2,1}^2K_{2,3}(G_1\gamma_1- G_2)^2+K_{2,1}^2C_{1,2}\gamma_1G_1(G_1\gamma_1- G_2)\right]}\\
=& e^{2\left[\inv{1+\lambda_2^2}\frac{\log(2\pi\lambda_2)}{2}\left(\alpha_1\sqrt{\frac{1+\lambda_2^2}{1+\lambda_1^2}}+\alpha_2\right)^2+\inv{1+\lambda_2^2}\frac{\log(\lambda_1/\lambda_2)}{2}\alpha_1\sqrt{\frac{1+\lambda_2^2}{1+\lambda_1^2}}\left(\alpha_1\sqrt{\frac{1+\lambda_2^2}{1+\lambda_1^2}}+\alpha_2\right)\right]}\\
= & e^{\log(2\pi\lambda_2)\left(\frac{\alpha_1}{\sqrt{1+\lambda_2^2}}+\frac{\alpha_2}{\sqrt{1+\lambda_2^2}}\right)^2+\left[\log\left(\frac{\lambda_1}{\lambda_2}\right)\right]\frac{\alpha_1}{\sqrt{1+\lambda_1^2}}\left(\frac{\alpha_1}{\sqrt{1+\lambda_1^2}}+\frac{\alpha_2}{\sqrt{1+\lambda_2^2}}\right)}\\
=& e^{\log\left[(2\pi\lambda_2)^{\left(\frac{\alpha_1}{\sqrt{1+\lambda_2^2}}+\frac{\alpha_2}{\sqrt{1+\lambda_2^2}}\right)^2}\left(\frac{\lambda_1}{\lambda_2}\right)^{\frac{\alpha_1}{\sqrt{1+\lambda_1^2}}\left(\frac{\alpha_1}{\sqrt{1+\lambda_1^2}}+\frac{\alpha_2}{\sqrt{1+\lambda_2^2}}\right)}\right]}\\
= & (2\pi\lambda_2)^{\left(\frac{\alpha_1}{\sqrt{1+\lambda_2^2}}+\frac{\alpha_2}{\sqrt{1+\lambda_2^2}}\right)^2}\left(\frac{\lambda_1}{\lambda_2}\right)^{\frac{\alpha_1}{\sqrt{1+\lambda_1^2}}\left(\frac{\alpha_1}{\sqrt{1+\lambda_1^2}}+\frac{\alpha_2}{\sqrt{1+\lambda_2^2}}\right)}\\
\tau_2 =& 2\left[K_{2,1}^2K_{2,2}(G_1\gamma_1- G_2)^2+K_{2,1}^2C_{1,1}G_1\gamma_1\left(G_1\gamma_1- G_2\right)\right]\\
=& 2\left[\left(\inv{\sqrt{1+\lambda_2^2}}\right)^2\times \inv{2}\times \left(\alpha_1\sqrt{\frac{1+\lambda_2^2}{1+\lambda_1^2}}+\alpha_2\right)^2+0\right], \quad \text{as $C_{1,1}=0$}; \\
= & \left(\frac{\alpha_1}{\sqrt{1+\lambda_1^2}}+\frac{\alpha_2}{\sqrt{1+\lambda_2^2}}\right)^2.
\end{align*}}
When $\lambda_1<0$, $\lambda_2>0$, we have $[G_1 = \alpha_1$, $G_2 = - \alpha_2$, $\gamma_1 = \sqrt{1+\lambda_2^2}$, $K_{2,1} = \inv{\sqrt{1+\lambda_2^2}}$, $K_{2,2} = \inv{2}$, $K_{2,3} = \frac{\log(2\pi\lambda_2)}{2}$, $C_{1,1}=-\inv{4}$, $C_{1,2} = -\frac{\log(2\lambda_2\sqrt{\pi})}{2}]$, {\adb
\begin{align*}
\notag  |B(u)| & \sim |K_{2,1}(G_1\gamma_1-G_2)|\sqrt{-2\log u} = 
\left|\alpha_1+\frac{\alpha_2}{\sqrt{1+\lambda_2^2}}\right|\sqrt{-2\log u}, \quad \text{as $u\to 0^+$}; \\
\theta = & K_{2,1}^2(G_1\gamma_1-G_2)^2
=  \left(\alpha_1+\frac{\alpha_2}{\sqrt{1+\lambda_2^2}}\right)^2;\\
\tau_1 =& e^{2\left[K_{2,1}^2K_{2,3}(G_1\gamma_1- G_2)^2+K_{2,1}^2C_{1,2}G_1\gamma_1(G_1\gamma_1- G_2)\right]}\\
=& e^{2\left[\frac{\log(2\pi\lambda_2)}{2}\left(\alpha_1+\frac{\alpha_2}{\sqrt{1+\lambda_2}}\right)^2-\frac{\log(2\lambda_2\sqrt{\pi})}{2}\alpha_1\left(\alpha_1+\frac{\alpha_2}{\sqrt{1+\lambda_2^2}}\right)\right]}\\
= & e^{\log\left[(2\pi\lambda_2)^{\left(\alpha_1+\frac{\alpha_2}{\sqrt{1+\lambda_2^2}}\right)^2}(2\lambda_2\sqrt{\pi})^{-\alpha_1\left(\alpha_1+\frac{\alpha_2}{\sqrt{1+\lambda_2^2}}\right)}\right]}\\
= & (2\pi\lambda_2)^{\left(\alpha_1+\frac{\alpha_2}{\sqrt{1+\lambda_2^2}}\right)^2}\left(2\lambda_2\sqrt{\pi}\right)^{-\alpha_1\left(\alpha_1+\frac{\alpha_2}{\sqrt{1+\lambda_2^2}}\right)}\\
\tau_2 = & 2\left[K_{2,1}^2K_{2,2}(G_1\gamma_1- G_2)^2+K_{2,1}^2C_{1,1}G_1\gamma_1(G_1\gamma_1- G_2)\right]\\
= & 2\left[\inv{2}\times \left(\alpha_1+\frac{\alpha_2}{\sqrt{1+\lambda_2^2}}\right)^2-\inv{4}\alpha_1\left(\alpha_1+\frac{\alpha_2}{\sqrt{1+\lambda_2^2}}\right)\right]\\
= & \left(\alpha_1+\frac{\alpha_2}{\sqrt{1+\lambda_2^2}}\right)^2-\frac{\alpha_1}{2}\left(\alpha_1+\frac{\alpha_2}{\sqrt{1+\lambda_2^2}}\right).
\end{align*}}
When $\lambda_1>0$, $\lambda_2<0$, we have $[G_1 = \alpha_1$, $G_2 = - \alpha_2$, $\gamma_1 = \inv{\sqrt{1+\lambda_1^2}}$, $K_{2,1} = 1$, $K_{2,2} = \inv{4}$, $K_{2,3} = \frac{\log\pi}{4}$, $C_{1,1}=\inv{4}$, $C_{1,2} = \frac{\log(2\lambda_1\sqrt{\pi})}{2}]$, {\adb
\begin{align*}
\notag  |B(u)| & \sim |K_{2,1}(G_1\gamma_1-G_2)|\sqrt{-2\log u} = 
\left|\frac{\alpha_1}{\sqrt{1+\lambda_1^2}}+ \alpha_2\right|\sqrt{-2\log u}, \quad \text{as $u\to 0^+$}; \\
\theta = & K_{2,1}^2\left(G_1\gamma_1- G_2\right)^2
= \left(\frac{\alpha_1}{\sqrt{1+\lambda_1^2}}+\alpha_2\right)^2;\\
\tau_1 =&  e^{2\left[K_{2,1}^2K_{2,3}(G_1\gamma_1- G_2)^2+K_{2,1}^2C_{1,2}G_1\gamma_1(G_1\gamma_1- G_2)\right]}\\
= & e^{2\left[1^2\times \frac{\log\pi}{4}\left(\frac{\alpha_1}{\sqrt{1+\lambda_1^2}}+\alpha_2\right)^2+1^2\times \frac{\log(2\lambda_1\sqrt{\pi})}{2}\times \frac{\alpha_1}{\sqrt{1+\lambda_1^2}}\left(\frac{\alpha_1}{\sqrt{1+\lambda_1^2}}+\alpha_2\right)\right]}\\
=& e^{\log\left[\pi^{\inv{2}\left(\frac{\alpha_1}{\sqrt{1+\lambda_1^2}}+\alpha_2\right)^2}\times (2\lambda_1\sqrt{\pi})^{\frac{\alpha_1}{\sqrt{1+\lambda_1^2}}\left(\frac{\alpha_1}{\sqrt{1+\lambda_1^2}}+\alpha_2\right)}\right]}\\
=& \pi^{\inv{2}\left(\frac{\alpha_1}{\sqrt{1+\lambda_1^2}}+\alpha_2\right)^2}\times (2\lambda_1\sqrt{\pi})^{\frac{\alpha_1}{\sqrt{1+\lambda_1^2}}\left(\frac{\alpha_1}{\sqrt{1+\lambda_1^2}}+\alpha_2\right)};\\
\tau_2 = & 2\left[K_{2,1}^2K_{2,2}(G_1\gamma_1- G_2)^2+K_{2,1}^2C_{1,1}G_1\gamma_1(G_1\gamma_1- G_2)\right]\\
= & 2\left[1^2\times \inv{4}\left(\frac{\alpha_1}{\sqrt{1+\lambda_1^2}}+\alpha_1\right)^2+1^2\times \inv{4}\times \frac{\alpha_1}{\sqrt{1+\lambda_1^2}}\left(\frac{\alpha_1}{\sqrt{1+\lambda_1^2}}+\alpha_2\right)\right]\\
=& \inv{2}\left(\frac{\alpha_1}{\sqrt{1+\lambda_1^2}}+\alpha_2\right)^2+\inv{2}\frac{\alpha_1}{\sqrt{1+\lambda_1^2}}\left(\frac{\alpha_1}{\sqrt{1+\lambda_1^2}}+\alpha_2\right).
\end{align*}}
When $\lambda_1=0$, $\lambda_2>0$, we have $[G_1 = \alpha_1$, $G_2 = -\alpha_2$, $\gamma_1 = \sqrt{1+\lambda_2^2}$, $K_{2,1} = \inv{\sqrt{1+\lambda_2^2}}$, $K_{2,2} = \inv{2}$, $K_{2,3} = \frac{\log(2\pi\lambda_2)}{2}$, $C_{1,1}=-\inv{4}$, $C_{1,2} = -\frac{\log(\lambda_2\sqrt{\pi})}{2}]$, {\adb
\begin{align*}
\notag  |B(u)| & \sim |K_{2,1}(G_1\gamma_1-G_2)|\sqrt{-2\log u} = 
\left|\alpha_1+\frac{\alpha_2}{\sqrt{1+\lambda_2^2}}\right|\sqrt{-2\log u}, \quad \text{as $u\to 0^+$}; \\
\theta = & K_{2,1}^2\left(G_1\gamma_1- G_2\right)^2
= \left(\alpha_1+\frac{\alpha_2}{\sqrt{1+\lambda_2^2}}\right)^2;\\
\tau_1 =&  e^{2\left[K_{2,1}^2K_{2,3}(G_1\gamma_1-G_2)^2+K_{2,1}^2C_{1,2}G_1\gamma_1(G_1\gamma_1-G_2)\right]}\\
=& e^{2\left[\frac{\log(2\pi\lambda_2)}{2}\left(\alpha_1+\frac{\alpha_2}{\sqrt{1+\lambda_2^2}}\right)^2- \frac{\log(\lambda_2\sqrt{\pi})}{2}\alpha_1\left(\alpha_1+\frac{\alpha_2}{\sqrt{1+\lambda_2^2}}\right)\right]}\\
=& e^{\log\left[(2\pi\lambda_2)^{\left(\alpha_1+\frac{\alpha_2}{\sqrt{1+\lambda_2^2}}\right)^2}(\lambda_2\sqrt{\pi})^{-\alpha_1\left(\alpha_1+\frac{\alpha_2}{\sqrt{1+\lambda_2^2}}\right)}\right]}\\
=& (2\pi\lambda_2)^{\left(\alpha_1+\frac{\alpha_2}{\sqrt{1+\lambda_2^2}}\right)^2}(\lambda_2\sqrt{\pi})^{-\alpha_1\left(\alpha_1+\frac{\alpha_2}{\sqrt{1+\lambda_2^2}}\right)};\\
\tau_2 = & 2\left[K_{2,1}^2K_{2,2}(G_1\gamma_1-G_2)^2+K_{2,1}^2C_{1,1}G_1\gamma_1(G_1\gamma_1-G_2)\right]\\
= & 2\left[\inv{2}\left(\alpha_1+\frac{\alpha_2}{\sqrt{1+\lambda_2^2}}\right)^2-\inv{4}\alpha_1\left(\alpha_1+\frac{\alpha_2}{\sqrt{1+\lambda_2^2}}\right)\right]\\
= & \left(\alpha_1+\frac{\alpha_2}{\sqrt{1+\lambda_2^2}}\right)^2-\frac{\alpha_1}{2}\left(\alpha_1+\frac{\alpha_2}{\sqrt{1+\lambda_2^2}}\right).
\end{align*}}
When $\lambda_1>0$, $\lambda_2=0$, we have $[G_1 = \alpha_1$, $G_2 = - \alpha_2$, $\gamma_1 = \inv{\sqrt{1+\lambda_1^2}}$, $K_{2,1} = 1$, $K_{2,2} = \inv{4}$, $K_{2,3} = \frac{\log(4\pi)}{4}$, $C_{1,1}= \inv{4}$, $C_{1,2} = \frac{\log(\lambda_1\sqrt{\pi})}{2}]$, {\adb
\begin{align*}
\notag  |B(u)| & \sim |K_{2,1}(G_1\gamma_1-G_2)|\sqrt{-2\log u} = 
\left|\frac{\alpha_1}{\sqrt{1+\lambda_1^2}}+ \alpha_2\right|\sqrt{-2\log u}, \quad \text{as $u\to 0^+$}; \\
\theta = & K_{2,1}^2\left(G_1\gamma_1-G_2\right)^2
=  1^2 \times \left(\frac{\alpha_1}{\sqrt{1+\lambda_1^2}}+\alpha_2\right)^2;\\
\tau_1 =&  e^{2\left[K_{2,1}^2K_{2,3}(G_1\gamma_1- G_2)^2+K_{2,1}^2C_{1,2}G_1\gamma_1(G_1\gamma_1- G_2)\right]}\\
= & e^{2\left[\frac{\log(4\pi)}{4}\left(\frac{\alpha_1}{\sqrt{1+\lambda_1^2}}+\alpha_2\right)^2+\frac{\log(\lambda_1\sqrt{\pi})}{2}\frac{\alpha_1}{\sqrt{1+\lambda_1^2}}\left(\frac{\alpha_1}{\sqrt{1+\lambda_1^2}}+\alpha_2\right)\right]}\\
= & e^{\log\left[(4\pi)^{\inv{2}\left(\frac{\alpha_1}{\sqrt{1+\lambda_1^2}}+\alpha_2\right)^2}(\lambda_1\sqrt{\pi})^{\frac{\alpha_1}{\sqrt{1+\lambda_1^2}}\left(\frac{\alpha_1}{\sqrt{1+\lambda_1^2}}+\alpha_2\right)}\right]}\\
= & (4\pi)^{\inv{2}\left(\frac{\alpha_1}{\sqrt{1+\lambda_1^2}}+\alpha_2\right)^2}(\lambda_1\sqrt{\pi})^{\frac{\alpha_1}{\sqrt{1+\lambda_1^2}}\left(\frac{\alpha_1}{\sqrt{1+\lambda_1^2}}+\alpha_2\right)};\\
\tau_2 = & 2\left[K_{2,1}^2K_{2,2}(G_1\gamma_1- G_2)^2+K_{2,1}^2C_{1,1}G_1\gamma_1(G_1\gamma_1- G_2)\right]\\
= & 2\left[1^2\times \inv{4}\left(\frac{\alpha_1}{\sqrt{1+\lambda_1^2}}+\alpha_2\right)^2+1^2\times \inv{4}\frac{\alpha_1}{\sqrt{1+\lambda_1^2}}\left(\frac{\alpha_1}{\sqrt{1+\lambda_1^2}}+\alpha_2\right)\right]\\
= & \inv{2}\left(\frac{\alpha_1}{\sqrt{1+\lambda_1^2}}+\alpha_2\right)^2+\inv{2}\frac{\alpha_1}{\sqrt{1+\lambda_1^2}}\left(\frac{\alpha_1}{\sqrt{1+\lambda_1^2}}+\alpha_2\right).
\end{align*}}
When $\lambda_1<0$, $\lambda_2<0$, we have $[G_1 = \alpha_1$, $G_2 = - \alpha_2$, $\gamma_1 = 1$, $K_{2,1} = 1$, $K_{2,2} = \inv{4}$, $K_{2,3} = \frac{\log\pi}{4}$, $C_{1,1}=0=C_{1,2}]$, {\adb
\begin{align*}
\notag  |B(u)| & \sim |K_{2,1}(G_1\gamma_1-G_2)|\sqrt{-2\log u} = 
\left|\alpha_1 + \alpha_2\right|\sqrt{-2\log u}, \quad \text{as $u\to 0^+$}; \\
\theta = & K_{2,1}^2\left(G_1\gamma_1- G_2\right)^2 = (\alpha_1+\alpha_2)^2;\\
\tau_1 =&  e^{2\left[K_{2,1}^2K_{2,3}(G_1\gamma_1- G_2)^2+K_{2,1}^2C_{1,2}G_1\gamma_1(G_1\gamma_1- G_2)\right]}\\
= & e^{2\left[\frac{\log\pi}{4}\left(\alpha_1+\alpha_2\right)^2+0\right]}, \quad \text{as $C_{1,2}=0$},\\
= & \pi^{\frac{(\alpha_1+\alpha_2)^2}{2}};\\
\tau_2 = & 2\left[K_{2,1}^2K_{2,2}(G_1\gamma_1- G_2)^2+K_{2,1}^2C_{1,1}G_1\gamma_1(G_1\gamma_1- G_2)\right]\\
= & 2\left[\inv{4}(\alpha_1+\alpha_2)^2+0\right] \quad \text{as $C_{1,1}=0$},\\
= & \inv{2}(\alpha_1+\alpha_2)^2.
\end{align*}}

When $\lambda_1<0$, $\lambda_2=0$, we have $[G_1 = \alpha_1$, $G_2 = - \alpha_2$, $\gamma_1 = 1$, $K_{2,1} = 1$, $K_{2,2} = \inv{4}$, $K_{2,3} = \frac{\log(4\pi)}{4}$, $C_{1,1}=0$, $C_{1,2} = -\frac{\log 2}{2}]$, {\adb
\begin{align*}
\notag  |B(u)| & \sim |K_{2,1}(G_1\gamma_1-G_2)|\sqrt{-2\log u} = 
\left|\alpha_1 + \alpha_2\right|\sqrt{-2\log u}, \quad \text{as $u\to 0^+$}; \\
\theta = & K_{2,1}^2\left(G_1\gamma_1- G_2\right)^2 =  (\alpha_1+\alpha_2)^2;\\
\tau_1 =&  e^{2\left[K_{2,1}^2K_{2,3}(G_1\gamma_1- G_2)^2+K_{2,1}^2C_{1,2}G_1\gamma_1(G_1\gamma_1- G_2)\right]}\\
= & e^{2\left[\frac{\log(4\pi)}{4}\left(\alpha_1+\alpha_2\right)^2 - \frac{\log 2}{2}\times \alpha_1(\alpha_1+\alpha_2)\right]}\\
= & e^{\log\left[(4\pi)^{\frac{(\alpha_1+\alpha_2)^2}{2}}\times 2^{-\alpha_1(\alpha_1+\alpha_2)}\right]}\\
=& (4\pi)^{\frac{(\alpha_1+\alpha_2)^2}{2}}\times 2^{-\alpha_1(\alpha_1+\alpha_2)};\\
\tau_2 = & 2\left[K_{2,1}^2K_{2,2}(G_1\gamma_1- G_2)^2+K_{2,1}^2C_{1,1}G_1\gamma_1(G_1\gamma_1- G_2)\right]\\
= & 2\left[\inv{4}\left(\alpha_1+\alpha_2\right)^2+0\right] = \inv{2}\left(\alpha_1+\alpha_2\right)^2, \quad \text{as $C_{1,1}=0$}.
\end{align*}}
When $\lambda_1=0$, $\lambda_2<0$, we have $[G_1 = \alpha_1$, $G_2 = - \alpha_2$, $\gamma_1 = 1$, $K_{2,1} = 1$, $K_{2,2} =\inv{4}$, $K_{2,3} = \frac{\log\pi}{4}$, $C_{1,1}=0$, $C_{1,2} = \frac{\log 2}{2}]$, {\adb
\begin{align*}
\notag  |B(u)| & \sim |K_{2,1}(G_1\gamma_1-G_2)|\sqrt{-2\log u} = 
\left|\alpha_1 + \alpha_2\right|\sqrt{-2\log u}, \quad \text{as $u\to 0^+$}; \\
\theta = & K_{2,1}^2\left(G_1\gamma_1- G_2\right)^2 =  \left(\alpha_1+\alpha_2\right)^2;\\
\tau_1 =&  e^{2\left[K_{2,1}^2K_{2,3}(G_1\gamma_1- G_2)^2+K_{2,1}^2C_{1,2}G_1\gamma_1(G_1\gamma_1- G_2)\right]}\\
= & e^{2\left[\frac{\log\pi}{4}\left(\alpha_1+\alpha_2\right)^2+\frac{\log 2}{2}\alpha_1(\alpha_1+\alpha_2)\right]}\\
=& e^{\log\left[\pi^{\frac{(\alpha_1+\alpha_2)^2}{2}}\times 2^{\alpha_1(\alpha_1+\alpha_2)}\right]}\\
=& \pi^{\frac{(\alpha_1+\alpha_2)^2}{2}}\times 2^{\alpha_1(\alpha_1+\alpha_2)};\\ 
\tau_2 = & 2\left[K_{2,1}^2K_{2,2}(G_1\gamma_1- G_2)^2+K_{2,1}^2C_{1,1}G_1\gamma_1(G_1\gamma_1- G_2)\right]\\
= & 2\left[\inv{4}(\alpha_1+\alpha_2)^2+0\right]\quad \text{as $C_{1,1}=0$},\\
= & \frac{(\alpha_1+\alpha_2)^2}{2}.
\end{align*}}
Finally, when $\lambda_1 = \lambda_2=0 \Rightarrow \alpha_1 = \alpha_2 = 0$, we have $F_1^{-1}(u) = F_2^{-1}(u) = \Phi^{-1}(u)$, and 
\begin{equation*}
|B(u)| = |\alpha_1 + \alpha_2| |\Phi^{-1}(u)| = |\alpha_1+ \alpha_2|\sqrt{-2\log u}|\to 0, \quad \text{as $u \to 0^+$}.
\end{equation*}
\end{proof}

\section*{Proof of Theorem \ref{Theorem: Summation}}

\begin{proof}
This means that when $\lambda_1<0$, $\lambda_2<0$ ($\Rightarrow A_1(u)\to -\infty$ \& $A_2(u)\to -\infty$ as well as $\alpha_1+\alpha_2<0 \Rightarrow B(u) = \alpha_1F_1^{-1}(u)+\alpha_2F_2^{-1}(u) \sim (\alpha_1+\alpha_2)F_1^{-1}(u) \to \infty$), we have {\adb
\begin{align*}
& \frac{\partial C(u,u)}{\partial u} \\
=& P\left(Z_1\leq F_1^{-1}(u)|Z_2 =F_2^{-1}(u)\right) + P\left(Z_2\leq F_2^{-1}(u)|Z_1 =F_1^{-1}(u)\right)\\
\sim & \frac{\sqrt{1-\rho^2}e^{-\frac{A_1^2(u)}{2(1-\rho^2)}}}{\sqrt{2\pi}|A_1(u)|\Phi(\lambda_2F_2^{-1}(u))} + \frac{\sqrt{1-\rho^2}e^{-\frac{A_2^2(u)}{2(1-\rho^2)}}}{\sqrt{2\pi}|A_2(u)|\Phi(\lambda_1F_1^{-1}(u))}, \quad \text{by Theorem \ref{Thm: main result}(c)};\\[0.4cm]
\sim & \frac{\pi^{\inv{2}\left(\frac{1-\rho}{1+\rho}\right)}u^{\frac{1-\rho}{1+\rho}}(-\log u)^{\inv{2}\left(\frac{1-\rho}{1+\rho}\right)}}{\sqrt{2\pi}\times 1\times \left|\frac{1-\rho}{\sqrt{1-\rho^2}}\right|\sqrt{-2\log u}\times 1} +
\frac{\pi^{\inv{2}\left(\frac{1-\rho}{1+\rho}\right)}u^{\frac{1-\rho}{1+\rho}}(-\log u)^{\inv{2}\left(\frac{1-\rho}{1+\rho}\right)}}{\sqrt{2\pi}\times 1\times \left|\frac{1-\rho}{\sqrt{1-\rho^2}}\right|\sqrt{-2\log u}\times 1}, \quad \text{by (\ref{eqn: detailed A in RV}) \&  Lemma  \ref{Lemma: exp A in RV};} \\
=& \sqrt{\frac{1+\rho}{1-\rho}}\pi^{-\frac{\rho}{1+\rho}}u^{\frac{1-\rho}{1+\rho}}(-\log u)^{-\frac{\rho}{1+\rho}}.
\end{align*}}
Here we also used the fact that $\Phi(\lambda_i F_i^{-1}(u)) \sim 1$ when $\lambda_i<0$.

When $\lambda_1=0$, $\lambda_2<0$ ($\Rightarrow A_1(u) \to -\infty$ \& $A_2(u)\to -\infty$). We still have $\alpha_1+\alpha_2<0 \Rightarrow B(u) = \alpha_1F_1^{-1}(u)+\alpha_2F_2^{-1}(u) \sim (\alpha_1+\alpha_2)F_1^{-1}(u) \to \infty$ as 
\begin{equation*}
\lambda_1= 0 \text{ \& } \lambda_2 <0  \quad \Rightarrow \quad \alpha_1+\rho\alpha_2=0 \text{ \& } \alpha_2+\rho\alpha_1<0 \quad \Rightarrow \quad (\alpha_1+\alpha_2)(1+\rho)<0\quad  \Rightarrow \quad \alpha_1+\alpha_2<0.
\end{equation*}
This means that {\adb
\begin{align*}
& \frac{dC(u,u)}{du} \\
=& P\left(Z_1\leq F_1^{-1}(u)|Z_2 =F_2^{-1}(u)\right) + P\left(Z_2\leq F_2^{-1}(u)|Z_1 =F_1^{-1}(u)\right)\\
\sim & \frac{\sqrt{1-\rho^2}e^{-\frac{A_1^2(u)}{2(1-\rho^2)}}}{\sqrt{2\pi}|A_1(u)|\Phi(\lambda_2F_2^{-1}(u))} + \frac{\sqrt{1-\rho^2}e^{-\frac{A_2^2(u)}{2(1-\rho^2)}}}{\sqrt{2\pi}|A_2(u)|\Phi(\lambda_1F_1^{-1}(u))}, \quad \text{by Theorem \ref{Thm: main result}(c)};\\[0.4cm]
\sim & \frac{\pi^{\inv{2}\left(\frac{1-\rho}{1+\rho}\right)}2^{\inv{1+\rho}}u^{\frac{1-\rho}{1+\rho}}(-\log u)^{\inv{2}\left(\frac{1-\rho}{1+\rho}\right)}}{\sqrt{2\pi}\times 1\times \left|\frac{1-\rho}{\sqrt{1-\rho^2}}\right|\sqrt{-2\log u}\times 1} +
\frac{(4\pi)^{\inv{2}\left(\frac{1-\rho}{1+\rho}\right)} 2^{-\inv{1+\rho}}u^{\frac{1-\rho}{1+\rho}}(-\log u)^{\inv{2}\left(\frac{1-\rho}{1+\rho}\right)}}{\sqrt{2\pi}\left|\frac{1-\rho}{\sqrt{1-\rho^2}}\right|\sqrt{-2\log u}\times \inv{2}}, \quad \text{by (\ref{eqn: detailed A in RV}) and Lemma \ref{Lemma: exp A in RV};}\\
=& \sqrt{\frac{1+\rho}{1-\rho}}u^{\frac{1-\rho}{1+\rho}}(-\log u)^{-\frac{\rho}{1+\rho}}\left(\pi^{\inv{2}\left(\frac{1-\rho}{1+\rho}-1\right)}2^{\inv{1+\rho}-1}+\pi^{\inv{2}\left(\frac{1-\rho}{1+\rho}-1\right)}2^{\frac{1-\rho}{1+\rho}-\inv{1+\rho}}\right)\\
= & \sqrt{\frac{1+\rho}{1-\rho}} u^{\frac{1-\rho}{1+\rho}}(-\log u)^{-\frac{\rho}{1+\rho}}\left(\pi^{-\frac{\rho}{1+\rho}}2^{-\frac{\rho}{1+\rho}}+\pi^{-\frac{\rho}{1+\rho}}2^{-\frac{\rho}{1+\rho}}\right)\\
= & \pi^{-\frac{\rho}{1+\rho}}2^{\inv{1+\rho}}\sqrt{\frac{1+\rho}{1-\rho}} u^{\frac{1-\rho}{1+\rho}}(-\log u)^{-\frac{\rho}{1+\rho}}.
\end{align*}}
Here we also used the fact that $\Phi(\lambda_i F_i^{-1}(u)) \sim 1$ when $\lambda_i<0$ and $\Phi(\lambda_i F_i^{-1}(u)) = 1/2$ when $\lambda_i=0$.

When $\lambda_1<0$, $\lambda_2>0$ ($\Rightarrow A_1(u)\to -\infty$, $A_2(u)\to -\infty$) and 
\begin{eqnarray*}
B(u) &=&  (\alpha_1F_1^{-1}(u)+\alpha_2F_2^{-1}(u))\\
& \sim & (\alpha_1+\frac{\alpha_2}{\sqrt{1+\lambda_2^2}})F_1^{-1}(u) \\
& \to & \begin{cases}
-\infty, & \text{if $\alpha_1+\frac{\alpha_2}{\sqrt{1+\lambda_2^2}}>0$;}\\
0, & \text{if $\alpha_1+\frac{\alpha_2}{\sqrt{1+\lambda_2^2}}=0$;}\\
\infty, & \text{if $\alpha_1+\frac{\alpha_2}{\sqrt{1+\lambda_2^2}}<0$.}
\end{cases}
\end{eqnarray*}
This means that if $\alpha_1+\frac{\alpha_2}{\sqrt{1+\lambda_2^2}}>0$, we have $B(u)\to -\infty$ and both $\beta_1$, $\beta_2>0$ by Lemma \ref{Lemma: beta1>0}. As a result,
{\adb
\begin{align*}
& \frac{dC(u,u)}{du} \\
=& P\left(Z_1\leq F_1^{-1}(u)|Z_2 =F_2^{-1}(u)\right) + P\left(Z_2\leq F_2^{-1}(u)|Z_1 =F_1^{-1}(u)\right)\\
\sim & \frac{e^{-\frac{A_1^2(u)}{2(1-\rho^2)}-\inv{2}B^2(u)}}{2\pi\sqrt{1-\rho^2}|F_2^{-1}(u)|\Phi(\lambda_2F_2^{-1}(u))|B(u)|\beta_1}\\
& \quad + \frac{e^{-\frac{A_2^2(u)}{2(1-\rho^2)}-\inv{2}B^2(u)}}{2\pi\sqrt{1-\rho^2}|F_1^{-1}(u)|\Phi(\lambda_1F_1^{-1}(u))|B(u)|\beta_2} , \quad \text{by Theorem \ref{Thm: main result}(a)};\\
\sim & \frac{e^{-\frac{A_1^2(u)}{2(1-\rho^2)}-\inv{2}B^2(u)}}{2\pi\sqrt{1-\rho^2}|F_2^{-1}(u)|\Phi(\lambda_2F_2^{-1}(u)) |B(u)|\beta_1}\\
& \quad + \frac{e^{-\frac{A_2^2(u)}{2(1-\rho^2)}-\inv{2}B^2(u)}}{2\pi\sqrt{1-\rho^2}|F_1^{-1}(u)|\times 1\times |B(u)|\beta_2}, \quad \text{as $\Phi(\lambda_1F_1^{-1}(u)) \sim 1$ with $\lambda_1<0$};\\
=& \frac{|\lambda_2|e^{-\frac{A_1^{2}(u)}{2(1-\rho^2)}-\inv{2}B^2(u)}}{{2\pi\sqrt{1-\rho^2}|\lambda_2F_2^{-1}(u)|\Phi(\lambda_2F_2^{-1}(u))|B(u)|\beta_1}}+ \frac{e^{-\frac{A_2^2(u)}{2(1-\rho^2)}-\inv{2}B^2(u)}}{2\pi\sqrt{1-\rho^2}|F_1^{-1}(u)B(u)|\beta_2}\\
%\frac{|\lambda_2|e^{-\frac{A_1^{2}(u)}{2(1-\rho^2)}-\inv{2}B^2(u)+\inv{2}\lambda_2^2(F_2^{-1}(u))^2}}{{\sqrt{2\pi(1-\rho^2)}|B(u)|\beta_1}}\\
\sim & \frac{\lambda_2}{2\pi\sqrt{1-\rho^2}|\alpha_1+\frac{\alpha_2}{\sqrt{1+\lambda_2^2}}|\sqrt{-2\log u}\times \beta_1}\times (2\pi\lambda_2)^{\inv{1-\rho^2}\left(1-\frac{\rho}{\sqrt{1+\lambda_2^2}}\right)^2}\left(2\lambda_2\sqrt{\pi}\right)^{-\inv{1-\rho^2}\left(1-\frac{\rho}{\sqrt{1+\lambda_2^2}}\right)}\\
& \quad \times u^{\inv{1-\rho^2}\left(1-\frac{\rho}{\sqrt{1+\lambda_2^2}}\right)^2}(-\log u)^{\inv{1-\rho^2}\left(1-\frac{\rho}{\sqrt{1+\lambda_2^2}}\right)^2-\inv{2(1-\rho^2)}\left(1-\frac{\rho}{\sqrt{1+\lambda_2^2}}\right)}\\
& \quad \times (2\pi\lambda_2)^{\left(\alpha_1+\frac{\alpha_2}{\sqrt{1+\lambda_2^2}}\right)^2}\left(2\lambda_2\sqrt{\pi}\right)^{-\alpha_1\left(\alpha_1+\frac{\alpha_2}{\sqrt{1+\lambda_2^2}}\right)}\\
& \quad \times u^{\left(\alpha_1+\frac{\alpha_2}{\sqrt{1+\lambda_2^2}}\right)^2}(-\log u)^{\left(\alpha_1+\frac{\alpha_2}{\sqrt{1+\lambda_2^2}}\right)^2-\inv{2}\alpha_1\left(\alpha_1+\frac{\alpha_2}{\sqrt{1+\lambda_2^2}}\right)} \\
& \quad \times \sqrt{2\pi}u^{-\frac{\lambda_2^2}{1+\lambda_2^2}}|\log u|^{-\frac{\lambda_2^2}{1+\lambda_2^2}}(2\pi\lambda_2)^{-\frac{\lambda_2^2}{1+\lambda_2^2}}\\
& \, +\inv{2\pi\sqrt{1-\rho^2}\times \sqrt{-2\log u}\times \left|\alpha_1+\frac{\alpha_2}{\sqrt{1+\lambda_2^2}}\right|\sqrt{-2\log u}\times \beta_2}\\
& \quad \times \pi^{\inv{2(1-\rho^2)}\left(\inv{\sqrt{1+\lambda_2^2}}-\rho\right)^2}(2\lambda_2\sqrt{\pi})^{\inv{2(1-\rho^2)\sqrt{1+\lambda_2^2}}\left(\inv{\sqrt{1+\lambda_2^2}}-\rho\right)} \\
& \quad \times  u^{\inv{1-\rho^2}\left(\inv{\sqrt{1+\lambda_2^2}}-\rho\right)^2}(-\log u)^{\inv{2(1-\rho^2)}\left(\inv{\sqrt{1+\lambda_2^2}}-\rho\right)^2+\inv{2(1-\rho^2)\sqrt{1+\lambda_2^2}}\left(\inv{\sqrt{1+\lambda_2^2}}-\rho\right)}\\
& \quad \times (2\pi\lambda_2)^{\left(\alpha_1+\frac{\alpha_2}{\sqrt{1+\lambda_2^2}}\right)^2}\left(2\lambda_2\sqrt{\pi}\right)^{-\alpha_1\left(\alpha_1+\frac{\alpha_2}{\sqrt{1+\lambda_2^2}}\right)}\\
& \quad \times u^{\left(\alpha_1+\frac{\alpha_2}{\sqrt{1+\lambda_2^2}}\right)^2}
(-\log u)^{\left(\alpha_1+\frac{\alpha_2}{\sqrt{1+\lambda_2^2}}\right)^2-\inv{2}\alpha_1\left(\alpha_1+\frac{\alpha_2}{\sqrt{1+\lambda_2^2}}\right)}, \quad \text{by using Lemmas \ref{Lemma: detailed F inverse}, \ref{Lemma: Phi in RV}--\ref{Lemma: exp B in RV} \& (\ref{eqn: detailed B in RV})};\\
= & u^{\inv{1-\rho^2}\left(1-\frac{\rho}{\sqrt{1+\lambda_2^2}}\right)^2+\left(\alpha_1+\frac{\alpha_2}{\sqrt{1+\lambda_2^2}}\right)^2-\frac{\lambda_2^2}{1+\lambda_2^2}}\\
& \quad \times (-\log u)^{\inv{1-\rho^2}\left(1-\frac{\rho}{\sqrt{1+\lambda_2^2}}\right)^2-\inv{2(1-\rho^2)}\left(1-\frac{\rho}{\sqrt{1+\lambda_2^2}}\right)+\left(\alpha_1+\frac{\alpha_2}{\sqrt{1+\lambda_2^2}}\right)^2-\frac{\alpha_1}{2}\left(\alpha_1+\frac{\alpha_2}{\sqrt{1+\lambda_2^2}}\right)-\frac{\lambda_2^2}{1+\lambda_2^2}-\inv{2}}\\
&\quad \times \frac{\lambda_2(2\pi\lambda_2)^{\inv{1-\rho^2}\left(1-\frac{\rho}{\sqrt{1+\lambda_2^2}}\right)^2+\left(\alpha_1+\frac{\alpha_2}{\sqrt{1+\lambda_2^2}}\right)^2-\frac{\lambda_2^2}{1+\lambda_2^2}}}{2\sqrt{\pi(1-\rho^2)}\left|\alpha_1+\frac{\alpha_2}{\sqrt{1+\lambda_2^2}}\right|\beta_1\times (2\lambda_2\sqrt{\pi})^{\inv{1-\rho^2}\left(1-\frac{\rho}{\sqrt{1+\lambda_2^2}}\right)+\alpha_1\left(\alpha_1+\frac{\alpha_2}{\sqrt{1+\lambda_2^2}}\right)}}\\
& \, +
u^{\inv{1-\rho^2}\left(\inv{\sqrt{1+\lambda_2^2}}-\rho\right)^2+ \left(\alpha_1+\frac{\alpha_2}{\sqrt{1+\lambda_2^2}}\right)^2}\\
& \quad \times \left(-\log u\right)^{\inv{2(1-\rho^2)}\left(\inv{\sqrt{1+\lambda_2^2}}-\rho\right)^2+\inv{2(1-\rho^2)\sqrt{1+\lambda_2^2}}\left(\inv{\sqrt{1+\lambda_2^2}}-\rho\right)+\left(\alpha_1+\frac{\alpha_2}{\sqrt{1+\lambda_2^2}}\right)^2-\inv{2}\alpha_1\left(\alpha_1+\frac{\alpha_2}{\sqrt{1+\lambda_2^2}}\right)-1}\\
& \quad \times \frac{(2\pi\lambda_2)^{\left(\alpha_1+\frac{\alpha_2}{\sqrt{1+\lambda_2^2}}\right)^2}\left(2\lambda_2\sqrt{\pi}\right)^{\inv{2(1-\rho^2)\sqrt{1+\lambda_2^2}}\left(\inv{\sqrt{1+\lambda_2^2}}-\rho\right)-\alpha_1\left(\alpha_1+\frac{\alpha_2}{\sqrt{1+\lambda_2^2}}\right)}\times \pi^{\inv{2(1-\rho^2)}\left(\inv{\sqrt{1+\lambda_2^2}}-\rho\right)^2}
}{4\pi\sqrt{1-\rho^2}\left|\alpha_1+\frac{\alpha_2}{\sqrt{1+\lambda_2^2}}\right|\times \beta_2}\\
=& u^{\inv{1-\rho^2}\left(\inv{\sqrt{1+\lambda_2^2}}-\rho\right)^2+ \left(\alpha_1+\frac{\alpha_2}{\sqrt{1+\lambda_2^2}}\right)^2} \\
& \quad \times \left(-\log u\right)^{\inv{2(1-\rho^2)}\left(\inv{\sqrt{1+\lambda_2^2}}-\rho\right)^2+\inv{2(1-\rho^2)\sqrt{1+\lambda_2^2}}\left(\inv{\sqrt{1+\lambda_2^2}}-\rho\right)+\left(\alpha_1+\frac{\alpha_2}{\sqrt{1+\lambda_2^2}}\right)^2-\inv{2}\alpha_1\left(\alpha_1+\frac{\alpha_2}{\sqrt{1+\lambda_2^2}}\right)-1}\\
& \quad \times \Biggl[ \frac{(2\pi\lambda_2)^{\left(\alpha_1+\frac{\alpha_2}{\sqrt{1+\lambda_2^2}}\right)^2}\left(2\lambda_2\sqrt{\pi}\right)^{\inv{2(1-\rho^2)\sqrt{1+\lambda_2^2}}\left(\inv{\sqrt{1+\lambda_2^2}}-\rho\right)-\alpha_1\left(\alpha_1+\frac{\alpha_2}{\sqrt{1+\lambda_2^2}}\right)}\times \pi^{\inv{2(1-\rho^2)}\left(\inv{\sqrt{1+\lambda_2^2}}-\rho\right)^2}
}{4\pi\sqrt{1-\rho^2}\left|\alpha_1+\frac{\alpha_2}{\sqrt{1+\lambda_2^2}}\right|\times \beta_2}\\
& \quad + \frac{\lambda_2(2\pi\lambda_2)^{\inv{1-\rho^2}\left(1-\frac{\rho}{\sqrt{1+\lambda_2^2}}\right)^2+\left(\alpha_1+\frac{\alpha_2}{\sqrt{1+\lambda_2^2}}\right)^2-\frac{\lambda_2^2}{1+\lambda_2^2}}}{2\sqrt{\pi(1-\rho^2)}\left|\alpha_1+\frac{\alpha_2}{\sqrt{1+\lambda_2^2}}\right|\beta_1\times (2\lambda_2\sqrt{\pi})^{\inv{1-\rho^2}\left(1-\frac{\rho}{\sqrt{1+\lambda_2^2}}\right)+\alpha_1\left(\alpha_1+\frac{\alpha_2}{\sqrt{1+\lambda_2^2}}\right)}}\Biggr], 
%& \frac{e^{-\frac{A_2^2(u)}{2(1-\rho^2)}-\inv{2}B^2(u)}}{2\pi\sqrt{1-\rho^2}|F_1^{-1}(u)|\Phi(\lambda_1F_1^{-1}(u))|B(u)|\beta_2}
\end{align*}}
as apparently both terms have the same regularly varying index: {\adb
\begin{align}
%& -\frac{A_1^2(u)}{(1-\rho^2)}+ \lambda_2^2(F_2^{-1}(u))+ \frac{A_2^2(u)}{1-\rho^2} \\
\notag & \inv{(1-\rho^2)}\left(1-\frac{\rho}{\sqrt{1+\lambda_2^2}}\right)^2+\left(\alpha_1+\frac{\alpha_2}{\sqrt{1+\lambda_2^2}}\right)^2-\frac{\lambda_2^2}{1+\lambda_2^2} \\
\notag =& \frac{1}{1-\rho^2}\left[1-\frac{2\rho}{\sqrt{1+\lambda_2^2}}+\frac{\rho^2}{1+\lambda_2^2}-\frac{(1-\rho^2)\lambda_2^2}{1+\lambda_2^2}\right]+\left(\alpha_1+\frac{\alpha_2}{\sqrt{1+\lambda_2^2}}\right)^2\\
\notag =& \frac{1}{1-\rho^2}\left[\frac{1+\lambda_2^2+\rho^2-\lambda_2^2+\rho^2\lambda_2^2}{1+\lambda_2^2}-\frac{2\rho}{\sqrt{1+\lambda_2^2}}\right]+\left(\alpha_1+\frac{\alpha_2}{\sqrt{1+\lambda_2^2}}\right)^2\\
\notag =& \frac{1}{1-\rho^2}\left[\frac{1}{1+\lambda_2^2}-\frac{2\rho}{\sqrt{1+\lambda_2^2}}+\rho^2\right]+\left(\alpha_1+\frac{\alpha_2}{\sqrt{1+\lambda_2^2}}\right)^2\\
=& \inv{1-\rho^2}\left(\inv{\sqrt{1+\lambda_2^2}}-\rho\right)^2+ \left(\alpha_1+\frac{\alpha_2}{\sqrt{1+\lambda_2^2}}\right)^2, \label{eq: same regularly varying index}
\end{align}}
and the same power for the slowly varying function:
{\adb
\begin{align}
\notag & \inv{1-\rho^2}\left(1-\frac{\rho}{\sqrt{1+\lambda_2^2}} \right)^2-\inv{2\left(1-\rho^2\right)}\left( 1-\frac{\rho}{\sqrt{1+\lambda_2^2}}\right)\\
\notag  & \quad +\left(\alpha_1+\frac{\alpha_2}{\sqrt{1+\lambda_2^2}} \right)^2-\frac{\alpha_1}{2}\left( \alpha_1+\frac{\alpha_2}{\sqrt{1+\lambda_2^2}} \right)-\frac{\lambda_2^2}{1+\lambda_2^2}-\inv{2}\\
\notag & \quad  - \Biggl[\inv{2(1-\rho^2)}\left(\inv{\sqrt{1+\lambda_2^2}}-\rho\right)^2+\inv{2(1-\rho^2)\sqrt{1+\lambda_2^2}}\left(\inv{\sqrt{1+\lambda_2^2}}-\rho\right)\\
\notag & \quad +\left(\alpha_1+\frac{\alpha_2}{\sqrt{1+\lambda_2^2}}\right)^2-\frac{\alpha_1}{2}\left(\alpha_1+\frac{\alpha_2}{\sqrt{1+\lambda_2^2}}\right)-1\Biggr]\\
\notag =& \inv{1-\rho^2}\left(1-\frac{\rho}{\sqrt{1+\lambda_2^2}} \right)^2
- \frac{\lambda_2^2}{1+\lambda_2^2}-\inv{2(1-\rho^2)}\left(1-\frac{\rho}{\sqrt{1+\lambda_2^2}}\right)
- \inv{2(1-\rho^2)}\left(\inv{\sqrt{1+\lambda_2^2}}-\rho\right)^2\\
\notag & \quad -\inv{2\left(1-\rho^2\right)\sqrt{1+\lambda_2^2}}\left(\inv{\sqrt{1+\lambda_2^2}}-\rho\right) +\inv{2}\\
\notag =& \inv{1-\rho^2}\left(\inv{\sqrt{1+\lambda_2^2}}-\rho\right)^2-\inv{2(1-\rho^2)}\left(1-\frac{\rho}{\sqrt{1+\lambda_2^2}}\right)
- \inv{2(1-\rho^2)}\left(\inv{\sqrt{1+\lambda_2^2}}-\rho\right)^2\\
\notag & \quad -\inv{2\left(1-\rho^2\right)\sqrt{1+\lambda_2^2}}\left(\inv{\sqrt{1+\lambda_2^2}}-\rho\right) +\inv{2}, \quad \text{by (\ref{eq: same regularly varying index});}\\
\notag =& \inv{2(1-\rho^2)}\left(\inv{\sqrt{1+\lambda_2^2}}-\rho\right)^2-\inv{2(1-\rho^2)}\left(1-\frac{\rho}{\sqrt{1+\lambda_2^2}}\right)-\inv{2\left(1-\rho^2\right)\sqrt{1+\lambda_2^2}}\left(\inv{\sqrt{1+\lambda_2^2}}-\rho\right) +\inv{2}\\
\notag =& \inv{2(1-\rho^2)}\left(\inv{1+\lambda_2^2}-\frac{2\rho}{\sqrt{1+\lambda_2^2}}+\rho^2-1+\frac{\rho}{\sqrt{1+\lambda_2^2}}- \inv{1+\lambda_2^2}+\frac{\rho}{\sqrt{1+\lambda_2^2}}\right)+\inv{2}\\
= & 0. \label{eq: same slowly varying index}
\end{align}
}
This means that 
{\adb
\begin{align*}
\frac{dC(u,u)}{du} \sim & u^{\inv{1-\rho^2}\left(\inv{\sqrt{1+\lambda_2^2}}-\rho\right)^2+ \left(\alpha_1+\frac{\alpha_2}{\sqrt{1+\lambda_2^2}}\right)^2}\\
& \quad \times \left(-\log u\right)^{\inv{2(1-\rho^2)}\left(\inv{\sqrt{1+\lambda_2^2}}-\rho\right)^2+\inv{2(1-\rho^2)\sqrt{1+\lambda_2^2}}\left(\inv{\sqrt{1+\lambda_2^2}}-\rho\right)+\left(\alpha_1+\frac{\alpha_2}{\sqrt{1+\lambda_2^2}}\right)^2-\inv{2}\alpha_1\left(\alpha_1+\frac{\alpha_2}{\sqrt{1+\lambda_2^2}}\right)-1}\\
& \quad \times \frac{(2\lambda_2)^{\left(\alpha_1+\frac{\alpha_2}{\sqrt{1+\lambda_2^2}} \right)^2 -\alpha_1\left(\alpha_1+\frac{\alpha_2}{\sqrt{1+\lambda_2^2}}\right)+\inv{2(1-\rho^2)\sqrt{1+\lambda_2^2}}\left( \inv{\sqrt{1+\lambda_2^2}}-\rho \right)}}{4\sqrt{1-\rho^2}\left| \alpha_1+\frac{\alpha_2}{\sqrt{1+\lambda_2^2}} \right|}\\
& \quad \times \pi^{\left(\alpha_1+\frac{\alpha_2}{\sqrt{1+\lambda_2^2}} \right)^2 -\frac{\alpha_1}{2}\left(\alpha_1+\frac{\alpha_2}{\sqrt{1+\lambda_2^2}}\right)+\inv{4(1-\rho^2)\sqrt{1+\lambda_2^2}}\left( \inv{\sqrt{1+\lambda_2^2}}-\rho \right)+\inv{2(1-\rho^2)}\left( \inv{\sqrt{1+\lambda_2^2}}-\rho \right)^2-1}\\
& \quad \times \left[ \inv{|\beta_2|}+\frac{(2\lambda_2\sqrt{\pi})^{\inv{2(1-\rho^2)\sqrt{1+\lambda_2^2}}\left( \inv{\sqrt{1+\lambda_2^2}}-\rho \right)}}{|\beta_1|} \right],
\end{align*}
}
as $u\to 0^+$ if one  attempt to pull the common components of the constants together. 

Moving on to the case when $B(u)\to 0$ i.e. $\alpha_1+\frac{\alpha_2}{\sqrt{1+\lambda_2^2}}=0$, we have
{\adb 
\begin{align}
\notag  & \frac{dC(u,u)}{du}\\
\notag =&  P\left(Z_1\leq F_1^{-1}(u)|Z_2 =F_2^{-1}(u)\right) + P\left(Z_2\leq F_2^{-1}(u)|Z_1 =F_1^{-1}(u)\right)\\
\sim & \frac{\sqrt{1-\rho^2}e^{-\frac{A_1^2(u)}{2(1-\rho^2)}}}{2\sqrt{2\pi}|A_1(u)|\Phi(\lambda_2F_2^{-1}(u))}
+ 
\frac{\sqrt{1-\rho^2}e^{-\frac{A_2^2(u)}{2(1-\rho^2)}}}{2\sqrt{2\pi}|A_2(u)|\Phi(\lambda_1F_1^{-1}(u))}, \quad \text{by Theorem \ref{Thm: main result}(b);}
\label{partial C(u,u) when B(u) to 0}\\
\notag \sim & \frac{\sqrt{1-\rho^2}|\lambda_2F_2^{-1}(u)|e^{-\frac{A_1^2(u)}{2(1-\rho^2)}}}{2\sqrt{2\pi}|A_1(u)||\lambda_2F_2^{-1}(u)|\Phi(\lambda_2F_2^{-1}(u))}
+\frac{\sqrt{1-\rho^2}e^{-\frac{A_2^2(u)}{2(1-\rho^2)}}}{2\sqrt{2\pi}|A_2(u)|\times 1}, \quad \text{as $\Phi(\lambda_1F_1^{-1}(u))\sim 1$ when $\lambda_1<0$;}\\
\notag \sim & \frac{\sqrt{1-\rho^2}\left| \lambda_2 \right| \inv{\sqrt{1+\lambda_2^2}}\sqrt{-2\log u}\times (2\pi\lambda_2)^{\inv{1-\rho^2}\left( 1-\frac{\rho}{\sqrt{1+\lambda_2^2}} \right)^2}(2\lambda_2\sqrt{\pi})^{-\inv{1-\rho^2}\left( 1-\frac{\rho}{\sqrt{1+\lambda_2^2}} \right)}}{2\sqrt{2\pi}\left| 1-\frac{\rho}{\sqrt{1+\lambda_2^2}} \right|\times\sqrt{-2\log u}}\\
\notag & \quad \times u^{\inv{1-\rho^2}\left( 1-\frac{\rho}{\sqrt{1+\lambda_2^2}} \right)^2}\left( -\log u \right)^{\inv{1-\rho^2}\left( 1-\frac{\rho}{\sqrt{1+\lambda_2^2}} \right)^2-\inv{2(1-\rho^2)}\left( 1-\frac{\rho}{\sqrt{1+\lambda_2^2}} \right)}\\
\notag & \quad \times \sqrt{2\pi}u^{-\frac{\lambda_2^2}{1+\lambda_2^2}}\left( -\log u \right)^{-\frac{\lambda_2^2}{1+\lambda_2^2}}(2\pi\lambda_2)^{-\frac{\lambda_2^2}{1+\lambda_2^2}}\\
\notag & \, + \sqrt{1-\rho^2}\pi^{\inv{2(1-\rho^2)}\left( \inv{\sqrt{1+\lambda_2^2}}-\rho \right)^2}(2\lambda_2\sqrt{\pi})^{\inv{2(1-\rho^2)\sqrt{1+\lambda_2^2}}\left( \inv{\sqrt{1+\lambda_2^2}}-\rho \right)}\\
\notag & \quad \times \frac{u^{\inv{1-\rho^2}\left( \inv{\sqrt{1+\lambda_2^2}}-\rho \right)^2}(-\log u)^{\inv{2(1-\rho^2)}\left( \inv{\sqrt{1+\lambda_2^2}}-\rho \right)^2+\inv{2(1-\rho^2)\sqrt{1+\lambda_2^2}}\left( \inv{\sqrt{1+\lambda_2^2}}-\rho \right)}}{2\sqrt{2\pi}\left| \inv{\sqrt{1+\lambda_2^2}}-\rho \right|\sqrt{-2\log u}}, \\
\notag & \quad \text{by (\ref{eqn: detailed A in RV}) \& Lemmas \ref{Lemma: detailed F inverse}, \ref{Lemma: Phi in RV} \& \ref{Lemma: exp A in RV};} \\
\notag =& \frac{\sqrt{1-\rho^2}|\lambda_2|(2\pi\lambda_2)^{\inv{1-\rho^2}\left(1-\frac{\rho}{\sqrt{1+\lambda_2^2}} \right)^2-\frac{\lambda_2^2}{1+\lambda_2^2}}(2\lambda_2\sqrt{\pi})^{-\inv{1-\rho^2}\left( 1-\frac{\rho}{\sqrt{1+\lambda_2^2}} \right)}}{2\left| \sqrt{1+\lambda_2^2}-\rho \right|}\\
\notag & \quad \times u^{\inv{1-\rho^2}\left( 1-\frac{\rho}{\sqrt{1+\lambda_2^2}} \right)^2-\frac{\lambda_2^2}{1+\lambda_2^2}}(-\log u)^{\inv{1-\rho^2}\left( 1-\frac{\rho}{\sqrt{1+\lambda_2^2}} \right)^2-\inv{2(1-\rho^2)}\left( 1-\frac{\rho}{\sqrt{1+\lambda_2^2}} \right)-\frac{\lambda_2^2}{1+\lambda_2^2}}\\
\notag & \,+ \frac{\sqrt{1-\rho^2}\pi^{\inv{2(1-\rho^2)}\left( \inv{\sqrt{1+\lambda_2^2}}-\rho \right)^2-\inv{2}}(2\lambda_2\sqrt{\pi})^{\inv{2(1-\rho^2)\sqrt{1+\lambda_2^2}}\left( \inv{\sqrt{1+\lambda_2^2}}-\rho \right)}}{4\left| \inv{\sqrt{1+\lambda_2^2}}-\rho \right|}\\
\notag & \times u^{\inv{1-\rho^2}\left(\inv{\sqrt{1+\lambda_2^2}}-\rho \right)^2}(-\log u)^{\inv{2(1-\rho^2)}\left( \inv{\sqrt{1+\lambda_2^2}}-\rho \right)^2+\inv{2(1-\rho^2)\sqrt{1+\lambda_2^2}}\left( \inv{\sqrt{1+\lambda_2^2}}-\rho \right)-\inv{2}}\\
\notag =& \Biggl[\frac{\sqrt{1-\rho^2}|\lambda_2|(2\pi\lambda_2)^{\inv{1-\rho^2}\left( 1-\frac{\rho}{\sqrt{1+\lambda_2^2}} \right)^2-\frac{\lambda_2^2}{1+\lambda_2^2}}(2\lambda_2\sqrt{\pi})^{-\inv{1-\rho^2}\left( 1-\frac{\rho}{\sqrt{1+\lambda_2^2}} \right)}}{2\left| \sqrt{1+\lambda_2^2}-\rho \right|}\\
\notag & \quad +  \frac{\sqrt{1-\rho^2}\pi^{\inv{2(1-\rho^2)}\left( \inv{\sqrt{1+\lambda_2^2}}-\rho \right)^2-\inv{2}}(2\lambda_2\sqrt{\pi})^{\inv{2(1-\rho^2)\sqrt{1+\lambda_2^2}}\left( \inv{\sqrt{1+\lambda_2^2}}-\rho \right)}}{4\left| \inv{\sqrt{1+\lambda_2^2}}-\rho \right|} \Biggr]\\
\notag & \quad \times u^{\inv{1-\rho^2}\left(\inv{\sqrt{1+\lambda_2^2}}-\rho \right)^2}(-\log u)^{\inv{2(1-\rho^2)}\left( \inv{\sqrt{1+\lambda_2^2}}-\rho \right)^2+\inv{2(1-\rho^2)\sqrt{1+\lambda_2^2}}\left( \inv{\sqrt{1+\lambda_2^2}}-\rho \right)-\inv{2}},
\end{align}}
as from (\ref{eq: same regularly varying index}) we know that 
\begin{equation*}
\inv{1-\rho^2}\left( 1-\frac{\rho}{\sqrt{1+\lambda_2^2}} \right)^2-\frac{\lambda_2^2}{1+\lambda_2^2} = \inv{1-\rho^2}\left( \inv{\sqrt{1+\lambda_2^2}}-\rho \right)^2
\end{equation*}
so both terms share the same regularly varying index. The two terms also share the same power of their slowly varying functions:
{\adb
\begin{align*}
& \inv{1-\rho^2}\left( 1-\frac{\rho}{\sqrt{1+\lambda_2^2}} \right)^2-\inv{2(1-\rho^2)}\left( 1-\frac{\rho}{\sqrt{1+\lambda_2^2}} \right)-\frac{\lambda_2^2}{1+\lambda_2^2}\\
& \quad -\left[\inv{2(1-\rho^2)}\left( \inv{\sqrt{1+\lambda_2^2}}-\rho \right)^2+\inv{2(1-\rho^2)\sqrt{1+\lambda_2^2}}\left( \inv{\sqrt{1+\lambda_2^2}}-\rho \right)-\inv{2}\right] \\
=& 0, 
\end{align*}
by (\ref{eq: same slowly varying index}).
}
This means that 
{\adb
\begin{align*}
\frac{dC(u,u)}{du} & \sim \frac{\sqrt{1-\rho^2}}{4} (2\lambda_2)^{\inv{2(1-\rho^2)\sqrt{1+\lambda_2^2}}\left( \inv{\sqrt{1+\lambda_2^2}-\rho}\right)}\pi^{\inv{2(1-\rho^2)}\left( \inv{\sqrt{1+\lambda_2^2}}-\rho \right)^2-\inv{2}}\\
& \quad \times \left[ (2\lambda\pi)^{\inv{2(1-\rho^2)\sqrt{1+\lambda_2^2}}\left( \inv{\sqrt{1+\lambda_2^2}}-\rho \right)}+ \pi^{-\inv{4(1-\rho^2)\sqrt{1+\lambda_2^2}}\left( \inv{\sqrt{1+\lambda_2^2}}-\rho \right)}
 \right] \\
\notag & \quad \times u^{\inv{1-\rho^2}\left(\inv{\sqrt{1+\lambda_2^2}}-\rho \right)^2}(-\log u)^{\inv{2(1-\rho^2)}\left( \inv{\sqrt{1+\lambda_2^2}}-\rho \right)^2+\inv{2(1-\rho^2)\sqrt{1+\lambda_2^2}}\left( \inv{\sqrt{1+\lambda_2^2}}-\rho \right)-\inv{2}}
\end{align*}
}
as $u\to 0^+$ if we attempt to pull the common components in the constants together.

Finally, when $B(u) \to \infty$ i.e. $\alpha_1+\frac{\alpha_2}{\sqrt{1+\lambda_2^2}}<0$, we have {\adb
\begin{align*}
 & \frac{dC(u,u)}{du}\\
=&  P\left(Z_1\leq F_1^{-1}(u)|Z_2 =F_2^{-1}(u)\right) + P\left(Z_2\leq F_2^{-1}(u)|Z_1 =F_1^{-1}(u)\right)\\
\sim & \frac{\sqrt{1-\rho^2}e^{-\frac{A_1^2(u)}{2(1-\rho^2)}}}{\sqrt{2\pi}|A_1(u)|\Phi(\lambda_2F_2^{-1}(u))}
+ \frac{\sqrt{1-\rho^2}e^{-\frac{A_2^2(u)}{2(1-\rho^2)}}}{\sqrt{2\pi}|A_2(u)|\Phi(\lambda_1F_1^{-1}(u))}, \quad \text{by Theorem \ref{Thm: main result}(c);} \\
& \quad \text{[this is very similar to (\ref{partial C(u,u) when B(u) to 0}), except (\ref{partial C(u,u) when B(u) to 0}) has an extra $\inv{2}$ factor to both terms]}\\
\sim & \frac{\sqrt{1-\rho^2}|\lambda_2F_2^{-1}(u)|e^{-\frac{A_1^2(u)}{2(1-\rho^2)}}}{\sqrt{2\pi}|A_1(u)||\lambda_2F_2^{-1}(u)|\Phi(\lambda_2F_2^{-1}(u))}
+ \frac{\sqrt{1-\rho^2}e^{-\frac{A_2^2(u)}{2(1-\rho^2)}}}{\sqrt{2\pi}|A_2(u)|}, \quad \text{as $\Phi(\lambda_1F_1^{-1}(u))\sim1$ as $\lambda_1<0$;}\\
%\sim &  \frac{\sqrt{1-\rho^2}|\lambda_2F_2^{-1}(u)|e^{-\frac{A_1^2(u)}{2(1-\rho^2)}+\inv{2}\lambda_2^2(F_2^{-1}(u))^2}}{|A_1(u)|}+\frac{\sqrt{1-\rho^2}e^{-\frac{A_2^2(u)}{2(1-\rho^2)}}}{\sqrt{2\pi}|A_2(u)|}\\
\sim & \Biggl[\frac{\sqrt{1-\rho^2}|\lambda_2|(2\pi\lambda_2)^{\inv{1-\rho^2}\left( 1-\frac{\rho}{\sqrt{1+\lambda_2^2}} \right)^2-\frac{\lambda_2^2}{1+\lambda_2^2}}(2\lambda_2\sqrt{\pi})^{-\inv{1-\rho^2}\left( 1-\frac{\rho}{\sqrt{1+\lambda_2^2}} \right)}}{\left| \sqrt{1+\lambda_2^2}-\rho \right|}\\
& \quad +  \frac{\sqrt{1-\rho^2}\pi^{\inv{2(1-\rho^2)}\left( \inv{\sqrt{1+\lambda_2^2}}-\rho \right)^2-\inv{2}}(2\lambda_2\sqrt{\pi})^{\inv{2(1-\rho^2)\sqrt{1+\lambda_2^2}}\left( \inv{\sqrt{1+\lambda_2^2}}-\rho \right)}}{2\left| \inv{\sqrt{1+\lambda_2^2}}-\rho \right|} \Biggr]\\
& \quad \times u^{\inv{1-\rho^2}\left(\inv{\sqrt{1+\lambda_2^2}}-\rho \right)^2}(-\log u)^{\inv{2(1-\rho^2)}\left( \inv{\sqrt{1+\lambda_2^2}}-\rho \right)^2+\inv{2(1-\rho^2)\sqrt{1+\lambda_2^2}}\left( \inv{\sqrt{1+\lambda_2^2}}-\rho \right)-\inv{2}} \\
\notag & \quad \text{by (\ref{eqn: detailed A in RV}) \& Lemmas \ref{Lemma: detailed F inverse}, \ref{Lemma: Phi in RV} \& \ref{Lemma: exp A in RV};} \\
=& \frac{\sqrt{1-\rho^2}}{2} (2\lambda_2)^{\inv{2(1-\rho^2)\sqrt{1+\lambda_2^2}}\left( \inv{\sqrt{1+\lambda_2^2}-\rho}\right)}\pi^{\inv{2(1-\rho^2)}\left( \inv{\sqrt{1+\lambda_2^2}-\rho}-\rho \right)^2-\inv{2}}\\
& \quad \times \left[ (2\lambda\pi)^{\inv{2(1-\rho^2)\sqrt{1+\lambda_2^2}}\left( \inv{\sqrt{1+\lambda_2^2}}-\rho \right)}+ \pi^{-\inv{4(1-\rho^2)\sqrt{1+\lambda_2^2}}\left( \inv{\sqrt{1+\lambda_2^2}}-\rho \right)}
 \right] \\
\notag & \quad \times u^{\inv{1-\rho^2}\left(\inv{\sqrt{1+\lambda_2^2}}-\rho \right)^2}(-\log u)^{\inv{2(1-\rho^2)}\left( \inv{\sqrt{1+\lambda_2^2}}-\rho \right)^2+\inv{2(1-\rho^2)\sqrt{1+\lambda_2^2}}\left( \inv{\sqrt{1+\lambda_2^2}}-\rho \right)-\inv{2}}
\end{align*}}

%% start here
When $\lambda_1=0$, $\lambda_2>0$ ($\Rightarrow A_1(u)\to -\infty$, $A_2(u)\to -\infty$) and 
\begin{eqnarray*}
B(u) &=& \alpha_1F_1^{-1}(u)+\alpha_2F_2^{-1}(u)\\
& \sim & \left(\alpha_1+\frac{\alpha_2}{\sqrt{1+\lambda_2^2}}\right)F_1^{-1}(u) \\
& \to & \begin{cases}
-\infty, & \text{if $\alpha_1+\frac{\alpha_2}{\sqrt{1+\lambda_2^2}}>0$;}\\
0, & \text{if $\alpha_1+\frac{\alpha_2}{\sqrt{1+\lambda_2^2}}=0$;}\\
\infty, & \text{if $\alpha_1+\frac{\alpha_2}{\sqrt{1+\lambda_2^2}}<0$.}
\end{cases}
\end{eqnarray*}
As $\lambda_1=0$, $\lambda_2>0$,  we have $\beta_1>0, \beta_2>0$ by Lemma \ref{Lemma: beta1>0}. This means that 
{\adb
\begin{align*}
& \frac{dC(u,u)}{du} \\
=& P\left(Z_1\leq F_1^{-1}(u)|Z_2 =F_2^{-1}(u)\right) + P\left(Z_2\leq F_2^{-1}(u)|Z_1 =F_1^{-1}(u)\right)\\
\sim & \frac{e^{-\frac{A_1^2(u)}{2(1-\rho^2)}-\inv{2}B^2(u)}}{2\pi\sqrt{1-\rho^2}|F_2^{-1}(u)|\Phi(\lambda_2F_2^{-1}(u))|B(u)|\beta_1}\\
& \quad + \frac{e^{-\frac{A_2^2(u)}{2(1-\rho^2)}-\inv{2}B^2(u)}}{2\pi\sqrt{1-\rho^2}|F_1^{-1}(u)|\Phi(\lambda_1F_1^{-1}(u))|B(u)|\beta_2} , \quad \text{by Theorem \ref{Thm: main result}(a)};\\
\sim & \frac{\lambda_2 e^{-\frac{A_1^2(u)}{2(1-\rho^2)}-\inv{2}B^2(u)}}{2\pi\sqrt{1-\rho^2}|\lambda_2 F_2^{-1}(u)| \Phi(\lambda_2F_2^{-1}(u)) |B(u)|\beta_1}\\
& \quad + \frac{e^{-\frac{A_2^2(u)}{2(1-\rho^2)}-\inv{2}B^2(u)}}{2\pi\sqrt{1-\rho^2}|F_1^{-1}(u)|\times \inv{2}\times |B(u)|\beta_2}, \quad \text{as $\Phi(\lambda_1F_1^{-1}(u)) = \inv{2}$ with $\lambda_1=0$};\\
=& \frac{\lambda_2 e^{-\frac{A_1^2(u)}{2(1-\rho^2)}-\inv{2}B^2(u)}}{2\pi\sqrt{1-\rho^2}|\lambda_2 F_2^{-1}(u)| \Phi(\lambda_2F_2^{-1}(u)) |B(u)|\beta_1}+ \frac{e^{-\frac{A_2^2(u)}{2(1-\rho^2)}-\inv{2}B^2(u)}}{\pi\sqrt{1-\rho^2}|F_1^{-1}(u)B(u)|\beta_2}\\
%\frac{|\lambda_2|e^{-\frac{A_1^{2}(u)}{2(1-\rho^2)}-\inv{2}B^2(u)+\inv{2}\lambda_2^2(F_2^{-1}(u))^2}}{{\sqrt{2\pi(1-\rho^2)}|B(u)|\beta_1}}\\
\sim & \frac{\lambda_2}{2\pi\sqrt{1-\rho^2}|\alpha_1+\frac{\alpha_2}{\sqrt{1+\lambda_2^2}}|\sqrt{-2\log u}\times \beta_1}\times (2\pi\lambda_2)^{\inv{1-\rho^2}\left(1-\frac{\rho}{\sqrt{1+\lambda_2^2}}\right)^2}\left(\lambda_2\sqrt{\pi}\right)^{-\inv{1-\rho^2}\left(1-\frac{\rho}{\sqrt{1+\lambda_2^2}}\right)}\\
& \quad \times u^{\inv{1-\rho^2}\left(1-\frac{\rho}{\sqrt{1+\lambda_2^2}}\right)^2}(-\log u)^{\inv{1-\rho^2}\left(1-\frac{\rho}{\sqrt{1+\lambda_2^2}}\right)^2-\inv{2(1-\rho^2)}\left(1-\frac{\rho}{\sqrt{1+\lambda_2^2}}\right)}\\
& \quad \times (2\pi\lambda_2)^{\left(\alpha_1+\frac{\alpha_2}{\sqrt{1+\lambda_2^2}}\right)^2}\left(\lambda_2\sqrt{\pi}\right)^{-\alpha_1\left(\alpha_1+\frac{\alpha_2}{\sqrt{1+\lambda_2^2}}\right)}\\
& \quad \times u^{\left(\alpha_1+\frac{\alpha_2}{\sqrt{1+\lambda_2^2}}\right)^2}(-\log u)^{\left(\alpha_1+\frac{\alpha_2}{\sqrt{1+\lambda_2^2}}\right)^2-\inv{2}\alpha_1\left(\alpha_1+\frac{\alpha_2}{\sqrt{1+\lambda_2^2}}\right)} \\
& \quad \times \sqrt{2\pi} u^{-\frac{\lambda_2^2}{1+\lambda_2^2}}|\log u|^{-\frac{\lambda_2^2}{1+\lambda_2^2}}(2\pi\lambda_2)^{-\frac{\lambda_2^2}{1+\lambda_2^2}}\\
& \, +\inv{\pi\sqrt{1-\rho^2}\times \sqrt{-2\log u}\times \left|\alpha_1+\frac{\alpha_2}{\sqrt{1+\lambda_2^2}}\right|\sqrt{-2\log u}\times \beta_2}\\
& \quad \times (4\pi)^{\inv{2(1-\rho^2)}\left(\inv{\sqrt{1+\lambda_2^2}}-\rho\right)^2}(\lambda_2\sqrt{\pi})^{\inv{(1-\rho^2)\sqrt{1+\lambda_2^2}}\left(\inv{\sqrt{1+\lambda_2^2}}-\rho\right)} \\
& \quad \times  u^{\inv{1-\rho^2}\left(\inv{\sqrt{1+\lambda_2^2}}-\rho\right)^2}(-\log u)^{\inv{2(1-\rho^2)}\left(\inv{\sqrt{1+\lambda_2^2}}-\rho\right)^2+\inv{2(1-\rho^2)\sqrt{1+\lambda_2^2}}\left(\inv{\sqrt{1+\lambda_2^2}}-\rho\right)}\\
& \quad \times (2\pi\lambda_2)^{\left(\alpha_1+\frac{\alpha_2}{\sqrt{1+\lambda_2^2}}\right)^2}\left(\lambda_2\sqrt{\pi}\right)^{-\alpha_1\left(\alpha_1+\frac{\alpha_2}{\sqrt{1+\lambda_2^2}}\right)}\\
& \quad \times u^{\left(\alpha_1+\frac{\alpha_2}{\sqrt{1+\lambda_2^2}}\right)^2}
(-\log u)^{\left(\alpha_1+\frac{\alpha_2}{\sqrt{1+\lambda_2^2}}\right)^2-\inv{2}\alpha_1\left(\alpha_1+\frac{\alpha_2}{\sqrt{1+\lambda_2^2}}\right)} \\
& \quad \text{by (\ref{eqn: detailed B in RV}) \& Lemmas \ref{Lemma: detailed F inverse}, \ref{Lemma: Phi in RV}--\ref{Lemma: exp B in RV};} \\
= & u^{\inv{1-\rho^2}\left(1-\frac{\rho}{\sqrt{1+\lambda_2^2}}\right)^2+\left(\alpha_1+\frac{\alpha_2}{\sqrt{1+\lambda_2^2}}\right)^2-\frac{\lambda_2^2}{1+\lambda_2^2}}\\
& \quad \times (-\log u)^{\inv{1-\rho^2}\left(1-\frac{\rho}{\sqrt{1+\lambda_2^2}}\right)^2-\inv{2(1-\rho^2)}\left(1-\frac{\rho}{\sqrt{1+\lambda_2^2}}\right)+\left(\alpha_1+\frac{\alpha_2}{\sqrt{1+\lambda_2^2}}\right)^2-\frac{\alpha_1}{2}\left(\alpha_1+\frac{\alpha_2}{\sqrt{1+\lambda_2^2}}\right)-\frac{\lambda_2^2}{1+\lambda_2^2}-\inv{2}}\\
&\quad \times \frac{\lambda_2(2\pi\lambda_2)^{\inv{1-\rho^2}\left(1-\frac{\rho}{\sqrt{1+\lambda_2^2}}\right)^2+\left(\alpha_1+\frac{\alpha_2}{\sqrt{1+\lambda_2^2}}\right)^2-\frac{\lambda_2^2}{1+\lambda_2^2}}}{2\sqrt{\pi(1-\rho^2)}\left|\alpha_1+\frac{\alpha_2}{\sqrt{1+\lambda_2^2}}\right|\beta_1\times (\lambda_2\sqrt{\pi})^{\inv{1-\rho^2}\left(1-\frac{\rho}{\sqrt{1+\lambda_2^2}}\right)+\alpha_1\left(\alpha_1+\frac{\alpha_2}{\sqrt{1+\lambda_2^2}}\right)}}\\
& \, +
u^{\inv{1-\rho^2}\left(\inv{\sqrt{1+\lambda_2^2}}-\rho\right)^2+ \left(\alpha_1+\frac{\alpha_2}{\sqrt{1+\lambda_2^2}}\right)^2}\\
& \quad \times \left(-\log u\right)^{\inv{2(1-\rho^2)}\left(\inv{\sqrt{1+\lambda_2^2}}-\rho\right)^2+\inv{2(1-\rho^2)\sqrt{1+\lambda_2^2}}\left(\inv{\sqrt{1+\lambda_2^2}}-\rho\right)+\left(\alpha_1+\frac{\alpha_2}{\sqrt{1+\lambda_2^2}}\right)^2-\inv{2}\alpha_1\left(\alpha_1+\frac{\alpha_2}{\sqrt{1+\lambda_2^2}}\right)-1}\\
& \quad \times \frac{(2\pi\lambda_2)^{\left(\alpha_1+\frac{\alpha_2}{\sqrt{1+\lambda_2^2}}\right)^2} \left(\lambda_2\sqrt{\pi}\right)^{\inv{(1-\rho^2)\sqrt{1+\lambda_2^2}}\left(\inv{\sqrt{1+\lambda_2^2}}-\rho\right)-\alpha_1\left(\alpha_1+\frac{\alpha_2}{\sqrt{1+\lambda_2^2}}\right)}\times (4\pi)^{\inv{2(1-\rho^2)}\left(\inv{\sqrt{1+\lambda_2^2}}-\rho\right)^2}
}{2\pi\sqrt{1-\rho^2}\left|\alpha_1+\frac{\alpha_2}{\sqrt{1+\lambda_2^2}}\right|\times \beta_2}\\
=& u^{\inv{1-\rho^2}\left(\inv{\sqrt{1+\lambda_2^2}}-\rho\right)^2+ \left(\alpha_1+\frac{\alpha_2}{\sqrt{1+\lambda_2^2}}\right)^2}\\
& \quad \times \left(-\log u\right)^{\inv{2(1-\rho^2)}\left(\inv{\sqrt{1+\lambda_2^2}}-\rho\right)^2+\inv{2(1-\rho^2)\sqrt{1+\lambda_2^2}}\left(\inv{\sqrt{1+\lambda_2^2}}-\rho\right)+\left(\alpha_1+\frac{\alpha_2}{\sqrt{1+\lambda_2^2}}\right)^2-\inv{2}\alpha_1\left(\alpha_1+\frac{\alpha_2}{\sqrt{1+\lambda_2^2}}\right)-1}\\
& \quad \times \Biggl[ \frac{(2\pi\lambda_2)^{\left(\alpha_1+\frac{\alpha_2}{\sqrt{1+\lambda_2^2}}\right)^2}\left(\lambda_2 \sqrt{\pi}\right)^{\inv{(1-\rho^2)\sqrt{1+\lambda_2^2}}\left(\inv{\sqrt{1+\lambda_2^2}}-\rho\right)-\alpha_1\left(\alpha_1+\frac{\alpha_2}{\sqrt{1+\lambda_2^2}}\right)}\times (4\pi)^{\inv{2(1-\rho^2)}\left(\inv{\sqrt{1+\lambda_2^2}}-\rho\right)^2}
}{2\pi\sqrt{1-\rho^2}\left|\alpha_1+\frac{\alpha_2}{\sqrt{1+\lambda_2^2}}\right|\times \beta_2}\\
& \quad + \frac{\lambda_2(2\pi\lambda_2)^{\inv{1-\rho^2}\left(1-\frac{\rho}{\sqrt{1+\lambda_2^2}}\right)^2+\left(\alpha_1+\frac{\alpha_2}{\sqrt{1+\lambda_2^2}}\right)^2-\frac{\lambda_2^2}{1+\lambda_2^2}}}{2\sqrt{\pi(1-\rho^2)}\left|\alpha_1+\frac{\alpha_2}{\sqrt{1+\lambda_2^2}}\right|\beta_1\times (\lambda_2\sqrt{\pi})^{\inv{1-\rho^2}\left(1-\frac{\rho}{\sqrt{1+\lambda_2^2}}\right)+\alpha_1\left(\alpha_1+\frac{\alpha_2}{\sqrt{1+\lambda_2^2}}\right)}}\Biggr]
%& \frac{e^{-\frac{A_2^2(u)}{2(1-\rho^2)}-\inv{2}B^2(u)}}{2\pi\sqrt{1-\rho^2}|F_1^{-1}(u)|\Phi(\lambda_1F_1^{-1}(u))|B(u)|\beta_2}
\end{align*}}
as both terms have the same regularly varying index and share the same power for the slowly varying functions similar to the case $\lambda_1<0$, $\lambda_2>0$ with $B(u)\to -\infty$.

This means that 
{\adb
\begin{align*}
\frac{dC(u,u)}{du} \sim & 
u^{\inv{1-\rho^2}\left(\inv{\sqrt{1+\lambda_2^2}}-\rho\right)^2+ \left(\alpha_1+\frac{\alpha_2}{\sqrt{1+\lambda_2^2}}\right)^2}\\
& \quad \times \left(-\log u\right)^{\inv{2(1-\rho^2)}\left(\inv{\sqrt{1+\lambda_2^2}}-\rho\right)^2+\inv{2(1-\rho^2)\sqrt{1+\lambda_2^2}}\left(\inv{\sqrt{1+\lambda_2^2}}-\rho\right)+\left(\alpha_1+\frac{\alpha_2}{\sqrt{1+\lambda_2^2}}\right)^2-\inv{2}\alpha_1\left(\alpha_1+\frac{\alpha_2}{\sqrt{1+\lambda_2^2}}\right)-1}
\\
\quad & \times \frac{(2\pi\lambda_2)^{\left(\alpha_1+\frac{\alpha_2}{\sqrt{1+\lambda_2^2}}\right)^2}(\lambda_2\sqrt{\pi})^{\inv{1-\rho^2)\sqrt{1+\lambda_2^2}}\left( \inv{\sqrt{1+\lambda_2^2}}-\rho \right)-\alpha_1\left(\alpha_1+\frac{\alpha_2}{\sqrt{1+\lambda_2^2}}\right)}\left( 2\sqrt{\pi}\right)^{\inv{1-\rho^2}\left( \inv{\sqrt{1+\lambda_2^2}}-\rho \right)^2}}{2\sqrt{\pi(1-\rho^2)}\left| \alpha_1+\frac{\alpha_2}{\sqrt{1+\lambda_2^2}} \right|}\\
&\quad \times \left[ \inv{\sqrt{\pi}\beta_2}+\inv{\beta_1} \right],
\end{align*}
}
as $u\to 0^+$ if we attempt to pull some common components out from the constant terms.

Moving on to the case when $B(u)\to 0$ i.e. $\alpha_1+\frac{\alpha_2}{\sqrt{1+\lambda_2^2}}=0$, we have 
{\adb
\begin{align*}
 & \frac{dC(u,u)}{du}\\
=&  P\left(Z_1\leq F_1^{-1}(u)|Z_2 =F_2^{-1}(u)\right) + P\left(Z_2\leq F_2^{-1}(u)|Z_1 =F_1^{-1}(u)\right)\\
\sim & \frac{\sqrt{1-\rho^2}e^{-\frac{A_1^2(u)}{2(1-\rho^2)}}}{2\sqrt{2\pi}|A_1(u)|\Phi(\lambda_2F_2^{-1}(u))}
+ \frac{\sqrt{1-\rho^2}e^{-\frac{A_2^2(u)}{2(1-\rho^2)}}}{2\sqrt{2\pi}|A_2(u)|\Phi(\lambda_1F_1^{-1}(u))}, \quad \text{by Theorem \ref{Thm: main result}(b);} \\
\sim & \frac{\sqrt{1-\rho^2}|\lambda_2F_2^{-1}(u)|e^{-\frac{A_1^{2}(u)}{2(1-\rho^2)}}}{2\sqrt{2\pi}|A_1(u)|\lambda_2F_2^{-1}(u)|\Phi(\lambda_2F_2^{-1}(u))} + \frac{\sqrt{1-\rho^2}e^{-\frac{A_2^{2}(u)}{2(1-\rho^2)}}}{2\sqrt{2\pi}|A_2(u)|\times \inv{2}}, \quad \text{as $\lambda_1 =0 \Rightarrow \Phi(\lambda_1F_1^{-1}(u)) = \inv{2}$;}\\
\sim &  \frac{\sqrt{1-\rho^2}|\lambda_2F_2^{-1}(u)|e^{-\frac{A_1^{2}(u)}{2(1-\rho^2)}}}{2\sqrt{2\pi}|A_1(u)|\lambda_2F_2^{-1}(u)|\Phi(\lambda_2F_2^{-1}(u))} +\frac{\sqrt{1-\rho^2}e^{-\frac{A_2^2(u)}{2(1-\rho^2)}}}{\sqrt{2\pi}|A_2(u)|}\\
\sim & \frac{\sqrt{1-\rho^2}\left| \lambda_2 \right| \inv{\sqrt{1+\lambda_2^2}}\sqrt{-2\log u}\times (2\pi\lambda_2)^{\inv{1-\rho^2}\left( 1-\frac{\rho}{\sqrt{1+\lambda_2^2}} \right)^2}(\lambda_2\sqrt{\pi})^{-\inv{1-\rho^2}\left( 1-\frac{\rho}{\sqrt{1+\lambda_2^2}} \right)}}{2\sqrt{2\pi}\left| 1-\frac{\rho}{\sqrt{1+\lambda_2^2}} \right|\times\sqrt{-2\log u}}\\
& \quad \times u^{\inv{1-\rho^2}\left( 1-\frac{\rho}{\sqrt{1+\lambda_2^2}} \right)^2}\left( -\log u \right)^{\inv{1-\rho^2}\left( 1-\frac{\rho}{\sqrt{1+\lambda_2^2}} \right)^2-\inv{2(1-\rho^2)}\left( 1-\frac{\rho}{\sqrt{1+\lambda_2^2}} \right)}\\
& \quad \times \sqrt{2\pi} u^{-\frac{\lambda_2^2}{1+\lambda_2^2}}\left( -\log u \right)^{-\frac{\lambda_2^2}{1+\lambda_2^2}}(2\pi\lambda_2)^{-\frac{\lambda_2^2}{1+\lambda_2^2}}\\
& \, + \sqrt{1-\rho^2}(4\pi)^{\inv{2(1-\rho^2)}\left( \inv{\sqrt{1+\lambda_2^2}}-\rho \right)^2}(\lambda_2\sqrt{\pi})^{\inv{(1-\rho^2)\sqrt{1+\lambda_2^2}}\left( \inv{\sqrt{1+\lambda_2^2}}-\rho \right)}\\
& \quad \times \frac{u^{\inv{1-\rho^2}\left( \inv{\sqrt{1+\lambda_2^2}}-\rho \right)^2}(-\log u)^{\inv{2(1-\rho^2)}\left( \inv{\sqrt{1+\lambda_2^2}}-\rho \right)^2+\inv{2(1-\rho^2)\sqrt{1+\lambda_2^2}}\left( \inv{\sqrt{1+\lambda_2^2}}-\rho \right)}}{\sqrt{2\pi}\left| \inv{\sqrt{1+\lambda_2^2}}-\rho \right|\sqrt{-2\log u}}\\
& \quad \text{by (\ref{eqn: detailed A in RV}) \& Lemmas \ref{Lemma: detailed F inverse}, \ref{Lemma: Phi in RV} \& \ref{Lemma: exp A in RV};} \\
=& \frac{\sqrt{1-\rho^2}|\lambda_2|(2\pi\lambda_2)^{\inv{1-\rho^2}\left( 1-\frac{\rho}{\sqrt{1+\lambda_2^2}} \right)^2-\frac{\lambda_2^2}{1+\lambda_2^2}}(\lambda_2\sqrt{\pi})^{-\inv{1-\rho^2}\left( 1-\frac{\rho}{\sqrt{1+\lambda_2^2}} \right)}}{2\left| \sqrt{1+\lambda_2^2}-\rho \right|}\\
& \quad \times u^{\inv{1-\rho^2}\left( 1-\frac{\rho}{\sqrt{1+\lambda_2^2}} \right)^2-\frac{\lambda_2^2}{1+\lambda_2^2}}(-\log u)^{\inv{1-\rho^2}\left( 1-\frac{\rho}{\sqrt{1+\lambda_2^2}} \right)^2-\inv{2(1-\rho^2)}\left( 1-\frac{\rho}{\sqrt{1+\lambda_2^2}} \right)-\frac{\lambda_2^2}{1+\lambda_2^2}}\\
& \,+ \frac{\sqrt{1-\rho^2}(4\pi)^{\inv{2(1-\rho^2)}\left( \inv{\sqrt{1+\lambda_2^2}}-\rho \right)^2}(\lambda_2\sqrt{\pi})^{\inv{(1-\rho^2)\sqrt{1+\lambda_2^2}}\left( \inv{\sqrt{1+\lambda_2^2}}-\rho \right)}}{2\sqrt{\pi}\left| \inv{\sqrt{1+\lambda_2^2}}-\rho \right|}\\
& \times u^{\inv{1-\rho^2}\left(\inv{\sqrt{1+\lambda_2^2}}-\rho \right)^2}(-\log u)^{\inv{2(1-\rho^2)}\left( \inv{\sqrt{1+\lambda_2^2}}-\rho \right)^2+\inv{2(1-\rho^2)\sqrt{1+\lambda_2^2}}\left( \inv{\sqrt{1+\lambda_2^2}}-\rho \right)-\inv{2}}\\
=& \Biggl[\frac{\sqrt{1-\rho^2}|\lambda_2|(2\pi\lambda_2)^{\inv{1-\rho^2}\left( 1-\frac{\rho}{\sqrt{1+\lambda_2^2}} \right)^2-\frac{\lambda_2^2}{1+\lambda_2^2}}(\lambda_2\sqrt{\pi})^{-\inv{1-\rho^2}\left( 1-\frac{\rho}{\sqrt{1+\lambda_2^2}} \right)}}{2\left| \sqrt{1+\lambda_2^2}-\rho \right|}\\
& \quad +  \frac{\sqrt{1-\rho^2}(4\pi)^{\inv{2(1-\rho^2)}\left( \inv{\sqrt{1+\lambda_2^2}}-\rho \right)^2}(\lambda_2\sqrt{\pi})^{\inv{(1-\rho^2)\sqrt{1+\lambda_2^2}}\left( \inv{\sqrt{1+\lambda_2^2}}-\rho \right)}}{2\sqrt{\pi}\left| \inv{\sqrt{1+\lambda_2^2}}-\rho \right|} \Biggr]\\
& \quad \times u^{\inv{1-\rho^2}\left(\inv{\sqrt{1+\lambda_2^2}}-\rho \right)^2}(-\log u)^{\inv{2(1-\rho^2)}\left( \inv{\sqrt{1+\lambda_2^2}}-\rho \right)^2+\inv{2(1-\rho^2)\sqrt{1+\lambda_2^2}}\left( \inv{\sqrt{1+\lambda_2^2}}-\rho \right)-\inv{2}},
\end{align*}}
as both terms share the same regularly varying index and slowly varying functions (similar to the case $\lambda_1<0$, $\lambda_2>0$ with $B(u)\to 0$.) This means that 
{\adb
\begin{align*}
\frac{dC(u,u)}{du} = & 
\frac{\sqrt{1-\rho^2}}{2\sqrt{\pi}}\left(\lambda_2\sqrt{\pi}\right)^{\inv{(1-\rho^2)\sqrt{1+\lambda_2^2}}\left( \inv{\sqrt{1+\lambda_2^2}}-\rho \right)}(2\sqrt{\pi})^{\inv{1-\rho^2}\left( \inv{\sqrt{1+\lambda_2^2}}-\rho \right)^2}\\
& \quad \times  \left[ \inv{\sqrt{1+\lambda_2^2}-\rho}+\inv{\left| \inv{\sqrt{1+\lambda_2^2}}-\rho \right|}
 \right]\\
 & \quad \times u^{\inv{1-\rho^2}\left(\inv{\sqrt{1+\lambda_2^2}}-\rho \right)^2}(-\log u)^{\inv{2(1-\rho^2)}\left( \inv{\sqrt{1+\lambda_2^2}}-\rho \right)^2+\inv{2(1-\rho^2)\sqrt{1+\lambda_2^2}}\left( \inv{\sqrt{1+\lambda_2^2}}-\rho \right)-\inv{2}}
\end{align*}
}
as $u\to 0^+$ if we attempt to pull some common components of the constants together.

Finally, when $B(u) \to \infty$ i.e. $\alpha_1+\frac{\alpha_2}{\sqrt{1+\lambda_2^2}}<0$, we have {\adb
\begin{align*}
 & \frac{d C(u,u)}{d u}\\
=&  P\left(Z_1\leq F_1^{-1}(u)|Z_2 =F_2^{-1}(u)\right) + P\left(Z_2\leq F_2^{-1}(u)|Z_1 =F_1^{-1}(u)\right)\\
\sim & \frac{\sqrt{1-\rho^2}e^{-\frac{A_1^2(u)}{2(1-\rho^2)}}}{\sqrt{2\pi}|A_1(u)|\Phi(\lambda_2F_2^{-1}(u))}
+ \frac{\sqrt{1-\rho^2}e^{-\frac{A_2^2(u)}{2(1-\rho^2)}}}{\sqrt{2\pi}|A_2(u)|\Phi(\lambda_1F_1^{-1}(u))}, \quad \text{by Theorem \ref{Thm: main result}(c);} \\
\sim &  \frac{\sqrt{1-\rho^2}|\lambda_2F_2^{-1}(u)|e^{-\frac{A_1^2(u)}{2(1-\rho^2)}}}{\sqrt{2\pi}|A_1(u)||\lambda_2F_2^{-1}(u)|\Phi(\lambda_2F_2^{-1}(u))}+\frac{\sqrt{1-\rho^2}e^{-\frac{A_2^2(u)}{2(1-\rho^2)}}}{\sqrt{2\pi}|A_2(u)|\times\inv{2}}, \quad \text{as $\lambda_1 = 0\Rightarrow \Phi(\lambda_1F_1^{-1}(u))=\inv{2}$};\\
\sim & \frac{\sqrt{1-\rho^2}|\lambda_2| \inv{\sqrt{1+\lambda_2^2}}\sqrt{-2\log u} \left( 2\pi\lambda_2 \right)^{\inv{1-\rho^2}\left( 1-\frac{\rho}{\sqrt{1+\lambda_2^2}} \right)^2}
\left( \lambda_2\sqrt{\pi} \right)^{-\inv{1-\rho^2}\left( 1-\frac{\rho}{\sqrt{1+\lambda_2^2}} \right)}
}{\sqrt{2\pi}\left| 1-\frac{\rho}{\sqrt{1+\lambda_2^2}} \right|\sqrt{-2\log u}}\\
& \quad \times u^{\inv{1-\rho^2}\left( 1-\frac{\rho}{\sqrt{1+\lambda_2^2}} \right)^2}
\left(-\log u \right)^{\inv{1-\rho^2}\left( 1-\frac{\rho}{\sqrt{1+\lambda_2^2}}\right)^2-\inv{2(1-\rho^2)}
\left( 1-\frac{\rho}{\sqrt{1+\lambda_2^2}} \right)}\\
& \quad \times \sqrt{2\pi}u^{-\frac{\lambda_2^2}{1+\lambda_2^2}}(-\log u)^{-\frac{\lambda_2^2}{1+\lambda_2^2}}(2\pi\lambda_2)^{-\frac{\lambda_2^2}{1+\lambda_2^2}}\\
& \quad + \sqrt{1-\rho^2}\left( 4\pi \right)^{\inv{2(1-\rho^2)}\left( \inv{\sqrt{1+\lambda_2^2}}-\rho \right)^2}\left(\lambda_2\sqrt{\pi} \right)^{\inv{(1-\rho^2)\sqrt{1-\lambda_2^2}}\left( \inv{\sqrt{1+\lambda_2^2}}-\rho \right)}\\
& \quad \times \frac{u^{\inv{1-\rho^2}\left( \inv{\sqrt{1+\lambda_2^2}}-\rho \right)^2}\left(-\log u \right)^{\inv{2(1-\rho^2)}\left( \inv{\sqrt{1+\lambda_2^2}}-\rho \right)^2+\inv{2(1-\rho^2)\sqrt{1+\lambda_2^2}}\left( \inv{\sqrt{1+\lambda_2^2}}-\rho \right)}}{\sqrt{\pi/2}\left| \inv{\sqrt{1+\lambda_2^2}}-\rho \right|\sqrt{-2\log u}}\\
& \quad \text{by (\ref{eqn: detailed A in RV}) \& Lemmas \ref{Lemma: detailed F inverse}, \ref{Lemma: Phi in RV} \& \ref{Lemma: exp A in RV};} \\
\sim & \Biggl[\frac{\sqrt{1-\rho^2}|\lambda_2|(2\pi\lambda_2)^{\inv{1-\rho^2}\left( 1-\frac{\rho}{\sqrt{1+\lambda_2^2}} \right)^2-\frac{\lambda_2^2}{1+\lambda_2^2}}(\lambda_2\sqrt{\pi})^{-\inv{1-\rho^2}\left( 1-\frac{\rho}{\sqrt{1+\lambda_2^2}} \right)}}{\left| \sqrt{1+\lambda_2^2}-\rho \right|}\\
& \quad +  \frac{\sqrt{1-\rho^2}(4\pi)^{\inv{2(1-\rho^2)}\left( \inv{\sqrt{1+\lambda_2^2}}-\rho \right)^2}(\lambda_2\sqrt{\pi})^{\inv{(1-\rho^2)\sqrt{1+\lambda_2^2}}\left( \inv{\sqrt{1+\lambda_2^2}}-\rho \right)}}{\sqrt{\pi}\left| \inv{\sqrt{1+\lambda_2^2}}-\rho \right|} \Biggr]\\
& \quad \times u^{\inv{1-\rho^2}\left(\inv{\sqrt{1+\lambda_2^2}}-\rho \right)^2}(-\log u)^{\inv{2(1-\rho^2)}\left( \inv{\sqrt{1+\lambda_2^2}}-\rho \right)^2+\inv{2(1-\rho^2)\sqrt{1+\lambda_2^2}}\left( \inv{\sqrt{1+\lambda_2^2}}-\rho \right)-\inv{2}}
\end{align*}}
as both terms are very similar to those when $B(u)\to 0$ except for some small changes to the constant. This means that 
{\adb
\begin{align*}
\frac{d C(u,u)}{d u} = & 
\frac{\sqrt{1-\rho^2}}{\sqrt{\pi}}\left(\lambda_2\sqrt{\pi}\right)^{\inv{(1-\rho^2)\sqrt{1+\lambda_2^2}}\left( \inv{\sqrt{1+\lambda_2^2}}-\rho \right)}(2\sqrt{\pi})^{\inv{1-\rho^2}\left( \inv{\sqrt{1+\lambda_2^2}}-\rho \right)^2}\\
& \quad \times  \left[ \inv{\sqrt{1+\lambda_2^2}-\rho}+\inv{\left| \inv{\sqrt{1+\lambda_2^2}}-\rho \right|}
 \right]\\
 & \quad \times u^{\inv{1-\rho^2}\left(\inv{\sqrt{1+\lambda_2^2}}-\rho \right)^2}(-\log u)^{\inv{2(1-\rho^2)}\left( \inv{\sqrt{1+\lambda_2^2}}-\rho \right)^2+\inv{2(1-\rho^2)\sqrt{1+\lambda_2^2}}\left( \inv{\sqrt{1+\lambda_2^2}}-\rho \right)-\inv{2}}
\end{align*}
}
as $u\to 0^+$ if we attempt to pull some common components of the constants together.

Finally, when $\lambda_1$, $\lambda_2>0$ ($\Rightarrow$ $\alpha_1+\alpha_2>0$; at least one of $\alpha_1$, $\alpha_2>0$; $B(u) \to -\infty$), we have $\beta_1>0$, $\beta_2>0$ by Lemma \ref{Lemma: beta1>0}. This means that 
{\adb
\begin{align*}
& \frac{d C(u,u)}{d u}\\
= & P\left( Z_1\leq F_1^{-1}(u)|Z_2 = F_2^{-1}(u) \right)+ P\left( Z_2\leq F_2^{-1}(u) | Z_1 = F_1^{-1}(u) \right)\\
\sim & \frac{e^{-\frac{A_1^{2}(u)}{2(1-\rho^2)}-\inv{2}B^2(u)}}{2\pi\sqrt{1-\rho^2}|F_2^{-1}(u)|\Phi\left(\lambda_2F_2^{-1}(u)\right)|B(u)|\beta_1}+\\
& \quad \frac{e^{-\frac{A_2^{2}(u)}{2(1-\rho^2)}-\inv{2}B^2(u)}}{2\pi\sqrt{1-\rho^2}|F_1^{-1}(u)|\Phi\left(\lambda_1F_1^{-1}(u)\right)|B(u)|\beta_2}, \quad \text{by Theorem \ref{Thm: main result}(a);}\\
\sim & \frac{\lambda_2 e^{-\frac{A_1^{2}(u)}{2(1-\rho^2)}-\inv{2}B^2(u)}}{2\pi\sqrt{1-\rho^2}|\lambda_2F_2^{-1}(u)| \Phi(\lambda_2F_2^{-1}(u))| |B(u)|\beta_1}+
\frac{\lambda_1 e^{-\frac{A_2^{2}(u)}{2(1-\rho^2)}-\inv{2}B^2(u)}}{2\pi \sqrt{1-\rho^2}|\lambda_1F_1^{-1}(u)|\Phi(\lambda_1F_1^{-1}(u)) |B(u)|\beta_2}\\
\sim & \lambda_2\left( 2\pi\lambda_2 \right)^{\inv{1-\rho^2}\left( \inv{\sqrt{1+\lambda_1^2}}-\frac{\rho}{\sqrt{1+\lambda_2^2}} \right)^2}\left( \frac{\lambda_1}{\lambda_2} \right)^{\inv{(1-\rho^2)\sqrt{1+\lambda_1^2}}\left( \inv{\sqrt{1+\lambda_1^2}}-\frac{\rho}{\sqrt{1+\lambda_2^2}} \right)}\\
& \quad \times  u^{\inv{1-\rho^2}\left( \inv{\sqrt{1+\lambda_1^2}}-\frac{\rho}{\sqrt{1+\lambda_2^2}} \right)^2}(-\log u)^{\inv{1-\rho^2}\left(\inv{\sqrt{1+\lambda_1^2}}-\frac{\rho}{\sqrt{1+\lambda_2^2}} \right)^2}\\
& \quad \times \left( 2\pi\lambda_2\right)^{\left( \frac{\alpha_1}{\sqrt{1+\lambda_1^2}}+\frac{\alpha_2}{\sqrt{1+\lambda_2^2}} \right)^2} \left( \frac{\lambda_1}{\lambda_2} \right)^{\frac{\alpha_1}{\sqrt{1+\lambda_1^2}}\left( \frac{\alpha_1}{\sqrt{1+\lambda_1^2}}+\frac{\alpha_2}{\sqrt{1+\lambda_2^2}} \right)}\\
& \quad \times u^{\left( \frac{\alpha_1}{\sqrt{1+\lambda_1^2}}+\frac{\alpha_2}{\sqrt{1+\lambda_2^2}} \right)^2}(-\log u)^{\left( \frac{\alpha_1}{\sqrt{1+\lambda_1^2}}+\frac{\alpha_2}{\sqrt{1+\lambda_2^2}}\right)^2}\\
& \quad \times \frac{\sqrt{2\pi} u^{-\frac{\lambda_2^2}{1+\lambda_2^2}}(-\log u)^{-\frac{\lambda_2^2}{1+\lambda_2^2}}(2\pi\lambda_2)^{-\frac{\lambda_2^2}{1+\lambda_2^2}}}{2\pi \sqrt{1-\rho^2}\left| \frac{\alpha_1}{\sqrt{1+\lambda_1^2}}+\frac{\alpha_2}{\sqrt{1+\lambda_2^2}} \right|\sqrt{-2\log u} \times \beta_1}\\
 & +\lambda_1\left( 2\pi\lambda_1 \right)^{\inv{1-\rho^2}\left( \inv{\sqrt{1+\lambda_2^2}}-\frac{\rho}{\sqrt{1+\lambda_1^2}} \right)^2}\left( \frac{\lambda_2}{\lambda_1} \right)^{\inv{(1-\rho^2)\sqrt{1+\lambda_2^2}}\left( \inv{\sqrt{1+\lambda_2^2}}-\frac{\rho}{\sqrt{1+\lambda_1^2}} \right)}\\
& \quad \times  u^{\inv{1-\rho^2}\left( \inv{\sqrt{1+\lambda_2^2}}-\frac{\rho}{\sqrt{1+\lambda_1^2}} \right)^2}(-\log u)^{\inv{1-\rho^2}\left(\inv{\sqrt{1+\lambda_2^2}}-\frac{\rho}{\sqrt{1+\lambda_1^2}} \right)^2}\\
& \quad \times \left( 2\pi\lambda_2\right)^{\left( \frac{\alpha_1}{\sqrt{1+\lambda_1^2}}+\frac{\alpha_2}{\sqrt{1+\lambda_2^2}} \right)^2} \left( \frac{\lambda_1}{\lambda_2} \right)^{\frac{\alpha_1}{\sqrt{1+\lambda_1^2}}\left( \frac{\alpha_1}{\sqrt{1+\lambda_1^2}}+\frac{\alpha_2}{\sqrt{1+\lambda_2^2}} \right)}\\
& \quad \times u^{\left( \frac{\alpha_1}{\sqrt{1+\lambda_1^2}}+\frac{\alpha_2}{\sqrt{1+\lambda_2^2}} \right)^2}(-\log u)^{\left( \frac{\alpha_1}{\sqrt{1+\lambda_1^2}}+\frac{\alpha_2}{\sqrt{1+\lambda_2^2}}\right)^2}\\
& \quad \times \frac{\sqrt{2\pi} u^{-\frac{\lambda_1^2}{1+\lambda_1^2}}(-\log u)^{-\frac{\lambda_1^2}{1+\lambda_1^2}}(2\pi\lambda_1)^{-\frac{\lambda_1^2}{1+\lambda_1^2}}}{2\pi \sqrt{1-\rho^2}\left| \frac{\alpha_1}{\sqrt{1+\lambda_1^2}}+\frac{\alpha_2}{\sqrt{1+\lambda_2^2}} \right|\sqrt{-2\log u} \times \beta_2}, \quad \text{by (\ref{eqn: detailed B in RV}) \& Lemmas \ref{Lemma: Phi in RV}--\ref{Lemma: exp B in RV};} \\
= & \frac{\lambda_2 (2\pi)^{\inv{1-\rho^2}\left( \inv{\sqrt{1+\lambda_1^2}}-\frac{\rho}{\sqrt{1+\lambda_2^2}} \right)^2+\left( \frac{\alpha_1}{\sqrt{1+\lambda_1^2}}+\frac{\alpha_2}{\sqrt{1+\lambda_2^2}} \right)^2-\frac{\lambda_2^2}{1+\lambda_2^2}-\inv{2}}}{\sqrt{2(1-\rho^2)}\left| \frac{\alpha_1}{\sqrt{1+\lambda_1^2}}+\frac{\alpha_2}{\sqrt{1+\lambda_2^2}} \right|\times \beta_1}\\
& \quad \times \lambda_2^{\inv{1-\rho^2}\left( \inv{\sqrt{1+\lambda_1^2}}-\frac{\rho}{\sqrt{1+\lambda_2^2}}\right)^2+\left(  \frac{\alpha_1}{\sqrt{1+\lambda_1^2}}+\frac{\alpha_2}{\sqrt{1+\lambda_2^2}}\right)^2-\frac{\lambda_2^2}{1+\lambda_2^2}}\\
& \quad \times \left( \frac{\lambda_1}{\lambda_2} \right)^{\inv{(1-\rho^2)\sqrt{1+\lambda_1^2}}\left( \inv{\sqrt{1+\lambda_1^2}}-\frac{\rho}{\sqrt{1+\lambda_2^2}} \right)+\frac{\alpha_1}{\sqrt{1+\lambda_1^2}}\left( \frac{\alpha_1}{\sqrt{1+\lambda_1^2}}+\frac{\alpha_2}{\sqrt{1+\lambda_2^2}} \right)}\\
& \quad \times u^{\inv{1-\rho^2}\left( \inv{\sqrt{1+\lambda_1^2}}-\frac{\rho}{\sqrt{1+\lambda_2^2}} \right)^2+\left( \frac{\alpha_1}{\sqrt{1+\lambda_1^2}}+\frac{\alpha_2}{\sqrt{1+\lambda_2^2}} \right)^2-\frac{\lambda_2^2}{1+\lambda_2^2}}\\
& \quad \times (-\log u)^{\inv{1-\rho^2}\left( \inv{\sqrt{1+\lambda_1^2}}-\frac{\rho}{\sqrt{1+\lambda_2^2}} \right)^2+\left( \frac{\alpha_1}{\sqrt{1+\lambda_1^2}}+\frac{\alpha_2}{\sqrt{1+\lambda_2^2}} \right)^2-\frac{\lambda_2^2}{1+\lambda_2^2}-\inv{2}}\\
& + \frac{\lambda_1 (2\pi)^{\inv{1-\rho^2}\left( \inv{\sqrt{1+\lambda_2^2}}-\frac{\rho}{\sqrt{1+\lambda_1^2}} \right)^2+\left( \frac{\alpha_1}{\sqrt{1+\lambda_1^2}}+\frac{\alpha_2}{\sqrt{1+\lambda_2^2}} \right)^2-\frac{\lambda_1^2}{1+\lambda_1^2}-\inv{2}}}{\sqrt{2(1-\rho^2)}\left| \frac{\alpha_1}{\sqrt{1+\lambda_1^2}}+\frac{\alpha_2}{\sqrt{1+\lambda_2^2}} \right|\beta_2}\\
& \quad \times \lambda_1^{\inv{1-\rho^2}\left( \inv{\sqrt{1+\lambda_2^2}}-\frac{\rho}{\sqrt{1+\lambda_1^2}} \right)^2-\frac{\lambda_1^2}{1+\lambda_1^2}}\lambda_2^{\left( \frac{\alpha_1}{\sqrt{1+\lambda_1^2}}+\frac{\alpha_2}{\sqrt{1+\lambda_2^2}} \right)^2}\\
& \quad \times \left( \frac{\lambda_1}{\lambda_2} \right)^{\frac{\alpha_1}{\sqrt{1+\lambda_1^2}}\left( \frac{\alpha_1}{\sqrt{1+\lambda_1^2}}+\frac{\alpha_2}{\sqrt{1+\lambda_2^2}} \right)-\inv{(1-\rho^2)\sqrt{1+\lambda_2^2}}\left( \inv{\sqrt{1+\lambda_2^2}}-\frac{\rho}{\sqrt{1+\lambda_1^2}} \right)}\\
&\quad \times u^{\inv{1-\rho^2}\left( \inv{\sqrt{1+\lambda_2^2}}-\frac{\rho}{\sqrt{1+\lambda_1^2}} \right)^2+\left( \frac{\alpha_1}{\sqrt{1+\lambda_1^2}}+\frac{\alpha_2}{\sqrt{1+\lambda_2^2}} \right)^2-\frac{\lambda_1^2}{1+\lambda_1^2}}\\
& \quad (-\log u)^{\inv{1-\rho^2}\left( \inv{\sqrt{1+\lambda_2^2}}-\frac{\rho}{\sqrt{1+\lambda_1^2}} \right)^2+\left( \frac{\alpha_1}{\sqrt{1+\lambda_1^2}}+\frac{\alpha_2}{\sqrt{1+\lambda_2^2}} \right)^2-\frac{\lambda_1^2}{1+\lambda_1^2}-\inv{2}}\\
= & \frac{(2\pi)^{\inv{1-\rho^2}\left( \inv{\sqrt{1+\lambda_1^2}}-\frac{\rho}{\sqrt{1+\lambda_2^2}} \right)^2+\left( \frac{\alpha_1}{\sqrt{1+\lambda_1^2}}+\frac{\alpha_2}{\sqrt{1+\lambda_2^2}} \right)^2-\frac{\lambda_2^2}{1+\lambda_2^2}-\inv{2}}}{\sqrt{2(1-\rho^2)}\left| \frac{\alpha_1}{\sqrt{1+\lambda_1^2}}+\frac{\alpha_2}{\sqrt{1+\lambda_2^2}} \right|\times \beta_1}\\
& \quad \times \lambda_2^{\inv{1-\rho^2}\left( \inv{\sqrt{1+\lambda_1^2}}-\frac{\rho}{\sqrt{1+\lambda_2^2}}\right)^2+\left(  \frac{\alpha_1}{\sqrt{1+\lambda_1^2}}+\frac{\alpha_2}{\sqrt{1+\lambda_2^2}}\right)^2-\frac{\lambda_2^2}{1+\lambda_2^2}+1} \\
& \quad \times \lambda_2^{-\inv{(1-\rho^2)\sqrt{1+\lambda_1^2}}\left( \inv{\sqrt{1+\lambda_1^2}}-\frac{\rho}{\sqrt{1+\lambda_2^2}}
\right)-\frac{\alpha_1}{\sqrt{1+\lambda_1^2}}\left( \frac{\alpha_1}{\sqrt{1+\lambda_1^2}}+\frac{\alpha_2}{\sqrt{1+\lambda_2^2}} \right)}\\
& \quad \times \lambda_1^{\inv{(1-\rho^2)\sqrt{1+\lambda_1^2}}\left( \inv{\sqrt{1+\lambda_1^2}}-\frac{\rho}{\sqrt{1+\lambda_2^2}} \right)+\frac{\alpha_1}{\sqrt{1+\lambda_1^2}}\left( \frac{\alpha_1}{\sqrt{1+\lambda_1^2}}+\frac{\alpha_2}{\sqrt{1+\lambda_2^2}} \right)}\\
& \quad \times u^{\inv{1-\rho^2}\left( \inv{\sqrt{1+\lambda_1^2}}-\frac{\rho}{\sqrt{1+\lambda_2^2}} \right)^2+\left( \frac{\alpha_1}{\sqrt{1+\lambda_1^2}}+\frac{\alpha_2}{\sqrt{1+\lambda_2^2}} \right)^2-\frac{\lambda_2^2}{1+\lambda_2^2}}\\
& \quad \times (-\log u)^{\inv{1-\rho^2}\left( \inv{\sqrt{1+\lambda_1^2}}-\frac{\rho}{\sqrt{1+\lambda_2^2}} \right)^2+\left( \frac{\alpha_1}{\sqrt{1+\lambda_1^2}}+\frac{\alpha_2}{\sqrt{1+\lambda_2^2}} \right)^2-\frac{\lambda_2^2}{1+\lambda_2^2}-\inv{2}}\\
& + \frac{(2\pi)^{\inv{1-\rho^2}\left( \inv{\sqrt{1+\lambda_2^2}}-\frac{\rho}{\sqrt{1+\lambda_1^2}} \right)^2+\left( \frac{\alpha_1}{\sqrt{1+\lambda_1^2}}+\frac{\alpha_2}{\sqrt{1+\lambda_2^2}} \right)^2-\frac{\lambda_1^2}{1+\lambda_1^2}-\inv{2}}}{\sqrt{2(1-\rho^2)}\left| \frac{\alpha_1}{\sqrt{1+\lambda_1^2}}+\frac{\alpha_2}{\sqrt{1+\lambda_2^2}} \right|\beta_2}\\
& \quad \times \lambda_1^{\inv{1-\rho^2}\left( \inv{\sqrt{1+\lambda_2^2}}-\frac{\rho}{\sqrt{1+\lambda_1^2}} \right)^2-\frac{\lambda_1^2}{1+\lambda_1^2}+1 +
\frac{\alpha_1}{\sqrt{1+\lambda_1^2}}\left( \frac{\alpha_1}{\sqrt{1+\lambda_1^2}}+\frac{\alpha_2}{\sqrt{1+\lambda_2^2}} \right)-\inv{(1-\rho^2)\sqrt{1+\lambda_2^2}}\left( \inv{\sqrt{1+\lambda_2^2}}-\frac{\rho}{\sqrt{1+\lambda_1^2}} \right)
} \\
& \quad \times \lambda_2^{\left( \frac{\alpha_1}{\sqrt{1+\lambda_1^2}}+\frac{\alpha_2}{\sqrt{1+\lambda_2^2}} \right)^2
- \frac{\alpha_1}{\sqrt{1+\lambda_1^2}}\left( \frac{\alpha_1}{\sqrt{1+\lambda_1^2}}+\frac{\alpha_2}{\sqrt{1+\lambda_2^2}} \right)+\inv{(1-\rho^2)\sqrt{1+\lambda_2^2}}\left( \inv{\sqrt{1+\lambda_2^2}}-\frac{\rho}{\sqrt{1+\lambda_1^2}} \right)}\\
&\quad \times u^{\inv{1-\rho^2}\left( \inv{\sqrt{1+\lambda_2^2}}-\frac{\rho}{\sqrt{1+\lambda_1^2}} \right)^2+\left( \frac{\alpha_1}{\sqrt{1+\lambda_1^2}}+\frac{\alpha_2}{\sqrt{1+\lambda_2^2}} \right)^2-\frac{\lambda_1^2}{1+\lambda_1^2}}\\
& \quad (-\log u)^{\inv{1-\rho^2}\left( \inv{\sqrt{1+\lambda_2^2}}-\frac{\rho}{\sqrt{1+\lambda_1^2}} \right)^2+\left( \frac{\alpha_1}{\sqrt{1+\lambda_1^2}}+\frac{\alpha_2}{\sqrt{1+\lambda_2^2}} \right)^2-\frac{\lambda_1^2}{1+\lambda_1^2}-\inv{2}}\\
= & \frac{(2\pi)^{\inv{1-\rho^2}\left( \inv{\sqrt{1+\lambda_1^2}}-\frac{\rho}{\sqrt{1+\lambda_2^2}} \right)^2+\left( \frac{\alpha_1}{\sqrt{1+\lambda_1^2}}+\frac{\alpha_2}{\sqrt{1+\lambda_2^2}} \right)^2-\frac{\lambda_2^2}{1+\lambda_2^2}-\inv{2}}}{\sqrt{2(1-\rho^2)}\left| \frac{\alpha_1}{\sqrt{1+\lambda_1^2}}+\frac{\alpha_2}{\sqrt{1+\lambda_2^2}} \right|\times \beta_1}\\
& \quad \times \lambda_2^{\inv{1-\rho^2}\left( \inv{\sqrt{1+\lambda_1^2}}-\frac{\rho}{\sqrt{1+\lambda_2^2}}\right)\left(-\frac{\rho}{\sqrt{1+\lambda_2^2}}\right)+
\frac{\alpha_2}{\sqrt{1+\lambda_2^2}}\left(\frac{\alpha_1}{\sqrt{1+\lambda_1^2}}+\frac{\alpha_2}{\sqrt{1+\lambda_2^2}}\right)-\frac{\lambda_2^2}{1+\lambda_2^2}+1} \\
& \quad \times \lambda_1^{\inv{(1-\rho^2)\sqrt{1+\lambda_1^2}}\left( \inv{\sqrt{1+\lambda_1^2}}-\frac{\rho}{\sqrt{1+\lambda_2^2}} \right)+\frac{\alpha_1}{\sqrt{1+\lambda_1^2}}\left( \frac{\alpha_1}{\sqrt{1+\lambda_1^2}}+\frac{\alpha_2}{\sqrt{1+\lambda_2^2}} \right)}\\
& \quad \times u^{\inv{1-\rho^2}\left( \inv{\sqrt{1+\lambda_1^2}}-\frac{\rho}{\sqrt{1+\lambda_2^2}} \right)^2+\left( \frac{\alpha_1}{\sqrt{1+\lambda_1^2}}+\frac{\alpha_2}{\sqrt{1+\lambda_2^2}} \right)^2-\frac{\lambda_2^2}{1+\lambda_2^2}}\\
& \quad \times (-\log u)^{\inv{1-\rho^2}\left( \inv{\sqrt{1+\lambda_1^2}}-\frac{\rho}{\sqrt{1+\lambda_2^2}} \right)^2+\left( \frac{\alpha_1}{\sqrt{1+\lambda_1^2}}+\frac{\alpha_2}{\sqrt{1+\lambda_2^2}} \right)^2-\frac{\lambda_2^2}{1+\lambda_2^2}-\inv{2}}\\
& + \frac{(2\pi)^{\inv{1-\rho^2}\left( \inv{\sqrt{1+\lambda_2^2}}-\frac{\rho}{\sqrt{1+\lambda_1^2}} \right)^2+\left( \frac{\alpha_1}{\sqrt{1+\lambda_1^2}}+\frac{\alpha_2}{\sqrt{1+\lambda_2^2}} \right)^2-\frac{\lambda_1^2}{1+\lambda_1^2}-\inv{2}}}{\sqrt{2(1-\rho^2)}\left| \frac{\alpha_1}{\sqrt{1+\lambda_1^2}}+\frac{\alpha_2}{\sqrt{1+\lambda_2^2}} \right|\beta_2}\\
& \quad \times \lambda_1^{\inv{1-\rho^2}\left( \inv{\sqrt{1+\lambda_2^2}}-\frac{\rho}{\sqrt{1+\lambda_1^2}} \right)\left(-\frac{\rho}{\sqrt{1+\lambda_1^2}}\right)-\frac{\lambda_1^2}{1+\lambda_1^2}+1 +
\frac{\alpha_1}{\sqrt{1+\lambda_1^2}}\left( \frac{\alpha_1}{\sqrt{1+\lambda_1^2}}+\frac{\alpha_2}{\sqrt{1+\lambda_2^2}}\right)
} \\
& \quad \times \lambda_2^{
\frac{\alpha_2}{\sqrt{1+\lambda_2^2}}\left( \frac{\alpha_1}{\sqrt{1+\lambda_1^2}}+\frac{\alpha_2}{\sqrt{1+\lambda_2^2}} \right)+\inv{(1-\rho^2)\sqrt{1+\lambda_2^2}}\left( \inv{\sqrt{1+\lambda_2^2}}-\frac{\rho}{\sqrt{1+\lambda_1^2}} \right)}\\
&\quad \times u^{\inv{1-\rho^2}\left( \inv{\sqrt{1+\lambda_2^2}}-\frac{\rho}{\sqrt{1+\lambda_1^2}} \right)^2+\left( \frac{\alpha_1}{\sqrt{1+\lambda_1^2}}+\frac{\alpha_2}{\sqrt{1+\lambda_2^2}} \right)^2-\frac{\lambda_1^2}{1+\lambda_1^2}}\\
& \quad (-\log u)^{\inv{1-\rho^2}\left( \inv{\sqrt{1+\lambda_2^2}}-\frac{\rho}{\sqrt{1+\lambda_1^2}} \right)^2+\left( \frac{\alpha_1}{\sqrt{1+\lambda_1^2}}+\frac{\alpha_2}{\sqrt{1+\lambda_2^2}} \right)^2-\frac{\lambda_1^2}{1+\lambda_1^2}-\inv{2}}\\
= & \frac{(2\pi)^{\inv{1-\rho^2}\left( \inv{\sqrt{1+\lambda_2^2}}-\frac{\rho}{\sqrt{1+\lambda_1^2}} \right)^2+\left( \frac{\alpha_1}{\sqrt{1+\lambda_1^2}}+\frac{\alpha_2}{\sqrt{1+\lambda_2^2}} \right)^2-\frac{\lambda_1^2}{1+\lambda_1^2}-\inv{2}}}{\sqrt{2(1-\rho^2)}\left| \frac{\alpha_1}{\sqrt{1+\lambda_1^2}}+\frac{\alpha_2}{\sqrt{1+\lambda_2^2}} \right|\times \beta_1}\\
& \quad \times \lambda_1^{\inv{(1-\rho^2)\sqrt{1+\lambda_1^2}}\left( \inv{\sqrt{1+\lambda_1^2}}-\frac{\rho}{\sqrt{1+\lambda_2^2}} \right)+\frac{\alpha_1}{\sqrt{1+\lambda_1^2}}\left( \frac{\alpha_1}{\sqrt{1+\lambda_1^2}}+\frac{\alpha_2}{\sqrt{1+\lambda_2^2}} \right)}\\
& \quad \times \lambda_2^{
\frac{\alpha_2}{\sqrt{1+\lambda_2^2}}\left( \frac{\alpha_1}{\sqrt{1+\lambda_1^2}}+\frac{\alpha_2}{\sqrt{1+\lambda_2^2}} \right)+\inv{(1-\rho^2)\sqrt{1+\lambda_2^2}}\left( \inv{\sqrt{1+\lambda_2^2}}-\frac{\rho}{\sqrt{1+\lambda_1^2}} \right)}\\
& \quad \times u^{\inv{1-\rho^2}\left( \inv{\sqrt{1+\lambda_2^2}}-\frac{\rho}{\sqrt{1+\lambda_1^2}} \right)^2+\left( \frac{\alpha_1}{\sqrt{1+\lambda_1^2}}+\frac{\alpha_2}{\sqrt{1+\lambda_2^2}} \right)^2-\frac{\lambda_1^2}{1+\lambda_1^2}}\\
& \quad \times (-\log u)^{\inv{1-\rho^2}\left( \inv{\sqrt{1+\lambda_2^2}}-\frac{\rho}{\sqrt{1+\lambda_1^2}} \right)^2+\left( \frac{\alpha_1}{\sqrt{1+\lambda_1^2}}+\frac{\alpha_2}{\sqrt{1+\lambda_2^2}} \right)^2-\frac{\lambda_1^2}{1+\lambda_1^2}-\inv{2}}\\
& + \frac{(2\pi)^{\inv{1-\rho^2}\left( \inv{\sqrt{1+\lambda_2^2}}-\frac{\rho}{\sqrt{1+\lambda_1^2}} \right)^2+\left( \frac{\alpha_1}{\sqrt{1+\lambda_1^2}}+\frac{\alpha_2}{\sqrt{1+\lambda_2^2}} \right)^2-\frac{\lambda_1^2}{1+\lambda_1^2}-\inv{2}}}{\sqrt{2(1-\rho^2)}\left| \frac{\alpha_1}{\sqrt{1+\lambda_1^2}}+\frac{\alpha_2}{\sqrt{1+\lambda_2^2}} \right|\beta_2}\\
& \quad \times \lambda_1^{\inv{(1-\rho^2)\sqrt{1+\lambda_1^2}}\left( \inv{\sqrt{1+\lambda_1^2}}-\frac{\rho}{\sqrt{1+\lambda_2^2}} \right)+\frac{\alpha_1}{\sqrt{1+\lambda_1^2}}\left( \frac{\alpha_1}{\sqrt{1+\lambda_1^2}}+\frac{\alpha_2}{\sqrt{1+\lambda_2^2}} \right)}\\
& \quad \times \lambda_2^{
\frac{\alpha_2}{\sqrt{1+\lambda_2^2}}\left( \frac{\alpha_1}{\sqrt{1+\lambda_1^2}}+\frac{\alpha_2}{\sqrt{1+\lambda_2^2}} \right)+\inv{(1-\rho^2)\sqrt{1+\lambda_2^2}}\left( \inv{\sqrt{1+\lambda_2^2}}-\frac{\rho}{\sqrt{1+\lambda_1^2}} \right)}\\
&\quad \times u^{\inv{1-\rho^2}\left( \inv{\sqrt{1+\lambda_2^2}}-\frac{\rho}{\sqrt{1+\lambda_1^2}} \right)^2+\left( \frac{\alpha_1}{\sqrt{1+\lambda_1^2}}+\frac{\alpha_2}{\sqrt{1+\lambda_2^2}} \right)^2-\frac{\lambda_1^2}{1+\lambda_1^2}}\\
& \quad (-\log u)^{\inv{1-\rho^2}\left( \inv{\sqrt{1+\lambda_2^2}}-\frac{\rho}{\sqrt{1+\lambda_1^2}} \right)^2+\left(\frac{\alpha_1}{\sqrt{1+\lambda_1^2}}+\frac{\alpha_2}{\sqrt{1+\lambda_2^2}} \right)^2-\frac{\lambda_1^2}{1+\lambda_1^2}-\inv{2}},
\end{align*}
}
as 
{\adb
\begin{align*}
& \inv{1-\rho^2}\left(\inv{\sqrt{1+\lambda_1^2}}-\frac{\rho}{\sqrt{1+\lambda_2^2}} \right)^2-\frac{\lambda_2^2}{1+\lambda_2^2}\\
=& \inv{1-\rho^2}\left( \inv{1+\lambda_1^2}-\frac{2\rho}{\sqrt{(1+\lambda_1^2)(1+\lambda_2^2)}}+\frac{\rho^2}{1+\lambda_2^2} \right) - \frac{\lambda_2^2}{1+\lambda_2^2}\\
=& \left[\inv{1-\rho^2}\left( \inv{1+\lambda_2^2} - \frac{2\rho}{\sqrt{(1+\lambda_1^2)(1+\lambda_2^2)}} + \frac{\rho^2}{1+\lambda_1^2}\right)- \frac{\lambda_1^2}{1+\lambda_1^2}\right] \\
& \quad + \frac{(1-\rho^2)\lambda_1^2 +1- \rho^2}{(1-\rho^2)(1+\lambda_1^2)}- \frac{(1-\rho^2)\lambda_2^2+(1-\rho^2)}{(1-\rho^2)(1+\lambda_2^2)}\\
=& \inv{1-\rho^2}\left( \inv{\sqrt{1+\lambda_2^2}}-\frac{\rho}{\sqrt{1+\lambda_1^2}} \right)^2-\frac{\lambda_1^2}{1+\lambda_1^2};
\end{align*}
}
and 
{\adb
\begin{align*}
& -\frac{\rho}{(1-\rho^2)\sqrt{1+\lambda_2}} \left(\inv{\sqrt{1+\lambda_1^2}}-\frac{\rho}{\sqrt{1+\lambda_2^2}} \right)-\frac{\lambda_2^2}{1+\lambda_2^2}+1 \\
= & -\frac{\rho}{(1-\rho^2)\sqrt{(1+\lambda_2^2)(1+\lambda_1^2)}} + \frac{\rho^2}{(1+\lambda_2^2)(1-\rho^2)} + \frac{-\lambda_2^2(1-\rho^2)+ (1+\lambda_2^2)-\rho^2(1+\lambda_2^2)}{(1+\lambda_2^2)(1-\rho^2)} \\
= & \inv{(1-\rho^2)\sqrt{1+\lambda_2^2}}\left(\inv{\sqrt{1+\lambda_2^2}}- \frac{\rho}{\sqrt{1+\lambda_1^2}}\right);\\
\Rightarrow \quad & -\frac{\rho}{(1-\rho^2)\sqrt{1+\lambda_1}} \left(\inv{\sqrt{1+\lambda_2^2}}-\frac{\rho}{\sqrt{1+\lambda_1^2}} \right)-\frac{\lambda_1^2}{1+\lambda_1^2}+1 \\
=& \inv{(1-\rho^2)\sqrt{1+\lambda_1^2}}\left(\inv{\sqrt{1+\lambda_1^2}}- \frac{\rho}{\sqrt{1+\lambda_2^2}}\right).
\end{align*}
}
This means the two summands are asymptotically the same, apart from differing factors $ 1/\beta_1, 1/\beta_2.$ As a result, 
{\adb
\begin{align*}
& \frac{d C(u,u)}{d u}\\
=& \left[ \inv{\beta_1}+\inv{\beta_2}\right] \times \frac{(2\pi)^{\inv{1-\rho^2}\left( \inv{\sqrt{1+\lambda_2^2}}-\frac{\rho}{\sqrt{1+\lambda_1^2}} \right)^2+\left( \frac{\alpha_1}{\sqrt{1+\lambda_1^2}}+\frac{\alpha_2}{\sqrt{1+\lambda_2^2}} \right)^2-\frac{\lambda_1^2}{1+\lambda_1^2}-\inv{2}}}{\sqrt{2(1-\rho^2)}\left| \frac{\alpha_1}{\sqrt{1+\lambda_1^2}}+\frac{\alpha_2}{\sqrt{1+\lambda_2^2}} \right|}\\
& \quad \times 
\lambda_1^{\inv{(1-\rho^2)\sqrt{1+\lambda_1^2}}\left( \inv{\sqrt{1+\lambda_1^2}}-\frac{\rho}{\sqrt{1+\lambda_2^2}} \right)+\frac{\alpha_1}{\sqrt{1+\lambda_1^2}}\left( \frac{\alpha_1}{\sqrt{1+\lambda_1^2}}+\frac{\alpha_2}{\sqrt{1+\lambda_2^2}}\right)}\\
& \quad \times 
\lambda_2^{\inv{(1-\rho^2)\sqrt{1+\lambda_2^2}}\left( \inv{\sqrt{1+\lambda_2^2}}-\frac{\rho}{\sqrt{1+\lambda_1^2}} \right)+\frac{\alpha_2}{\sqrt{1+\lambda_2^2}}\left( \frac{\alpha_1}{\sqrt{1+\lambda_1^2}}+\frac{\alpha_2}{\sqrt{1+\lambda_2^2}}\right)}\\
&\quad \times u^{\inv{1-\rho^2}\left( \inv{\sqrt{1+\lambda_2^2}}-\frac{\rho}{\sqrt{1+\lambda_1^2}} \right)^2+\left( \frac{\alpha_1}{\sqrt{1+\lambda_1^2}}+\frac{\alpha_2}{\sqrt{1+\lambda_2^2}} \right)^2-\frac{\lambda_1^2}{1+\lambda_1^2}}\\
& \quad (-\log u)^{\inv{1-\rho^2}\left( \inv{\sqrt{1+\lambda_2^2}}-\frac{\rho}{\sqrt{1+\lambda_1^2}} \right)^2+\left( \frac{\alpha_1}{\sqrt{1+\lambda_1^2}}+\frac{\alpha_2}{\sqrt{1+\lambda_2^2}} \right)^2-\frac{\lambda_1^2}{1+\lambda_1^2}-\inv{2}}
\end{align*}
} 
\end{proof} Finally,  we reconcile the the above general result when  $\lambda_1 >0, \lambda_2 >0  $, with the case expressed in Theorem \ref{thm:fungTailAsymptoticsBivariate2016} 

When $\alpha_1 = \alpha_2 = \alpha >0$, s0 that $\lambda_1= \lambda_2 = \lambda$
{\allowdisplaybreaks
\begin{align*}
& \inv{1-\rho^2}\left(\inv{\sqrt{1+\lambda_2^2}}-\frac{\rho}{\sqrt{1+\lambda_1^2}} \right)^2+\left( \frac{\alpha_1}{\sqrt{1+\lambda_1^2}}+\frac{\alpha_2}{\sqrt{1+\lambda_2^2}} \right)^2-\frac{\lambda_1^2}{1+\lambda_1^2}\\
= & \inv{1-\rho^2}\left(\frac{1-\rho}{\sqrt{1+\lambda^2}}\right)^2 + \left( \frac{2\alpha}{\sqrt{1+\lambda^2}}\right)^2-\frac{\lambda^2}{1+\lambda^2}\\
= & \frac{1-\rho}{(1 +\rho)(1+\lambda^2)} + \frac{4\alpha^2}{1+\lambda^2} - \frac{\lambda^2}{1+\lambda^2}\\
=& \frac{1+\alpha^2(1-\rho^2)}{1+2\alpha^2(1+\rho)}\left[ \frac{1-\rho}{1+\rho} + 4\alpha^2 - \frac{\alpha^2(1+\rho)^2}{1+\alpha^2(1-\rho^2)}\right], \quad \text{as $\lambda = \frac{\alpha(1+\rho)}{\sqrt{1+\alpha ^2(1-\rho^2)}}$;}\\
=& \frac{1}{1+2\alpha^2(1+\rho)}\left[\frac{(1-\rho)(1+\alpha^2(1-\rho^2)) + 4\alpha^2(1+\rho)(1+\alpha^2(1-\rho^2)) - \alpha^2(1+\rho)^3}{1+\rho}\right]\\
=& \frac{1}{1+2\alpha^2(1+\rho)}\left\{ 
\frac{(1-\rho) + (1+\rho)\left[\alpha^2(1-2\rho+\rho^2)+4\alpha^2-\alpha^2(1+2\rho+\rho^2)\right] + 4\alpha^4(1+\rho)^2(1-\rho)}{1+\rho}\right\}\\
=& \frac{1}{1+2\alpha^2(1+\rho)}\left\{ 
\frac{(1-\rho) + (1+\rho)\left[4\alpha^2(1-\rho)\right]+ 4\alpha^4(1+\rho)^2(1-\rho)}{1+\rho}\right\}\\
=& \frac{(1-\rho)\left[1+2\alpha^2(1+\rho)\right]^2}{(1+2\alpha^2(1+\rho))(1+\rho)}\\
=& \frac{(1-\rho)(1+2\alpha^2(1+\rho))}{1+\rho} = \beta^2, \quad \text{where $\beta$ is as defined in Theorem \ref{thm:fungTailAsymptoticsBivariate2016}}.
\end{align*}
}
Further, 
{\allowdisplaybreaks
\begin{align*}
    & \inv{(1-\rho^2)\sqrt{1+\lambda_1^2}}\left( \inv{\sqrt{1+\lambda_1^2}}-\frac{\rho}{\sqrt{1+\lambda_2^2}} \right)+\frac{\alpha_1}{\sqrt{1+\lambda_1^2}}\left( \frac{\alpha_1}{\sqrt{1+\lambda_1^2}}+\frac{\alpha_2}{\sqrt{1+\lambda_2^2}}\right)\\
    =& \inv{(1+\rho)(1+\lambda^2)}+ \frac{2\alpha^2}{1+\lambda^2}\\
    =& \frac{1+\alpha^2(1-\rho^2)}{1+2\alpha^2(1+\rho)}\left[
    \frac{1 + 2\alpha^2(1+\rho)}{1+\rho}
    \right] \\
    =& \frac{2+2\alpha^2(1-\rho^2)}{2(1+\rho)} \\
    =& \frac{1-\rho + 2\alpha^2(1-\rho)(1+\rho) + 1+\rho}{2(1+\rho)}\\
    =& \frac{\beta^2}{2} + \inv{2}.
\end{align*}

}
This means that 
{\adb
\begin{align*}
& \frac{d C(u,u)}{d u}\\
=& \left[ \inv{\beta_1}+\inv{\beta_2}\right] \times \frac{(2\pi)^{\inv{1-\rho^2}\left( \inv{\sqrt{1+\lambda_2^2}}-\frac{\rho}{\sqrt{1+\lambda_1^2}} \right)^2+\left( \frac{\alpha_1}{\sqrt{1+\lambda_1^2}}+\frac{\alpha_2}{\sqrt{1+\lambda_2^2}} \right)^2-\frac{\lambda_1^2}{1+\lambda_1^2}-\inv{2}}}{\sqrt{2(1-\rho^2)}\left| \frac{\alpha_1}{\sqrt{1+\lambda_1^2}}+\frac{\alpha_2}{\sqrt{1+\lambda_2^2}} \right|}\\
& \quad \times 
\lambda_1^{\inv{(1-\rho^2)\sqrt{1+\lambda_1^2}}\left( \inv{\sqrt{1+\lambda_1^2}}-\frac{\rho}{\sqrt{1+\lambda_2^2}} \right)+\frac{\alpha_1}{\sqrt{1+\lambda_1^2}}\left( \frac{\alpha_1}{\sqrt{1+\lambda_1^2}}+\frac{\alpha_2}{\sqrt{1+\lambda_2^2}}\right)}\\
& \quad \times 
\lambda_2^{\inv{(1-\rho^2)\sqrt{1+\lambda_2^2}}\left( \inv{\sqrt{1+\lambda_2^2}}-\frac{\rho}{\sqrt{1+\lambda_1^2}} \right)+\frac{\alpha_2}{\sqrt{1+\lambda_2^2}}\left( \frac{\alpha_1}{\sqrt{1+\lambda_1^2}}+\frac{\alpha_2}{\sqrt{1+\lambda_2^2}}\right)}\\
&\quad \times u^{\inv{1-\rho^2}\left( \inv{\sqrt{1+\lambda_2^2}}-\frac{\rho}{\sqrt{1+\lambda_1^2}} \right)^2+\left( \frac{\alpha_1}{\sqrt{1+\lambda_1^2}}+\frac{\alpha_2}{\sqrt{1+\lambda_2^2}} \right)^2-\frac{\lambda_1^2}{1+\lambda_1^2}}\\
& \quad (-\log u)^{\inv{1-\rho^2}\left( \inv{\sqrt{1+\lambda_2^2}}-\frac{\rho}{\sqrt{1+\lambda_1^2}} \right)^2+\left( \frac{\alpha_1}{\sqrt{1+\lambda_1^2}}+\frac{\alpha_2}{\sqrt{1+\lambda_2^2}} \right)^2-\frac{\lambda_1^2}{1+\lambda_1^2}-\inv{2}}\\
=&  \left[ \inv{\beta_1}+\inv{\beta_2}\right]\frac{(2\pi)^{\beta^2-\inv{2}}}{\sqrt{2(1-\rho^2)}
\left| \frac{2\alpha}{\sqrt{1+\lambda^2}}\right| } \lambda^{\beta^2+1} u^{\beta^2}(-\log u)^{\beta^2-\inv{2}}\\
=& \left\{\left[ \inv{\beta_1}+\inv{\beta_2}\right]\frac{\lambda}{2\sqrt{(1-\rho^2)}\left|\frac{2\alpha}{\sqrt{1+\lambda^2}}\right|}\right\}
\frac{(2\pi\lambda)^{\beta^2}}{\sqrt{\pi}}u^{\beta^2}(-\log u)^{\beta^2-\inv{2}}\\
\end{align*}
} 
and the results will collapse back the one in Theorem \ref{thm:fungTailAsymptoticsBivariate2016} if 
{\adb
\begin{align*}
    \left[ \inv{\beta_1}+\inv{\beta_2}\right]\frac{\lambda}{2\sqrt{(1-\rho^2)}\left|\frac{2\alpha}{
    \sqrt{1+\lambda^2}}\right|} = \inv{\beta(1+\beta^2)}.
\end{align*}
}
Since $\alpha_1 = \alpha_2 = \alpha$, we have 
\begin{align*}
    \beta_1 = \beta_2 = \frac{1-\rho}{1-\rho^2} + \alpha(\alpha+\alpha) = \inv{1+\rho}+ 2\alpha^2 = \frac{1+2\alpha^2(1+\rho)}{1+\rho}, 
\end{align*}
and 
{\adb
\begin{align*}
    & \left[ \inv{\beta_1}+\inv{\beta_2}\right]\frac{\lambda}{2\sqrt{(1-\rho^2)}\left|\frac{2\alpha}{
    \sqrt{1+\lambda^2}}\right|} \\
    =& \frac{\lambda\sqrt{1+\lambda^2}(1+\rho)}{(1+2\alpha^2(1+\rho))\sqrt{1-\rho^2}(2\alpha)}\\
    =& \frac{\lambda\sqrt{1+\lambda^2}}{2\alpha\sqrt{1+2\alpha^2(1+\rho)}}\times \left(\frac{\sqrt{1+\rho}}{\sqrt{(1-\rho)(1+2\alpha^2(1+\rho)}}\right) \\
    =& \frac{\frac{\alpha(1+\rho)}{\sqrt{1+\alpha^2(1-\rho^2)}}\times \frac{\sqrt{1+2\alpha^2(1+\rho)}}{\sqrt{1+\alpha^2(1-\rho^2)}}}{2\alpha \sqrt{1+2\alpha^2(1+\rho)}}\times \inv{\beta}\\
    =& \frac{1+\rho}{\beta(2+2\alpha^2(1-\rho^2))} \\
    =& \frac{1}{\beta\left[\frac{(1-\rho)(1+2\alpha^2(1+\rho))}{1+\rho}+1\right]} \\
    =& \inv{\beta(1+\beta^2)}.
\end{align*}
}

\bibliography{library}

\end{document}